\documentclass[reqno]{amsart}
\usepackage{graphicx}
\usepackage{amssymb}
\usepackage{subcaption}
\usepackage{float}
\usepackage[]{algorithm2e}

\usepackage{algorithmic}
\usepackage{hyperref}
\hypersetup{
colorlinks=true,
linkcolor=red,
filecolor=magenta, 
urlcolor=cyan,
} 
\usepackage{color}
\usepackage{bm}

\newtheorem{theorem}{Theorem}[section]

\theoremstyle{definition}

\newtheorem{example}[theorem]{Example}

\theoremstyle{remark}
\newtheorem{remark}[theorem]{Remark}

\DeclareMathOperator*{\argmin}{arg\,min}
\numberwithin{equation}{section}

\newcommand{\lprod}[2]{\Big(#1,#2\Big)_2}

\newcommand{\mbb}[1]{\mathbb{#1}}
\newcommand{\mcal}[1]{\mathcal{#1}}

\title{Interpretation of Plug-and-Play (PnP) algorithms from a different angle}
\author{Abinash Nayak}
\usepackage{fancyhdr}
\begin{document}


\address{Visiting Assistant Professor, Department of Mathematics, University of Alabama at Birmingham, University Hall, Room 4005, 1402 10th Avenue South, Birmingham AL 35294-1241, (p) 205.934.2154, (f) 205.934.9025}
\email{nash101@uab.edu; avinashnike01@gmail.com}


\date{\today}


\begin{abstract}
It's well-known that inverse problems are ill-posed and to solve them meaningfully, one has to employ regularization methods. Traditionally, the most popular regularization approaches are Variational-type approaches, i.e., penalized/constrained functional minimization. In recent years, the classical regularization approaches have been replaced by the so called plug-and-play (PnP) algorithms, which copies the proximal gradient minimization processes, such as ADMM or FISTA, but with any general denoiser. However, unlike the traditional proximal gradient methods, the theoretical underpinnings and convergence results have been insufficient for these PnP-algorithms. Hence, the results from these algorithms, though empirically outstanding, are not well-defined, in the sense of, being a minimizer of a Variational problem, or in some other forms. In this paper, we address this question of ``well-definedness", but from a different angle. We explain these algorithms from the viewpoint of a semi-iterative regularization method. In addition, we expand the family of regularized solutions, corresponding to the classical semi-iterative methods, to a much larger class, which encompasses these algorithms, as well as, enhance the recovery process. We conclude with several numerical results which validate the developed theories and reflect the improvements over the traditional PnP-algorithms, such as ADMM-PnP and FISTA-PnP.  
\end{abstract}
\maketitle

\section{\textbf{Introduction}}
\subsection{Inverse Problems and Regularization:}
Mathematically, an inverse problem is often expressed as the problem of estimating a (\textit{source}) $\hat{x}$ which satisfies, for a given (\textit{effect}) $b$, the following matrix (or operator) equation
\begin{equation}\label{Inv. prob.}
A\hat{x} = b,
\end{equation}
where the matrix $A \in \mbb{R}^{n\times m}$ and the vectors $\hat{x} \in \mbb{R}^{n}$, $b \in \mbb{R}^{m}$ are the discrete approximations of an infinite dimensional model describing the underlying physical phenomenon. The inverse problem \eqref{Inv. prob.} is usually ill-posed, in the sense of violating any one of the Hadamard's conditions for well-posedness: (i) Existence of a solution (for $b$ not in the range of $A$), (ii) uniqueness of the solution (for non-trivial null-space of $A$) and, (iii) continuous dependence on the data (for ill-conditioned A). Conditions (i) and (ii) can be circumvented by relaxing the definition of a solution for \eqref{Inv. prob.}, for example, finding the least square solution or the minimal norm solution (i.e., the pseudo-inverse solution $x^\dagger$). The most (practically) significant condition is condition-(iii), since failing of this leads to an absurd (unstable) solution. That is, for an (injective) A and an exact $b$ (noiseless), the solution of \eqref{Inv. prob.} can be approximated by the (LS) least-square solution ($x^\dagger$), i.e., $x^\dagger$ is the minimizer of the following least-square functional
\begin{equation}\label{LS functional}
    F(x) = ||Ax - b||_2^2,
\end{equation}
where as, for a noisy data $b_\delta$ (which is practically true all the time) such that $||b - b_\delta|| \leq \delta$, the simple least-square solution ($x^\dagger_\delta$) with respect to $b_\delta$ in \eqref{LS functional} fails to approximate the true solution, i.e., $||x^\dagger_\delta - x^\dagger|| >> \delta$, which in turn implies, $||x^\dagger_\delta - \hat{x}|| >> \delta$, due to the ill-posedness of the inverse problem \eqref{Inv. prob.}. To counter such instabilities or ill-posedness of inverse problems, regularization methods have to be employed, which are broadly divided into two types. 
\subsection{Variational (or penalized) regularization}
\mbox{ }

Such approaches, also known as Tikhonov-type regularization, are probably the most well known regularization techniques for solving linear as well as nonlinear inverse problems (see \cite{Engl+Hanke+Neubauer,Bakushinsky+Goncharsky, Groetsch, Baumeister,Morozov_a}), where, instead of minimizing the simple least-square functional \eqref{LS functional}, one minimizes a penalized (or constrained) functional:
\begin{equation}\label{Gen. Tik. fun.}
F(x;\mcal{D},\lambda,\mcal{R}) =  \mcal{D}(Ax,b_\delta) + \lambda \mcal{R}(x),
\end{equation}
where $\mcal{D}$ is called the data-fidelity term (imposing data-consistency), $\mcal{R}$ is the regularization term (imposing certain structures, based on some prior knowledge of the solution $\hat{x}$) and, $\lambda \geq 0$ is the regularization parameter that balances the trade-off between them, depending on the noise level $\delta$, i.e., $\lambda=\lambda(\delta)$. The formulation \eqref{Gen. Tik. fun.} also has a Bayesian interpretation, where the minimization of $F(x;\mcal{D},\lambda,\mcal{R})$ corresponds to the maximum-a-posteriori (MAP) estimate of $\hat{x}$ given $b_\delta$, where the likelihood of $b_\delta$ is proportional to $exp(-\mcal{D}(x))$ and the prior distribution on $\hat{x}$ is proportional to $exp(-\mcal{R}(x))$. Classically, $\mcal{D}(Ax,b_\delta) = ||{Ax - b_\delta}||_p^p$ and $\mcal{R}(x) = ||Lx - \bar{x}||_q^q$, where $L$ is a regularization matrix with the null spaces of $A$ and $L$ intersecting trivially, and $p$, $q$ determine the involved norms. For large scale problems, the minimization of \eqref{Gen. Tik. fun.} is done iteratively, i.e., starting from an initial guess $x_0^\delta$ and step-sizes $\tau_k > 0$, the minimizer of \eqref{Gen. Tik. fun.} is approximated via the sequence 
\begin{align}\label{Tik. Iter.}
    x_{k+1}^\delta = x_k^\delta - \tau_k \nabla_x F(x_k^\delta; \mcal{D},\lambda, \mcal{R}),
\end{align}
Hence, one has to address the following two technicalities,
\begin{enumerate}
    \item Compute, if possible, the gradient $\nabla_x F(x_k^\delta; \mcal{D},\lambda, \mcal{R})$ at each k. It's easier when both $\mcal{D}(x)$ and $\mcal{R}(x)$ are differentiable, as then, $\nabla_x F(x_k^\delta; \mcal{D}(x_k^\delta),\lambda, \mcal{R}) = \nabla_x \mcal{D}(x_k^\delta) + \lambda \nabla_x \mcal{R}(x_k^\delta)$. However, when they are non-differentiable, then it's not that straight forward. In such scenarios, for convex $\mcal{R}$, one has to make use of convex optimization techniques, such as proximal gradient methods, to circumvent the differentiability issue, which is discussed in later sections.  
    \item Convergence of the sequence of iterates $\{x_k^\delta\}$ as $k \rightarrow \infty$. For strictly convex $\mcal{D}$ and $\mcal{R}$, we have $x_k^\delta \xrightarrow{k \rightarrow \infty} x^\delta_*$, where $x^\delta_*$ is the global minimizer of \eqref{Gen. Tik. fun.}. Again, if they are non-convex (which is the case for most modern $\mcal{R}$), then one has to analyze the behaviour of the sequence when $k\rightarrow \infty$, for example, convergence to a local minimizer or a saddle-point etc. 
\end{enumerate}

\subsection{Related Works}
\mbox{ }

The above two queries have been studied extensively for classical regularizers, such as TV or sparsifying wavelet transformation, and the associated results can be found in the literature, for example in \cite{Osher_Burger_Goldfarb_Xu_Yin, Ma_Yin_Zhang_Chakraborty,Yang_Zhang_Yin, Ravishankar_Bresler, Liao_Sapiro}. Typically, for a non-smooth $\mcal{R}$, which is proper, closed, and convex, the differentiability issue is circumvented by using a proximal operator (see \cite{Beck_Teboulle, Boyd_Parikh_Chu_Peleato_Eckstein, Chambolle_Pock}, and references therein), given by
\begin{equation}
    \mbox{Prox}_{\lambda \mcal{R}}(v) = \argmin_x \; \lambda \mcal{R}(x) + \frac{1}{2} ||x - v||_2^2.
\end{equation}
Basically, for smooth $\mcal{D}(x)$ and non-smooth $\mcal{R}(x)$, the minimization problem corresponding to \eqref{Gen. Tik. fun.} can be solved via two first-order iterative methods:
\begin{enumerate}
    \item Forward-backward splitting (FBS), also known as Iterative shrinkage/soft thresholding algorithm (ISTA) and has a faster variant Fast ISTA (FISTA), where each minimization step is divided into two sub-steps, given by
    \begin{align}\label{FBS min.}
        z_{k+1}^\delta &= x_k^\delta - \tau_k \nabla_x \mcal{D}(x_k^\delta) \;\; \longleftarrow \text{ data-consistency step}\\
        x_{k+1}^\delta &= \text{Prox}_{\lambda \tau_k \mcal{R}}(z_{k+1}^\delta) \;\; \longleftarrow \text{ data-denoising step}
    \end{align}    
    where $\tau_k \geq 0$ is the step-size at the k$^{th}$ iteration.

    \item Alternating direction method of multipliers (ADMM), where three sequences are alternatively updated as follows, 
    \begin{align}\label{ADMM min.}
        x_{k+1}^\delta &= \text{Prox}_{\frac{1}{\rho}\mcal{D}}(z_k^\delta - u_k^\delta) \;\; \longleftarrow \text{ data-consistency step}\\
        z_{k+1}^\delta &= \text{Prox}_{\frac{\lambda}{\rho}\mcal{R}}(x_{k+1}^\delta + u_{k}^\delta) \;\; \longleftarrow \text{ data-denoising step}\\
        u_{k+1}^\delta &= u_k^\delta + x_{k+1}^\delta - z_{k+1}^\delta, \;\; \longleftarrow \text{ noise update step}
    \end{align}
    where $\rho > 0$ is the Lagrangian parameter, which only effects the speed of convergence and not the solution (minimizer) of \eqref{Gen. Tik. fun.}.
\end{enumerate}
From the above two expressions, one can observe that, each method comprises of two fundamental steps: (1) data-consistency and (2) data-denoising. This motivated, authors in \cite{Venkatakrishnan_Bouman_Wohlberg}, to replace the Prox$_{\sigma \mcal{R}}$ operator in the denoising step of ADMM by an off-the-shelf denoiser $H_\sigma$, with denoising strength corresponding to noise level $\sigma$, and termed it as the PnP-algorithm (plug-and-play method). However, note that, once the proximal operator is replaced by any general denoiser, then the Variational problem \eqref{Gen. Tik. fun.} breaks down, as not all denoisers can be expressed as a proximal operator of some function $\mcal{R}$. Hence, certain natural questions arise, such as, where does this new sequence of iterates converge (if they do)? And/or, how to assign a meaning to the recovered solution?

Though empirical results show the convergence of these PnP-algorithms, there is no proof of it, for any general denoisers. However, under certain assumptions and restrictions (such as boundedness, nonexpansiveness, etc.) on the denoiser, there have been some convergence proof, see \cite{Chan_Wang_Elgendy, Buzzard_Chan_Sreehari_Bouman, Ryu_Liu_Wang_Chen_Zhangyang_Yin, Teodoro_Bioucas_Figueiredo, Sun_Wohlberg_Kamilov} and references therein. There are also some other variants of such PnP-methods, such as Regularization by Denoising (RED)\cite{Romano_Elad_Milanfar}, Regularization by Artifact-Removal (RARE)  \cite{Kamilov_Liu_Sun_Eldeniz_Gan_An}, etc.

\subsection*{Contribution of this paper}
\mbox{}

\begin{itemize}
\item In this paper, we try to answer the above questions, the convergence of the iterates and the meaning of the solution of a PnP-method, from a different angle. We show a connection between these PnP-algorithms and semi-iterative regularization methods. We then present a bound for the difference of a ``PnP-solution" and a ``regularized LS-solution", where the definitions of these terms are explained below.

\item PnP-algorithms, like Variational regularization, are very sensitive to the regularization/denoising parameter. Here, we present a method to significantly reduce the sensitivity of these algorithms to the denoiser's strength of denoising. This is crucial, since this can transform a bad denoiser (denoisers with too strong or too weak denoising) to a good one (an appropriate level of denoising), without actually altering the involved denoiser.

\item We also address the question of when to stop/terminate the iterations to recover ``an appropriate solution". Again, the meaning of an appropriate solution is described below, which depends on some selection criteria.

\item In addition, we compare the FBS-PnP algorithm with the ADMM-PnP algorithm and point out some of the pros and cons of these algorithms. Note that, although in the traditional setting \eqref{Gen. Tik. fun.} both these algorithms produce a similar result, the minimizer of \eqref{Gen. Tik. fun.}, but for their PnP versions they may differ depending on the scenarios, as the convergence is not known.

\item Finally, we conclude with presenting several numerical results in support of the theories developed in this paper, which also validates the efficiency and effectiveness of the modified algorithms suggested here. We perform computational experiments on the X-ray computed tomography problems, and compare the results obtained using our approach with the traditional ADMM-PnP and FBS-PnP algorithms, with BM3D as the denoiser.
\end{itemize}

\section{\textbf{Structure imposing descent directions}}
In the previous section we explained Variational or Tikhonov-type regularization methods and their variants. The second most popular regularization technique is known as \textit{Semi-iterative regularization method}.
\subsection{Semi-iterative regularization}
\mbox{ }

Landweber, in \cite{Landweber}, showed that 
minimizing the simple LS functional \eqref{LS functional} via the gradient direction, with a constant step-size, leads to the recovery of the solution of \eqref{Inv. prob.} in a semi-convergent manner, i.e., starting from an initial guess $x_0^\delta$, the sequence of iterates
\begin{equation}\label{LS descent}
x_{k}^\delta = x_{k-1}^\delta + d_{k-1}^\delta,
\end{equation}
where $d_{k-1}^\delta := -\tau \nabla_x \mcal{D}(x_{k-1}^\delta) = -\tau A^*(Ax_{k-1}^\delta - b_\delta)$, when $\mcal{D}(x) = ||Ax - b_\delta||_2^2$, for $0 \leq \tau < \frac{1}{2||A||_2^2} $, initially converges towards the true LS-solution $x^\dagger$, then diverges away from $x^\dagger$, i.e., the recovery-errors follow a semi-convergent nature. That is, the iteration index $k$ plays the role of a regularization parameter, where for smaller values of $k$ the recovered solution $x_k^\delta$ is over-regularized (over-smoothed) and for larger values of k, $x_k^\delta$ is under-regularized (over-fitted). Therefore, for ``an appropriate" index $k(\delta)$, one would obtain a regularized solution $x_{k(\delta)}^\delta$ which can approximate the true solution of \eqref{Inv. prob.}, where an appropriate index depends on the choices of the selection criteria, this is discussed in details in successive sections. Also, instead of a constant step-size in \eqref{LS descent}, one can have varying step-sizes $\tau_k \geq 0$ for each iteration and can even generalize the simple gradient descent method to much faster methods, such as Krylov-Subspace or Conjugate-gradient semi-iterative methods, for details and generalizations see \cite{Hanke1991, Landweber, Hanke_Neubauer_Scherzer, Engl+Hanke+Neubauer}. However, as one can notice, through this approach (or its extensions, which one improves the speed of convergence) one cannot make use of any prior knowledge of the solution $\hat{x}$, i.e., no structures can be imposed on the recovery process. 
\subsection{Regularized solution vs. Regularized solution family}\label{Reg. family sol.}
\mbox{ }

Before we jump into our interpretation of the PnP-algorithms, we make a quick inspection on the pros and cons of the aforementioned regularization methods, as this will help in better understanding of our formulation. The similarity between both these methods is that, first, one generates a ``family of regularized-solutions" (depending on a parameter) and then, an ``appropriate regularized solution ($x^\delta$)" is selected form that class, based on some a-priori or a-posteriori selection rules or criteria, depending on the noise-level $\delta$. In the dissimilarity, the constraint-regularization offers a much larger set of regularized solutions which depends not only on the choices of $\lambda$ but also on $\mcal{R}$ (the constraint/regularizer), i.e., the set (family) of regularized solutions in Tikhonov-type methods is
\begin{equation}\label{Tik. family sol.}
\mcal{T}_1 := \left\{ x^\delta(\mcal{D},\lambda,\mcal{R}): x^\delta(\mcal{D},\lambda,\mcal{R}) \mbox{ is a minimizer of \eqref{Gen. Tik. fun.}, } \lambda = \lambda(\delta) \geq 0 \right\}, 
\end{equation}
where as, the classical (simple) semi-iterative regularization can only provide a set of regularized solutions depending on the iteration index $k$, i.e.,
\begin{equation}\label{LS. family sol.}
\mcal{I}_2 := \left\{ x^\delta(\mcal{D},k) : x_{k}^\delta = x_{k-1}^\delta - \tau_{k-1} \nabla_x \mcal{D}(x_{k-1}^\delta), \;\; 1\leq k \leq k(\delta) \leq \infty  \right\}.
\end{equation}
However, on the plus side, a regularized solution in the family \eqref{LS. family sol.} can be obtained much faster (though without any structures in it) than a regularized solution in \eqref{Tik. family sol.}, as in $\mcal{T}_1$ one has to minimize completely the associated functional in \eqref{Gen. Tik. fun.}.

\begin{remark}
Note that, the formulations \eqref{Gen. Tik. fun.} and \eqref{LS descent} only yield corresponding families of regularized solutions $\mcal{T}_1$ and $\mcal{I}_2$, respectively, and not a solution of \eqref{Inv. prob.}. One estimate the true solution $\hat{x}$ through these families of solutions depending on some ``parameter choice or selection criteria". That is, based on a specific criterion, say $\mcal{S}_0$, one chooses a regularized solution from either the family \eqref{Tik. family sol.}, $x_0^\delta(\mcal{D},\lambda(\delta),\mcal{R}) \in \mcal{T}_1$, such that $\delta \rightarrow 0$ implies $\lambda(\delta;\mcal{S}_0) \rightarrow 0$ and $x_0^\delta(\mcal{D},\lambda(\delta),\mcal{R}) \rightarrow \hat{x}$, or from the family \eqref{LS. family sol.}, $x_0^\delta(\mcal{D},k(\delta)) \in \mcal{I}_1$, such that $\delta \rightarrow 0$ implies $k(\delta;\mcal{S}_0) \rightarrow \infty$ and $x_0^\delta(\mcal{D},k(\delta)) \rightarrow \hat{x}$. Of course, when the selection criterion changes, the regularized solution changes too and hence, the families \eqref{Tik. family sol.} and \eqref{LS. family sol.} of regularized solutions get extended by the parameter choice criterion, i.e.,
\begin{align}\label{Tik. family sol. ext.}
\mcal{T}_2 := \left\{ x^\delta(\mcal{D},\lambda,\mcal{R}; \mcal{S}) \right. :& \left. x^\delta(\mcal{D},\lambda,\mcal{R}) \mbox{ is a minimizer of \eqref{Gen. Tik. fun.}, } \lambda = \lambda(\delta,\mcal{S})\geq 0, \right. \notag \\
 & \left. \mbox{and $\mcal{S}$ is a selection criterion} \right\},
\end{align}
\begin{align}\label{LS. family sol. ext.}
\mcal{I}_2 := \left\{ x^\delta(\mcal{D},k;\mcal{S}) \right. \;  :& \; \left. x_{k}^\delta = x_{k-1}^\delta - \tau_{k-1} \nabla_x \mcal{D}(x_{k-1}^\delta), \;\; 1\leq k \leq k(\delta,\mcal{S}) \leq  \infty, \right. \notag\\
& \left. \mbox{ and $\mcal{S}$ is a selection criterion}  \right\}.
\end{align}
Note that, for practical problems, it's near impossible to obtain the optimal parameters for the respective families, since the true solution $\hat{x}$ is unknown. Although the problem of choosing an appropriate selection criterion is non-trivial, especially for large scale problems with unknown error, there has been a lot of studies in the literature in this regard, for example Cross-Validation (CV), Generalized Cross-Validation (GCV), Discrepancy principle (DP), L-curve, monotone error rule etc., see \cite{Morozov1966,Vainikko1982, Gfrerer1987, Hansen1992, Lawson_Hanson,Tautenhahn_Hamarik,Bauer_Hohage, Mathe}.
\end{remark}

\begin{remark}
Now, one may be tempted to fuse the above two families naively, for example, one can incorporate the iteration index $k$ as an additional regularization parameter in \eqref{Tik. family sol. ext.} and can define a regularized solution, $x^\delta(\mcal{D},\lambda,\mcal{R},k;\mcal{S})$, as the k-th iterate $x_{k}^\delta(\mcal{D},\lambda,\mcal{R})$ when minimizing \eqref{Gen. Tik. fun.}, as defined in \eqref{Tik. Iter.}, which satisfies the selection criterion $\mcal{S}$, for fixed $(\mcal{D},\lambda,\mcal{R};\mcal{S})$. However, under this definition, the new regularized solution $x_{k}^\delta(\mcal{D},\lambda,\mcal{R})$ is neither a solution (minimizer) corresponding to a Variational problem \eqref{Gen. Tik. fun.}, nor a solution corresponding to any classical semi-iterative regularization method, as the gradient $\nabla_x F(x_k^\delta;\mcal{D},\lambda,\mcal{R})$, defined in \eqref{Tik. Iter.} for all $k$, may not be a descent direction for minimizing the LS-functional, defined in \eqref{LS functional}. That is, $x^\delta(\mcal{D},\lambda,\mcal{R},k;\mcal{S})$ does not belong to either of the two classes, $\mcal{T}_2$ or $\mcal{I}_2$, and hence, we need to further expand these classes. Moreover, this definition does not work for PnP-algorithms, for any general denoiser, as not all PnP-algorithms can be formulated as a Variational problem. Nevertheless, the advantage of such a definition is that, one doesn't have to worry about the convergence of the iterates $x_{k}^\delta(\mcal{D},\lambda,\mcal{R})$, since, if it does not belong to the family $\mcal{I}_3$ then we are not looking for the solution of \eqref{Gen. Tik. fun.} anymore. In addition, this also speeds up the recovery process, as one does not have to iterate indefinitely for convergence. This make sense, to a certain extent, since we are interested in the solution of \eqref{Inv. prob.} and the  formulation \eqref{Gen. Tik. fun.} only aids us to approximate the true solution $\hat{x}$, through the regularization parameter $\lambda(\delta,\mcal{S})$. Hence, even when the iterates \eqref{Tik. Iter.} converge (for fixed $(\mcal{D},\lambda,\mcal{R})$) to $x^\delta_\lambda$, a minimizer of \eqref{Gen. Tik. fun.}, $x^\delta_\lambda$ is not the solution of \eqref{Inv. prob.}, rather, one typically generates a list of $\{x^\delta_\lambda\}_\lambda$ and an appropriate estimate  of $\hat{x}$ is then selected based on the selection criterion $\mcal{S}$. This naive idea lays the ground for our formulation to interpret the PnP-algorithms, but has to be polished properly for meaningful definitions, which is done in the following sections.
\end{remark}

\subsection{Structured iterations}\label{S. Structured iterations}
\mbox{ }

\begin{figure}
    \begin{subfigure}{0.495\textwidth}
        \includegraphics[width=\textwidth]{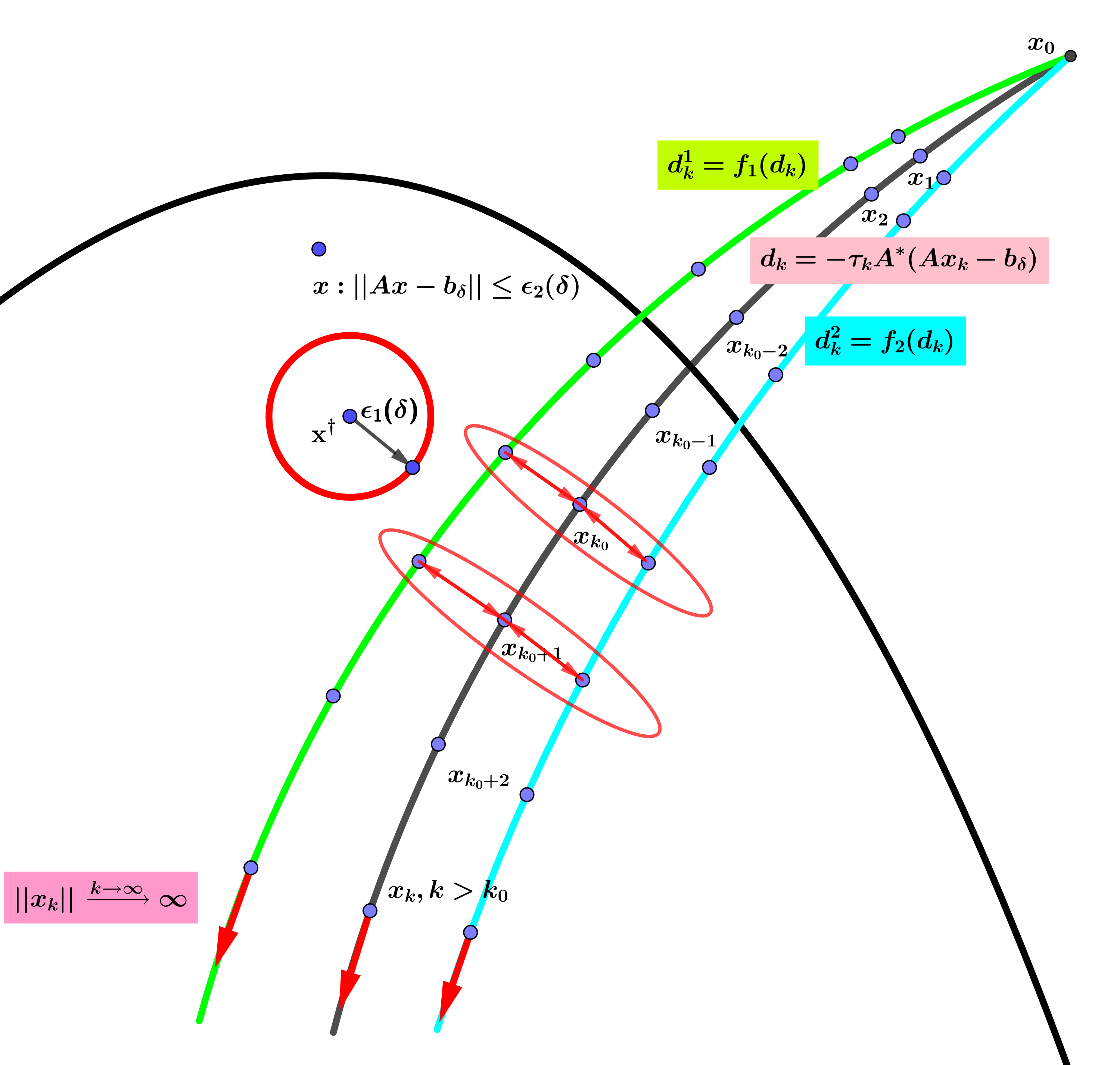}
        \caption{Semi-iterative regularization for different descent directions $d_k^\delta$.}
        \label{graph of h_k vs. d_k}
    \end{subfigure}
    \begin{subfigure}{0.495\textwidth}
        \includegraphics[width=\textwidth]{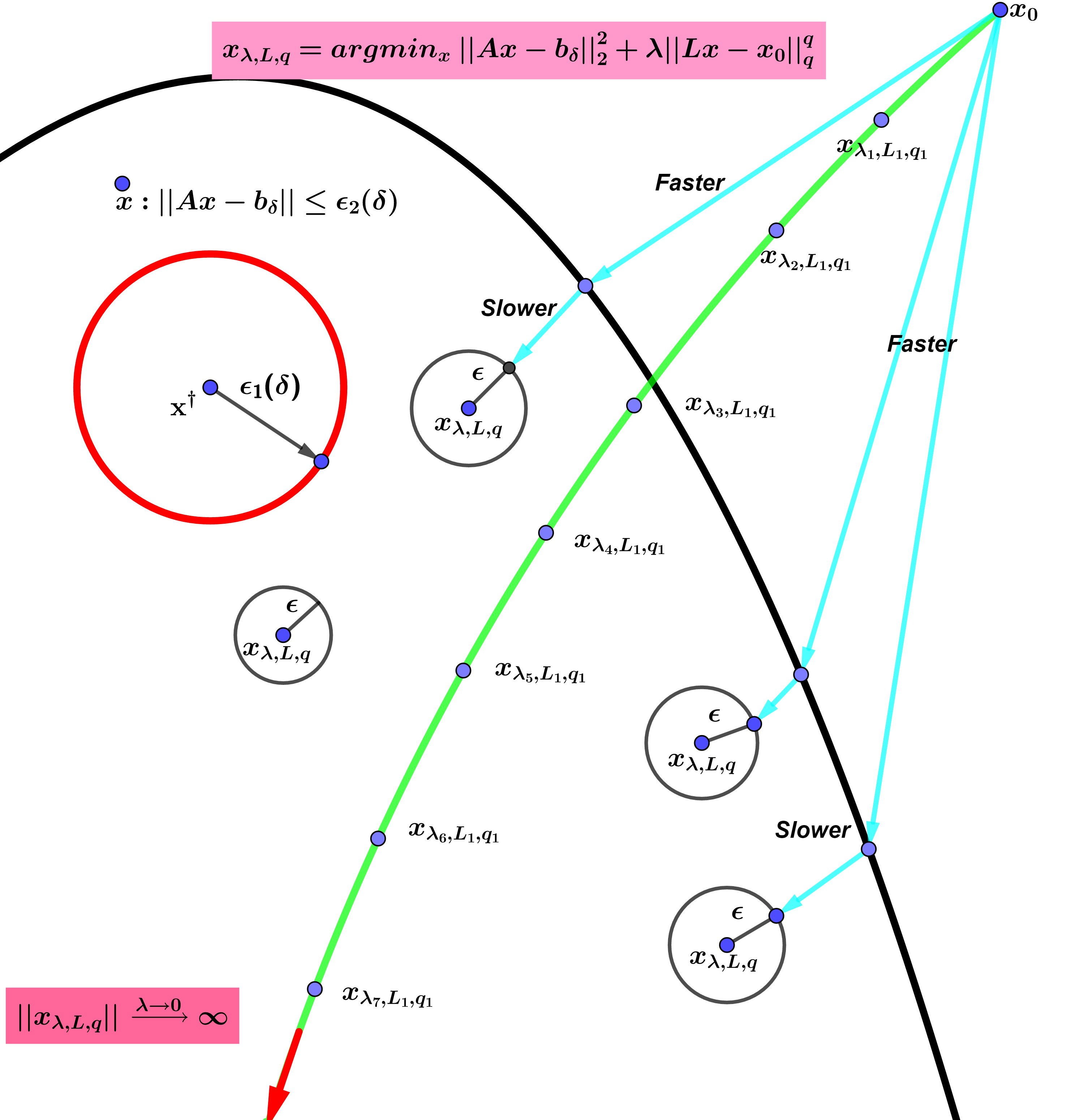}
        \caption{Tikhonov-type regularization for different regularizers $\mcal{R}$.}
        \label{graph of Tikhonov path}
    \end{subfigure}
    \caption{Semi-iterative vs. Tikhonov-type regularization.} 
    \label{Semi-iterative vs. Tikhonov regularization}
\end{figure}

Here, instead of expanding the family $\mcal{T}_2$ to incorporate the PnP-algorithms, we expand the family $\mcal{I}_2$ to address the above problems. First, note that, when minimizing the LS-functional \eqref{LS functional}, $d_k^\delta$ in \eqref{LS descent} needs to be a descent direction only and doesn't have to be the negative gradient (which is the steepest descent direction). Hence, any direction with the following property can be considered as a descent direction, for simplicity we choose $\mcal{D}(x)=||Ax - b_\delta||_2^2$, 
\begin{equation}\label{LS descend prop.}
	\lprod{d_k^\delta}{-\nabla_x \mcal{D}(x_k^\delta)} = \lprod{d_k^\delta}{-A^*(Ax_k^\delta - b_\delta)} > 0,
\end{equation}
where $\lprod{.}{.}$ is the associated $\ell_2$-product. It's also easy to verify that if
\begin{equation}\label{LS proximity prop.}
||e_k^\delta + \tau \nabla_x \mcal{D}(x_k^\delta)||_2^2 \; \leq \; \epsilon_k(\delta),
\end{equation}
for some $\epsilon_k(\delta) \geq 0$, then, starting from initial points $z_0^\delta$ and $x_0^\delta$,
\begin{align}
||z_k^\delta - x_k^\delta||_2^2 \; &\leq \; ||z_0^\delta - x_0^\delta||_2^2 \; + \; \sum_{i=0}^{k-1}||e_i^\delta + \tau \nabla_x \mcal{D}(x_i^\delta)||_2^2 \; \leq \; ||z_0^\delta - x_0^\delta||_2^2 \; + \; \sum_{i=0}^{k-1} \epsilon_i(\delta)
\end{align}
where $z_i^\delta = z_{i-1}^\delta + e_i^\delta$ and $x_i^\delta = x_{i-1}^\delta - \tau \nabla_x \mcal{D}(x_{i-1}^\delta)$. Hence, if for some $\epsilon(\delta) \geq 0$, $||z_0^\delta - x_0^\delta||_2^2 \leq \epsilon(\delta)$ and $\epsilon_i(\delta) \leq \epsilon(\delta)$, for $0\leq i \leq k-1$, then we have
\begin{equation}\label{LS z_k - x_k}
||z_k^\delta - x_k^\delta||_2^2 \; \leq \; C(k) \; \epsilon(\delta),
\end{equation}
where the (increasing) constant $C(k)$ depends only on the iteration index $k$. In other words, if any direction $d_k^\delta$ satisfies conditions \eqref{LS descend prop.} and \eqref{LS proximity prop.}, then the iterates $z_k^\delta$, formed based on it, also remains close to the LS iterates $x_k^\delta$, upto some constant depending on the iteration index $k$. And, since in a semi-iterative regularization $k(\delta)$ is bounded away from infinity, the new iterates $z_k^\delta$ are bounded at most $C(k(\delta))\epsilon(\delta)$ away from the LS iterates $x_k^\delta$. It is also easy to check that, for the residues or discrepancy terms $\mcal{D}(x) = ||Ax - b_\delta||_2^2$, we have
\begin{align}\label{LS residue diff.}
||Az_k^\delta - b_\delta||_2^2 \leq \; ||Ax_k^\delta - b_\delta||_2^2 \; + \; C(k)\epsilon(\delta)||A||_2^2, 
\end{align}
that is, if the residues of the LS iterates are small ($||Ax_k^\delta - b_\delta||_2^2 \approx \delta$), then so are the residues of $z_k^\delta$ with an error term $C(k)\epsilon(\delta)||A||_2^2$, where $A$ is usually bounded. This is helpful when the selection criterion, such as discrepancy principle, depends on the residue norms, i.e., terminate the iterations when the residue norms of the iterates $x_k^\delta$ is around certain limit ($||Ax_k^\delta - b_\delta||_2^2 \approx \delta$), and \eqref{LS residue diff.} implies that for those iterations the residue norms corresponding to the new iterates $z_k^\delta$ is also around that limit, with an error term, implying that the iterates $z_k^\delta$ have also fitted the data to that extent. Hence, when minimizing the LS-functional \eqref{LS functional}, replacing the vanilla descent directions ($-\nabla_x D(x_k^\delta)$) with an ``appropriate" descent direction $d_k^\delta$, one can impose certain structures in the recovery process through the iterates $z_k^\delta$, and thus, can recover a regularized structured solution $z^\delta_{k(\delta)}$. With this formulation, the family of regularized solution $\mcal{I}_2$, as defined in \eqref{LS. family sol. ext.}, is expanded to
\begin{align}\label{LS. family I_3}
\mcal{I}_3 := &\left\{ x^\delta(\mcal{D},k,d_k^\delta;\mcal{S}) \right. \;  : \; \left. x_{k}^\delta = x_{k-1}^\delta \;+\; d_{k-1}^\delta, \;\; 1\leq k \leq k(\delta,\mcal{S}) \leq  \infty, \right.\\
& \left. \; \mbox{s.t.}\; \lprod{d_{k-1}^\delta}{- \tau_{k-1} \nabla_x \mcal{D}(x_{k-1}^\delta)}>0 \mbox{ and $\mcal{S}$ is a selection criterion.}  \right\}. \notag
\end{align}
Figure \ref{Semi-iterative vs. Tikhonov regularization} elaborates the families of solutions, $\mcal{I}_3$ and $\mcal{T}_2$, visually.
Note that, although $\mcal{I}_3 \supseteq \mcal{I}_2$, $\mcal{I}_3$ is neither a super-set, nor a subset, nor equals to $\mcal{T}_2$. However, for a fixed $\mcal{S}_0$, any $x^\delta(\mcal{D},\lambda(\delta),\mcal{R}) \in \mcal{T}_2$ can be ``approximated" by a $x^\delta(\mcal{D},k(\delta),d_k^\delta) \in \mcal{I}_3$, where the extent of approximation depends on the choice of $\lambda(\delta)$, since, while minimizing \eqref{Gen. Tik. fun.} using the FBS approach \eqref{FBS min.}, the resulting direction at the k-1 step is given by
\begin{align}\label{FBS prox. dir.}
d_{k-1}^\delta &= \mbox{Prox}_{\lambda \mcal{R}}(x_k^\delta) \; - \; x_{k-1}^\delta \\
&= \left( \mbox{Prox}_{\lambda \mcal{R}}(x_k^\delta) - x_k^\delta \right) + \left( -\tau \nabla_x \mcal{D}(x_{k-1}^\delta) \right),  \notag
\end{align} 
and hence, if 
\begin{equation}\label{Prox(x) - x}
|| \mbox{Prox}_{\lambda \mcal{R}}(x_k^\delta) - x_k^\delta ||_2 \; < \; || -\tau \nabla_x \mcal{D}(x_{k-1}^\delta)||_2,
\end{equation}
we have
\begin{align}
\lprod{d_{k-1}^\delta}{-\tau \nabla_x \mcal{D}(x_{k-1}^\delta)} \geq \tau ||\nabla_x \mcal{D}(x_{k-1}^\delta)||\left( \tau ||\nabla_x \mcal{D}(x_{k-1}^\delta)|| - || \mbox{Prox}_{\lambda \mcal{R}}(x_k^\delta) - x_k^\delta || \right) > 0, \notag
\end{align}
i.e., $d_{k-1}^\delta$ is a descent direction if \eqref{Prox(x) - x} holds. But, as the residues $\mcal{D}(x_k^\delta)$ get smaller, so does the gradient $-\tau \ \nabla_x \mcal{D}(x_k^\delta)$ for a fixed $\tau$, 
 and hence, \eqref{Prox(x) - x} won't hold for such iterates, as well as, for large values of $\lambda$. Therefore, for a given $\lambda(\delta)$, the largest iterate $k(\delta)$ for which \eqref{Prox(x) - x} holds, yields a regularized solution $x^\delta(\mcal{D},k(\delta),d_k^\delta) \in \mcal{I}_3$ which can be close to $x^\delta(\mcal{D},\lambda(\delta),\mcal{R}) \in \mcal{T}_2$, depending on $\lambda(\delta)$. 

\begin{remark} \label{remark struc. iter.}
Note that, if a direction $d_k^\delta$ satisfies \eqref{LS descend prop.} for all $k$, then the corresponding recovery-errors for the iterates $z_k^\delta$ will also follow a semi-convergent nature. However, the parameter choice value $k(\delta)$, for a fixed selection criterion $\mcal{S}_0$, will be different for the respective iterates. Also, equation \eqref{LS z_k - x_k} does not suggest that $||z_{k(\delta)}^\delta - \hat{x}|| \leq ||x_{k(\delta)}^\delta - \hat{x}||$, where $\hat{x}$ is the true solution, i.e., $z_{k(\delta)}^\delta$ is a better estimate of $\hat{x}$ than $x_k^\delta$, but, if we believe that $d_k^\delta$ is a ``better" direction than $-\tau \nabla_x \mcal{D}(x_k^\delta)$, then we would expect that the iterates $z_k^\delta$ to approximate the true solution better than $x_k^\delta$. In any case, in the absence of the true solution, a ``better solution" for a fixed $\mcal{S}_0$ is the one that satisfies $\mcal{S}_0$ better. 

Moreover, for $z_{k(\delta)}^\delta$ to be a regularized solution, $z_{k(\delta)}^\delta \rightarrow \hat{x}$ when $\delta \rightarrow 0$, and this is true from \eqref{LS z_k - x_k} provided $C(k)\epsilon(\delta) \rightarrow 0$. That is, for a fixed selection criterion $\mcal{S}_0$, and hence $k(\delta),$ $\epsilon(\delta) = o(1/k^2(\delta))$. In other words, as the noise level decreases the descent direction should get closer to the negative gradient direction, i.e., when $\delta \rightarrow 0$, $d_k^\delta \rightarrow -\tau \nabla_x \mcal{D}(x_k^\delta) \Leftrightarrow z_k^\delta \rightarrow x_k^\delta$, and $k_i(\delta) \rightarrow \infty$, for $i=1,2$, which results in $x_{k_1(\delta)}^\delta \rightarrow \hat{x} \leftarrow z_{k_2(\delta)}^\delta$. 

The natural question that one can ask now is, when can a regularized structured solution be obtained without encountering the semi-convergence in the recovery errors? The answer is `never' for those structure imposing directions $d_k^\delta$ that satisfies \eqref{LS descend prop.} for all $k$, as then one minimizes the LS functional \eqref{LS functional} for noisy $b_\delta$, leading to the semi-convergence of the recovery errors. However, one can relax the condition \eqref{LS descend prop.} appropriately to circumvent the semi-convergence as well as recover a structured regularized solution, but then, it won't fall in the family $\mcal{I}_3$. Hence, one has to define an extension of the family $\mcal{I}_3$ to encompass such structured regularized solutions, details in the later sections.

\end{remark}

\section{\textbf{PnP-algorithms as structured iterations}}
The core of any PnP-algorithm is the \textit{data-consistency step} followed by the \textit{data-denoising step}, where the denoising operator doesn't need to be a proximal operator anymore. We first analyze the FBS-PnP algorithm and then the ADMM-PnP algorithm.

\subsection{Forward-backward splitting (FBS) - PnP}\label{Sec. FBS-PnP}
\mbox{ }

In this algorithm, for a fixed denoiser $H_\sigma$ (of denoising strength corresponding to noise level $\sigma$) and starting from an initial choice $z_0^\delta$, at any iteration step $k\geq 1$, we have 
\begin{align}\label{FBS-PnP}
z_{k-1}^\delta &\longmapsto \; x_k^\delta \; = \; z_{k-1}^\delta - \tau_k \nabla_x \mcal{D}(z_{k-1}^\delta) \; &\longleftarrow \mbox{ data-consistency step} \\
x_k^\delta &\longmapsto \; z_k^\delta \; = \; H_{\sigma_k}(x_k^\delta), \; &\longleftarrow \mbox{ data-denoising step},
\end{align}
where $H_{\sigma_k}$ is the updated kth denoiser, with the denoising strength corresponding to $\sigma_k = \tau_k \sigma$. Hence, an approach similar to \eqref{FBS prox. dir.}, provides a structure imposing direction $d_{k-1}^\delta$,  which can be defied as
\begin{align}\label{FBS-PnP direction}
d_{k-1}^\delta &:= z_k^\delta - z_{k-1}^\delta \\ 
&= \underbrace{-\tau \nabla_x \mcal{D}(z_{k-1}^\delta)}_{data-consistency} \; + \; \underbrace{\left(H_{\sigma_k}(x_k^\delta) - x_k^\delta\right)}_{data-denoising} \;, \notag
\end{align}
and the direction $d_{k-1}^\delta$, as defined in \eqref{FBS-PnP direction}, will be a descent direction provided it satisfies \eqref{LS descend prop.}. Therefore, for $H_{\sigma_k}$ satisfying
\begin{equation}\label{H(x) - x}
||H_{\sigma_k}(x_k^\delta) - x_k^\delta|| \; < \; ||-\tau \nabla_x \mcal{D}(x_{k-1}^\delta)||,
\end{equation}
we have $d_{k-1}^\delta$ satisfying \eqref{LS descend prop.}, since
\begin{align}
\lprod{d_{k-1}^\delta}{-\tau \nabla_x \mcal{D}(z_{k-1}^\delta)} \geq \tau ||\nabla_x \mcal{D}(z_{k-1}^\delta)||\left( \tau ||\nabla_x \mcal{D}(z_{k-1}^\delta)|| - ||H_{\sigma_k}(x_k^\delta) - x_k^\delta|| \right) > 0. \notag
\end{align}
In other words, \eqref{H(x) - x} suggests that if the change (denoising) in the denoising step is smaller than the change (improvement) in the data-consistency step, than the resulting direction $d_{k-1}^\delta$ serves as a \textit{structure imposing descent direction}, as defined in \S \ref{S. Structured iterations}, and hence, the iterates  $z_k^\delta$ recovers a structured regularized-solution $z_{k(\delta)}^\delta$, depending on the selection criteria $\mcal{S}$. This makes sense, since, if the strength of the denoiser exceeds the changes arising in the data-consistency step then it will over-denoise or over-smooth the iterates, leading to an over-smoothed recovery. Of course, if $||H_{\sigma_k}(x_k^\delta) - x_k^\delta||$ is very small then, from \eqref{FBS-PnP direction}, $d_{k-1}^\delta$ will be very close to $-\tau \nabla_x \mcal{D}(z_{k-1}^\delta)$ and thus, the structured iterates $z_k^\delta$ will also be close to the unstructured (noisy) iterates $x_k^\delta$, leading to an under-denoised or under-smoothed recovery. Note that, the condition \eqref{H(x) - x} is sufficient but not necessary for the directions $d_k^\delta$ to be a descent direction, i.e., $d_k^\delta$ can satisfy \eqref{LS descend prop.} even when it does not satisfy \eqref{H(x) - x}.

Now, similar to \eqref{LS z_k - x_k}, we would like to estimate a bound on the difference of the new denoised iterates $z_k^\delta$ with the original LS-iterates $y_k^\delta$ (generated from the negative gradients), i.e., starting from the same initial point $z_0^\delta = y_0^\delta \equiv 0$ and with the same step-size $\tau > 0$, we have $x_1^\delta = y_1^\delta$ and
\begin{align}\label{FBS-PnP z_k - y_k}
||z_1^\delta - y_1^\delta|| &= ||H_{\sigma_1}(x_1^\delta) - x_1^\delta||, \notag \\
||z_2^\delta - y_2^\delta|| &\leq ||H_{\sigma_2}(x_2^\delta) - x_2^\delta|| + (1+\tau||A^*A||)||H_{\sigma_1}(x_1^\delta) - x_1^\delta|| \notag \\
||z_3^\delta - y_3^\delta|| &\leq ||H_{\sigma_3}(x_3^\delta) - x_3^\delta|| + (1+\tau||A^*A||)||z_2^\delta - x_2^\delta|| \notag\\
&\leq ||H_{\sigma_3}(x_3^\delta) - x_3^\delta|| + (1+\tau||A^*A||)||H_{\sigma_2}(x_2^\delta) - x_2^\delta|| \notag \\
&\hspace{3cm}+ (1+\tau||A^*A||)^2||H_{\sigma_1}(x_1^\delta) - x_1^\delta|| \notag \\
&\vdots \notag \\
||z_k^\delta - y_k^\delta|| &\leq \sum_{i=0}^{k-1} (1 + \tau||A^*A||)^{i} ||H_{\sigma_{k-i}}(x_{k-i}^\delta) - x_{k-i}^\delta||.
\end{align}
That is, for the $k(\delta)$ iteration the denoised iterate $z_{k(\delta)}^\delta$ will be bounded away from the LS iterate $y_{k(\delta)}^\delta$ at most by the factor given in \eqref{FBS-PnP z_k - y_k}, which also depends on the denoising strength of the denoiser $H_{\sigma_k}$ at each iteration $k$. Again, this does not guarantee that $||z_{k(\delta)}^\delta - \hat{x}|| \leq ||y_{k(\delta)}^\delta - \hat{x}||$, as it can happen that the iterates may be over-smoothed (for stronger denoiser $H_{\sigma_k}$) and might have crossed the optimal solution. However, if we trust (tune) our denoiser to be ``appropriate" and directing the iterates $x_k^\delta$ towards the true solution $\hat{x}$ through $z_k^\delta$, i.e., towards a ``better" direction, then we would hope the shift in \eqref{FBS-PnP z_k - y_k} to be in the proper direction, and hence, leading to $||z_{k(\delta)}^\delta - \hat{x}|| \leq ||y_{k(\delta)}^\delta - \hat{x}||$, i.e., a better recovery. In \S \ref{Sec. weaker deno.} and \S \ref{Sec. stronger deno.} we provide some suggestions to tune a bad denoiser $H_\sigma$ appropriately, i.e., when the original denoiser $H_\sigma$ is too strong or too weak, then how can it be improved to perform the appropriate level of denoising to recover better estimates of the true solution, without actually altering the original denoiser.

In this section we showed that, if the denoiser $H_{\sigma_k}$ satisfies the descent condition \eqref{H(x) - x} at each iteration k, then the FBS-PnP algorithm will fall in the family $\mcal{I}_3$. The weakness, which also the strength, of the condition \eqref{H(x) - x} is that, the recovery errors corresponding to the iterates $z_k^\delta$ will also follow a semi-convergent path, and hence, the iterations have to be terminated at an appropriate early instance ($k(\delta)$), as explained in remark \ref{remark struc. iter.}. To overcome this inconvenience, as stated in remark \ref{remark struc. iter.}, the family $\mcal{I}_3$ needs to be further extended, which is done later.

\subsubsection{\textbf{Reinforcing weaker denoisers}}\label{Sec. weaker deno.}
\mbox{ }

In practice, one usually starts with a denoiser of certain strength, say $\sigma$, and at every iteration $k$, the denoising of $x_k^\delta$ is done via $H_{\sigma_k}$, where $\sigma_k = \sigma \tau_k$, where $\tau_k$ is the associated step-size for that iteration. As discussed in \S \ref{Sec. FBS-PnP} the FBS-PnP algorithm is guaranteed to fall in the regularized family $\mcal{I}_3$ provided the denoiser $H_{\sigma_k}$ satisfy \eqref{H(x) - x}. Now, for a weaker denoiser (i.e., $||H_{\sigma_k}(x_k^\delta) - x_k^\delta||$ is small), though the structure imposing direction $d_k^\delta$, as defined in \eqref{FBS-PnP direction}, will satisfy \eqref{H(x) - x}, the iterates $x_k^\delta$ won't be sufficiently denoised, i.e., the denoised iterates $z_k^\delta$ will still contain certain amount of noise. Of course, during the next minimization process, one can opt for an ``appropriate" stronger denoiser $H_{\sigma_k}$, if possible, to provide stronger denoising of the iterates $x_k^\delta$. There is also a second approach to augment the denoising strength of $d_k^\delta$, without altering the denoiser, for example, when there is only one denoiser to work with. Observe that, the satisfaction of the condition \eqref{H(x) - x}, by a denoiser $H_{\sigma_k}$, is also dependent on the step-size $\tau$. That is, even for a weaker denoiser $H_{\sigma_k}$, if the value of  $\tau$ is small enough, then $H_{\sigma_k}$ will fail to satisfy the condition \eqref{H(x) - x}, i.e., the denoising strength of the weaker denoiser will increase implicitly. Hence, by decreasing the step-size $\tau$ smaller than the previous choice, such that $||H_{\sigma_k}(x_k^\delta) - x_k^\delta||$ is only slightly smaller than $||-\tau \nabla_x \mcal{D}(z_{k-1}^\delta)||$ will, not only lead to a descent direction but also, the amount of denoising will be increased for the original denoiser $H_{\sigma}$, see Example \ref{Ex. weak deno.}. However, this does increase the computational time, as with smaller values of $\tau$ the descent rate also decreases.

\subsubsection{\textbf{Attenuating stronger denoisers}}\label{Sec. stronger deno.}
\mbox{ }

In contrast to a weak denoiser, a stronger denoisers $H_{\sigma}$ may not lead to the satisfaction of the condition \eqref{H(x) - x}, since, the extent of denoising $||H_{\sigma_k}(x_k^\delta) - x_k^\delta||$ may exceed the improvements coming from the data-consistency step $||x_k^\delta - z_{k-1}^\delta||$. Again, the obvious way to overcome this issue is by opting for a weaker denoiser $H_{\sigma}$, the next time, such that the condition \eqref{H(x) - x} is satisfied. However, this can lead to several attempts of numerous minimization process to find an appropriate denoiser, which can substantially increase the computational time. Note that, unlike the previous tweak, here one cannot increase the step-size to compensate the strong denoising in \eqref{H(x) - x}, as then the descent process will break down, for larger step-sizes. A cleverer approach of handling this issue is through relaxing the denoising strength of the strong denoiser $H_\sigma$ by incorporating an additional relaxation parameter $\alpha \in (0, 1]$, via
\begin{align}
z_k^\delta = z_k^\delta(\alpha) &:= \; x_k^\delta \; + \; \alpha \left( H_{\sigma_k}(x_k^\delta) \; - \; x_k^\delta \right) \\
&= (1 - \alpha)\;x_k^\delta \; + \; \alpha \; H_{\sigma_k}(x_k^\delta). \notag
\end{align}
Hence, if $||H_{\sigma_k}(x_k^\delta) - x_k^\delta|| > ||x_k^\delta - z_{k-1}^\delta|| = ||-\tau \nabla_x \mcal{D}(z_{k-1}^\delta)||$, then for any
\begin{equation}
\alpha \;\leq\; \frac{||x_k^\delta - z_{k-1}^\delta||}{||z_k^\delta - x_k^\delta||} \;=\; \frac{||-\tau \nabla_x \mcal{D}(z_{k-1}^\delta)||}{||H_{\sigma_k}(x_k^\delta) - x_k^\delta||} \; < 1,
\end{equation}
we will have the resulting direction
\begin{align}\label{Strong deno. alpha1}
d_k^\delta &= z_k^\delta - z_{k-1}^\delta \notag \\
&= (x_k^\delta - z_{k-1}^\delta) \; + \; \alpha \left( H_{\sigma_k}(x_k^\delta) - x_k^\delta \right)
\end{align}
satisfying the condition \eqref{LS descend prop.}. Therefore, if $\alpha << 1$ (small) then $d_k^\delta$ is close to $-\tau \nabla_x \mcal{D}(z_{k-1}^\delta)$, which leads to weak denoising, and if $\alpha \approx 1$ (close to 1) then it leads to strong denoising. This can be seen in Table \ref{Table Strong deno.}, for Example \ref{Ex. Strong deno.}. Usually, one can choose $\alpha = \gamma \; \frac{||-\tau \nabla_x \mcal{D}(z_{k-1}^\delta)||}{||H_{\sigma_k}(x_k^\delta) - x_k^\delta||}$, where the proximity parameter $\gamma \leq 1$ controls the closeness of $d_k^\delta(\alpha(\gamma))$ to $-\tau \nabla_x \mcal{D}(z_{k-1}^\delta)$ through the parameter $\alpha(\gamma)$.

An alternative way of estimating a proper $\alpha$ value is through minimizing $z_k^\delta(\alpha)$ based on the selection criterion, i.e.,
\begin{align}\label{Strong deno. alpha2}
\alpha_0 := &\argmin_{\alpha \in (0,1]} \;\;\; \mcal{S}(z_k^\delta(\alpha))\\
&\mbox{such that,} \;\; z_k^\delta(\alpha) =  x_k^\delta +  \alpha \left( H_{\sigma_k}(x_k^\delta) \; - \; x_k^\delta \right). \notag
\end{align}
Note that, the minimization problem \eqref{Strong deno. alpha2} may not be a strictly convex, i.e., there might not be a global minimizer $\alpha_0$. Nevertheless, this is simply a sub-problem intended to find an appropriate $\alpha$ value between 0 and 1, depending on the selection criterion $\mcal{S}$, and hence, even if the best $\alpha_0$ is not obtained (which is not a requirement), any $\alpha \in (0, 1)$ will reduce the denoising strength and empirical experiments (see Example \ref{Ex. Strong deno.}) show that it works fine. Also, the alpha value from \eqref{Strong deno. alpha2} doesn't guarantee that $d_k^\delta(\alpha_0) := z_k^\delta(\alpha_0) - z_{k-1}^\delta$ will satisfy \eqref{LS descend prop.}, but, it does attenuates the strong denoising of $H_{\sigma_k}$ and provides an iterate $z_k^\delta$ which satisfies the selection criterion the best. 

\begin{remark}
Note that, inspired from \eqref{Strong deno. alpha2}, one can also choose the denoising strength of the denoiser $H_{\sigma_k}$, at every iteration, depending on the selection criterion, i.e., by defining the denoised iterates as a function of the denoising strength $\sigma$
\begin{equation}
z_k^\delta(\sigma_k) := H_{\sigma_k}(x_k^\delta),
\end{equation}
one can choose the denoising strength $\sigma_k$, within a certain range, which best satisfies the selection criterion, i.e.,
\begin{align}\label{Min. over sigma}
\sigma_k := &\argmin_{\sigma \in [\sigma_{{min}},\; \sigma_{max}]} \;\;\; \mcal{S}(z_k^\delta(\sigma)) \\
&\mbox{such that, } \;\; z_k^\delta(\sigma) = H_\sigma(x_k^\delta). \notag
\end{align}
However, the denoising step can be computationally expensive step during the iterative process, and hence, performing a minimization over the denoising strength can notably increase the computational time. In addition, similar to \eqref{Strong deno. alpha2}, the minimization in \eqref{Min. over sigma} need not be even convex.
\end{remark}

\subsection{Fast Forward-backward splitting (Fast FBS) - PnP}
\mbox{ }
The plain FBS-PnP method tends to show very slow convergence rate, even with the semi-iterative approach, i.e., $k(\delta)$ is very large. To speed up the descent rate, faster algorithms are usually employed such as, the fast proximal gradient methods (FPGM) or FISTA, see \cite{Beck_Teboulle}. The iterates in this algorithm have an intermediate momentum step to accelerate the descent process, which is given by, starting from initial $z_0^\delta = x_0^\delta$ and $t_0=1$, for $k \geq 1$,
\begin{align} \label{Fast FBS-PnP}
z_{k-1}^\delta &\longmapsto \; x_k^\delta \; = \; z_{k-1}^\delta - \tau_k \nabla_x \mcal{D}(z_{k-1}^\delta) \; &\longleftarrow \mbox{ data-consistency step} \\
x_k^\delta &\longmapsto \; z_k^\delta \; = \; H_{\sigma_k}(x_k^\delta), \; &\longleftarrow \mbox{ data-denoising step} \\
z_k^\delta &\longmapsto \; z_k^\delta \; = \; z_k^\delta \; + \; \alpha_k (z_k^\delta - z_{k-1}^\delta), \; &\longleftarrow \mbox{ momentum-step}
\end{align}
where $\alpha_k = \frac{t_{k-1} - 1}{t_k}$ with $t_k = \frac{ (1 + \sqrt{1 + 4t_{k-1}^2})}{2}$.

Now, one can also perform similar analysis as done for the vanilla FBS-PnP method, we skip the details.

\section{\textbf{Expansion of the regularization family}}

Before we proceed to attach a meaning to the ADMM-PnP algorithms we would like to address the question raised in the last part of remark \ref{remark struc. iter.}, i.e., can the family of regularized solutions $\mcal{I}_3$ be generalized to soften the condition \eqref{LS descend prop.}, which leads to the semi-convergence of the recovery errors? To answer this question let's first interpret the family $\mcal{I}_3$ from a different perspective. Note that, the formulation of $\mcal{I}_3$ in \eqref{LS. family I_3} is also equivalent to 
\begin{gather}\label{LS. family I_3'}
\mcal{I}_3' := \left\{ x^\delta(\mcal{D},k,d_k^\delta;\mcal{S}) : \; x_{k}^\delta \;=\; x_{k-1}^\delta \;+ \; d_{k-1}^\delta, \mbox{ subject to } \mcal{D}(x_k^\delta) \geq \epsilon(\delta,\mcal{S}) \geq 0, \right. \notag \\
 \left. \mbox{such that $d_k^\delta$ satisfies \eqref{LS descend prop.} and $\mcal{S}$ is a selection criterion.}  \right\},
\end{gather}
where the lower bound $\epsilon(\delta,\mcal{S})$ for the discrepancy term $\mcal{D}$ now serves the purpose of $k(\delta)$ and determines the extent of regularization, for example, if $\mcal{S}$ is the Morozov's discrepancy principle (DP), then $\epsilon(\delta) \approx \eta \delta$, for $\eta > 1$, see \cite{Morozov_a, Morozov_b, Morozov1966}. Equivalently, the formulation in $\mcal{I}_3'$ suggests to minimize the data-discrepancy functional $\mcal{D}$, via a (descent) direction $d_k^\delta$ satisfying \eqref{LS descend prop.}, but constraint to $\mcal{D}(x_k^\delta) \geq \epsilon(\delta,\mcal{S})$, i.e., not to minimize $\mcal{D}$ completely, which avoids overfitting to the noisy data $b_\delta$. Therefore, the key is to not overfit the noisy data, which motivates us to generalize the constraint $\mcal{D}(x_k^\delta) \geq \epsilon(\delta)$, for larger iterations ($k > k(\delta)$), to
\begin{equation}\label{DP generalize}
	0 \leq \epsilon_1(\delta) \leq \mcal{D}(x_k^\delta) \leq \epsilon_2(\delta),
\end{equation}
such that $\epsilon_2(\delta) \xrightarrow{\delta \rightarrow 0} 0$ (so that $x_k^\delta \xrightarrow{\delta \rightarrow 0} x^\dagger$), where the lower bound $\epsilon_1(\delta)$ avoids overfitting of the noisy data and the upper bound $\epsilon_2(\delta)$ avoids underfitting of the data. Now, if we enforce the directions $d_k^\delta$ to satisfy \eqref{LS descend prop.} for $x_k^\delta$ satisfying $\mcal{D}(x_k^\delta) > \epsilon_2(\delta)$, we minimize the discrepancy-term $\mcal{D}$, where as, removing the constraint \eqref{LS descend prop.} on those directions $d_k^\delta$ for which $\mcal{D}(x_k^\delta) \leq \epsilon_2(\delta)$, but constraining them implicitly through $x_k^\delta$ such that $\epsilon_1(\delta) \leq \mcal{D}(x_k^\delta + d_k^\delta) \leq \epsilon_2(\delta)$, one is able to relax the descent criterion while (still) avoiding the overfitting to the noisy data. Hence, the new conditions for the directions $d_k^\delta$ are
\begin{gather}\label{descent dir. gen.}
\lprod{d_k^\delta}{-\nabla_x \mcal{D}(x_k^\delta)} > 0, \;\;\; \mbox{ for } \; \mcal{D}(x_k^\delta) > \epsilon_2(\delta) \; \mbox{ or } \; k \leq k(\delta), \mbox{ and}\\
\mcal{D}(x_k^\delta + d_k^\delta) \geq \epsilon_1(\delta), \;\;\;\ \mbox{ for }\; \mcal{D}(x_k^\delta) \leq \epsilon_2(\delta) \; \mbox{ or } \; k > k(\delta). \notag
\end{gather}
Note that, with this modification, the recovery errors for the iterates $x_k^\delta$ does not necessarily follow a semi-convergence trail, nor a convergence, anymore; it may end up oscillating/fluctuating around the optimal solution, see Example \ref{Ex. Strong deno.}. Again, the answer to ``does  it improve the recovery process?", depends on the choice of the recovery directions $d_k^\delta$ and the selection criteria $\mcal{S}$. The formal extension of the family $\mcal{I}_3'$ can be expressed as
\begin{align}\label{LS. family I_4}
\mcal{I}_4 := &\left\{ x^\delta(\mcal{D},k,d_k^\delta;\mcal{S}) \right. \;  : \; \left. x_{k}^\delta = x_{k-1}^\delta \;+\; d_{k-1}^\delta, \; \mbox{ s.t. $d_{k-1}^\delta$ satisfies \eqref{descent dir. gen.},} \right.\\
& \left. \; \mbox{ and $\mcal{S}$ is a selection criterion.}  \right\}. \notag
\end{align}
Observe that, by increasing the number of constraints from one, in equation \eqref{LS. family I_3'}, to two, in equation \eqref{LS. family I_4}, we are able to expand the family of regularized solutions, i.e., for $\epsilon_1(\delta) < \epsilon_2(\delta)$, we have $\mcal{I}_4 \supset \mcal{I}_3'$ and when $\epsilon_1(\delta) = \epsilon_2(\delta)$, then $\mcal{I}_4 = \mcal{I}_3'$. Moreover, the family $\mcal{I}_4 \supseteq \mcal{T}_2$ (the family of Tikhonov-regularization), since, the discrepancy of a regularized solution in $\mcal{T}_2$ is always positive, for any $\lambda(\delta)>0$, i.e., $\mcal{D}(x^\delta(\mcal{D},\lambda(\delta),\mcal{R})) > 0$, and the directions $d_k^\delta$ in $\mcal{I}_4$ need not satisfy \eqref{LS descend prop.} for all $k$. Moreover, the advantage of defining a regularized solution as in $\mcal{I}_4$ over $\mcal{T}_2$ is that, one does not have to worry about the convergence of $\{x_k^\delta\}$ anymore, since we are not interested in finding the minimizer of a Variational problem \eqref{Gen. Tik. fun.}, rather, we would like to estimate the solution of the inverse problem \eqref{Inv. prob.} iteratively via a sequence of well-defined iterates $x_k^\delta$, depending on the directions $d_k^\delta$, where the best estimate $x_{k(\delta,\mcal{S})}^\delta$ to the true solution $\hat{x}$ is determined via a selection criterion $\mcal{S}$.

Observe that, although the initial iterate $x_0^\delta$ is not explicitly mentioned in the families of the regularized solutions $\mcal{I}_1, \mcal{I}_2, \mcal{I}_3$ and $\mcal{I}_4$, it is inherently embedded in the recovery process, unlike in the families $\mcal{T}_1$ and $\mcal{T}_2$, for convex regularizers $\mcal{R}$, where $x_0^\delta$ does not matter. To explicitly state the importance of the initial iterates in these (non-Variational or (semi-) iterative) families, $\mcal{I}_4$ can be expanded to 
\begin{align}\label{LS. family I_5}
\mcal{I}_5 := &\left\{ x^\delta(\mcal{D},k,d_k^\delta,x_0^\delta;\mcal{S}) \right. : \left. x_{k}^\delta = x_{k-1}^\delta \;+\; d_{k-1}^\delta, \; \mbox{ starting from $x_0^\delta$}, \right. \notag\\
& \left. \; \mbox{s.t. $d_{k-1}^\delta$ satisfies \eqref{descent dir. gen.} and $\mcal{S}$ is a selection criterion.}  \right\}.
\end{align}
The significance of the initial iterate $x_0^\delta$ can be prominently seen in the ADMM-PnP algorithm, which is discussed later. However, if the descent directions depends continuously, as a function, on its iterates, then for a small change in the initial iterates the difference in the recovered regularized solutions will be bounded by a factor similar to \eqref{LS z_k - x_k}, i.e., if $||d_k^\delta(x_k^\delta) - d_k^\delta(z_k^\delta))||$ is small for small $||x_k^\delta - z_k^\delta||$ and $||x_0^\delta - z_0^\delta|| \leq \epsilon$, then $||z_k^\delta - x_k^\delta|| \leq C(k)\epsilon$. This implies stability in the recovery process with respect to any small perturbations in the initial choices.
\subsection{Differential equation's solutions as a regularization family}
\mbox{ }

As one might have guessed by now, all the above families of solutions can be categorized as the solution states corresponding to a discrete approximation of some differential equations, in particular, the evolution equation given by

\begin{align}\label{evolution equation}
    \frac{\partial}{\partial t} x^\delta(.,t) &= f(x^\delta(.,t)) \\
    x^\delta(.,t_0) &= x_0^\delta, \notag
\end{align}
where the function $f(x^\delta(.,t))$ determines the flow or the evolution of a family of estimated solutions $x^\delta(.,t)$ to the inverse problem \eqref{Inv. prob.}, and an appropriate estimate $x^\delta(.,T(\delta))$ for $\hat{x}$ is determined via a selection criterion. When $f(x^\delta(.,t)) = -\tau(t)\nabla_t \mcal{D}(x^\delta(.,t)) $, for $\tau(t) \geq 0$, we have the classical steepest-descent method in the continuous settings. For $f(x^\delta(.,t))$ satisfying $\lprod{f(x^\delta(.,t))}{-\nabla_t \mcal{D}(x^\delta(.,t))}>0$, it falls in the family $\mcal{I}_3$, where as, for $f(x^\delta(.,t))$ satisfying \eqref{descent dir. gen.}, we end in the family $\mcal{I}_5$, in their respective continuous versions. 
Here, one can observe that a regularized solution is not only a function of $f(x^\delta(.,t))$ (the flow direction), but also, on $x_0^\delta$ (the initial/starting state) and $T(\delta)$ (the stopping/terminal time), i.e., $x^\delta = x^\delta(.;f,x_0^\delta, T(\delta))$. In contrast, the minimization of the (convex) Variational problem \eqref{Gen. Tik. fun.} is independent of the starting point and the descent direction, when minimized in an iterative manner, since a good starting point ($x_0^\delta$) and a faster descent direction $(f(x^\delta(.,t)))$ only improves the descent rate to reach the minimizer of the \eqref{Gen. Tik. fun.}, but not the final solution. It only depends on the definition of the penalized/constrained functional ($\mcal{R}$) and the penalty/regularization parameter ($\lambda$) in \eqref{Gen. Tik. fun.}, i.e., $x^\delta(.;\mcal{R},\lambda)$. Therefore, based on a fixed selection criterion $\mcal{S}$, comparing a Variational regularized solution $x^\delta(.;\mcal{R},\lambda)$ to a (semi-) iterative regularized solution $x^\delta(.;f,T,x_0^\delta)$ we see that the analogous of $\mcal{R}$ is $f$ and $\lambda$ is $T$. However, the (true) dependence of the initial point $x_0^\delta$ is still missing in the Variational solution (assuming one minimizes the convex functional completely and the iterates converge to the global minimizer) and has an importance in the (semi-) iterative method. Nevertheless, the recovery will be stable to a small perturbation in the initial condition, if $f(x^\delta(.,t))$ depends continuously on $x^\delta(.,t)$, i.e., one can estimate a bound for the difference of the solution flow at any time $t$, which is given by
\begin{equation}
||x^\delta(.,t) - z^\delta(.,t)|| \leq c(t)\int_0^{t} ||f(x^\delta(.,s)) - f(z^\delta(.,s))|| ds \; + \; ||x_0^\delta  - z_0^\delta||,
\end{equation}
for some constant $c(t)$, which implies
\begin{equation}
||x^\delta(.,t) - z^\delta(.,t)|| \leq c(t)\int_0^{t} \epsilon^\delta(s) ds + ||x_0^\delta  - z_0^\delta|| \; \leq \; C(t) \epsilon,
\end{equation}
if $||f(x^\delta(.,s)) - f(z^\delta(.,s))|| \leq \epsilon^\delta(t)$, for $t_0 \leq s \leq t$, and  $||x_0^\delta - z_0^\delta|| \leq \epsilon $.

Note that, if $f(x^\delta(.,t))$ satisfies \eqref{LS descend prop.} or \eqref{descent dir. gen.} then it is implicitly dependent on the discrepancy term $\mcal{D}(x)$ through the inner product with its gradient $-\nabla_x \mcal{D}(x_k^\delta)$, i.e., $f(x^\delta(.,t)) = f(x^\delta(.,t);\mcal{D})$. Hence, to further generalize the class $\mcal{I}_5$, one can even remove any dependence (implicit or explicit) on the discrepancy term $\mcal{D}$, i.e., no knowledge of the forward operator $A$ is used. For example, one can generate the flow of solutions $x^\delta(.,t)$ driven solely by a collection of informative dataset, i.e., using machine learning or deep learning algorithms, for example, the AUTOMAP network in \cite{Zhu_Liu_Bruce_Rosen}. The advantage of this formulation is that, any short-comings or mismatch in the model $A$ can be compensated by replacing it with a data based model (say neural networks). However, this strength also becomes its weakness for certain (adversarial) noisy data  and pushes the recovery process further into the realm of instabilities (now being generated by the data model), see \cite{Antun_Renna_Poon_Adcock_Hansen}, if no ``proper regularization" is employed, where the proper regularization implies regularizing the recovery process rather than regularizing the neural network, this is a topic of discussion for another paper. Nevertheless, we end this section by defining the largest family of regularized solution to the inverse problem \eqref{Inv. prob.} as

\begin{align}\label{LS. family I_6}
\mcal{I}_6 := &\left\{ x^\delta(f,x_0^\delta,t;\mcal{S}) \right. \; : \; x^\delta(f,x_0^\delta,t) \; \mbox{ is a solution state of \eqref{evolution equation} at time t,}  \\
& \left. \mbox{s.t. $f(x^\delta(.,t)) \xrightarrow{\delta \rightarrow 0} -\tau(t)\nabla_x \mcal{D}(x(.,t))$ and $\mcal{S}$ is a selection criterion.}  \right\}. \notag
\end{align}

\subsection{Alternating direction method of multipliers (ADMM) - PnP}\label{Sec. ADMM-PnP}
\mbox{ }

Now we have developed all the theories needed to attach a meaning to the ADMM-PnP algorithms. Note that, a typical step in the ADMM-PnP algorithm is given by
    \begin{align}
        x_{k+1}^\delta &= \text{Prox}_{\frac{1}{\rho}\mcal{D}}(z_k^\delta - u_k^\delta) \;\; \longleftarrow \text{ data-consistency step}\label{ADMM-PnP data-cons.} \\
        &= \argmin_x \;\; \mcal{D}(x) + \rho ||x - (z_k^\delta - u_k^\delta)||_2^2 \notag\\
        z_{k+1}^\delta &= H_{\sigma_k}(x_{k+1}^\delta + u_{k}^\delta) \;\; \longleftarrow \text{ data-denoising step} \label{ADMM-PnP data-deno.}\\
        u_{k+1}^\delta &= u_k^\delta + x_{k+1}^\delta - z_{k+1}^\delta \;\; \longleftarrow \text{ noise update step} 
        \label{ADMM-PnP noise upd.}
    \end{align}
where if $H_{\sigma_k}$, for $\sigma_k = \sigma/\rho$, corresponds to the proximal map of a proper, closed and convex regularizer $\mcal{R}(x)$, then the Lagrangian parameter ($\rho>0$) only effects the speed of the convergence and not the solution of \eqref{Gen. Tik. fun.}, i.e., the iterates $x_k^\delta$ and $z_k^\delta$ will converge to the minimzer of the functional \eqref{Gen. Tik. fun.}, irrespective of the $\rho$-value. However, this can not be guaranteed for any general denoiser $H_{\sigma_k}$, in fact, there might not even exists an associated proximal map corresponding to any general denoiser, and hence, a regularization function $\mcal{R}(x)$ in \eqref{Gen. Tik. fun.}. In this section, we connect ADMM-PnP algorithms to the theories developed in the previous sections and also point out the significance of the Lagrangian parameter $\rho$ in this setting. Note that, at any step $k$, the minimizer ($x_{k+1}^\delta$) in the data-consistency step \eqref{ADMM-PnP data-cons.}  has a closed form solution, for $\mcal{D}(x) = ||Ax-b_\delta||_2^2$, given by
\begin{align}\label{ADMM-PnP x_k}
x_k^\delta = (A^*A + \rho I)^{-1}(A^*b_\delta + \rho (z_{k-1}^\delta - u_{k-1}^\delta)) \;\; \longleftarrow \; \mbox{Tikhonov-$\ell_2$ solution with $\lambda = \rho$}.
\end{align}
Therefore, the resulting direction from $x_k^\delta$ to $x_{k+1}^\delta$, for $k \geq 1$, is given by
\begin{equation}
d_k^\delta := x_{k+1}^\delta - x_k^\delta = (A^*A + \rho I)^{-1}\rho \left[ (z_k^\delta - z_{k-1}^\delta) - (u_{k}^\delta - u_{k-1}^\delta) \right],
\end{equation}
and it's a descent direction if it satisfy \eqref{LS descend prop.}, i.e.,
\begin{equation}
\lprod{d_k^\delta}{-\nabla_x\mcal{D}(x_k)} = \lprod{(A^*A + \rho I)^{-1}\rho \left[ (z_k^\delta - z_{k-1}^\delta) - (u_{k}^\delta - u_{k-1}^\delta) \right]}{-A^*(Ax_k^\delta - b_\delta)} > 0
\end{equation}
which, though is trivial for $k=0$ with $x_0^\delta = z_0^\delta = u_0^\delta \equiv 0$, as
\begin{equation}
\lprod{d_0^\delta}{-\nabla_x\mcal{D}(x_0)} = \lprod{(A^*A + \rho I)^{-1}A^*b_\delta}{A^*b_\delta} > 0,
\end{equation}
is non-trivial for $k \geq 1$. However, this is where we can either use the regularization family $\mcal{I}_5$, as $d_0^\delta$ satisfies \eqref{LS descend prop.}, or $\mcal{I}_6$, where $d_k^\delta$ doesn't need to satisfy any conditions, to attach a meaning to the iterative process of generating a regularized family of structured-solutions $\{x_k^\delta \}$, where an appropriate estimate ($x_{k(\delta)}^\delta$) to the true solution ($\hat{x}$) solution is determined by some selection criteria. Note that, for $\epsilon_1(\delta) \leq \epsilon_2(\delta)$, there is a $\rho(\delta)$ such that $\epsilon_1(\delta) \leq \mcal{D}(x_1^\delta(\rho(\delta))) \leq \epsilon_2(\delta)$, where $x_1^\delta(\rho(\delta)) = (A^*A+\rho(\delta)I)^{-1}A^*b_\delta$, and for $\delta \rightarrow 0$, we need $\rho(\delta) \rightarrow0$, since $\epsilon_2(\delta) \xrightarrow{\delta \rightarrow 0} 0$.

As mentioned earlier, the ADDM-minimization corresponding to a traditional regularizer $\mcal{R}(x)$ is independent of the Lagrangian parameter $\rho$, as the iterates converge. But, for any general denoiser $H_{\sigma}$ there is no proof of convergence and hence, the choice of $\rho$ will affect the recovered regularized solution. In addition, the initial choices $(x_0^\delta,z_0^\delta,u_0^\delta)$ also influence the recovery process in this case, if we consider the iterative process to be in the family $\mcal{I}_5$ or $\mcal{I}_6$ of regularized solutions. To show the nature of dependence on $\rho$, we first transform the minimization in \eqref{ADMM-PnP data-cons.} to an equivalent constraint minimization problem, given by
\begin{align}\label{ADMM x_k const.}
x_{k+1}^\delta :=
\begin{cases}
 \argmin_x \;\; &\mcal{D}(x) = ||Ax - b_\delta||_2^2, \\
 \mbox{subject to } \; &||x - (z_k^\delta - u_k^\delta)||_2^2 \; \leq \; \epsilon(\rho),
\end{cases}
\end{align}
for some associated $\epsilon(\rho) \sim 1/\rho$. That is, for larger values of $\rho$, $x_{k+1}^\delta$ is closer to $z_k^\delta - u_k^\delta$, and for smaller values of $\rho$, $x_{k+1}^\delta$ will be farther away from $z_k^\delta - u_k^\delta$, and closer to the noise corrupted solution $(A^*A + \rho I)^{-1}(A^*b_\delta + \rho (z_{k}^\delta - u_{k}^\delta))$. Hence, we divide the analysis into two parts:

\subsubsection{\textbf{When Lagrangian parameter is large}}
Assuming we start from $x_0^\delta = z_0^\delta = u_0^\delta \equiv 0$, after the $1^{st}$ step, we have a slightly improved $x_1^\delta$ towards the LS-solution (without being much corrupted by the noise), and hence, even for a weaker denoiser $H_{\sigma_1}$, one can denoise it to get $z_1^\delta$. On the $2^{nd}$ step, $x_2^\delta$ will again be improved only slightly towards the noisy LS-solution, because of the stronger constraint in \eqref{ADMM x_k const.}, which can then be effectively denoised by a weaker denoiser and so on. Therefore, starting from smooth initial choices, one can recover well denoised iterates $\{ z_k^\delta \}$, even for weaker denoisers, see Example \ref{Ex. ADMM_1}. This process is very similar to FBS-PnP algorithm, though the iterates here are not explicitly generated by any descent directions, rather, through a constrained minimization at each step.

\subsubsection{\textbf{When Lagrangian parameter is small}} In contrast to the above recovery process, here, even when the initial choices $(x_k^\delta, z_k^\delta, u_k^\delta)$ are smooth (say 0s), after the $1^{st}$ iteration we will have $x_1^\delta$ close to the noisy LS-solution ($x^\dagger_\delta$), i.e., $x_1^\delta$ too noisy, and hence, we need a much stronger denoiser to clean the corruptions; which in turn can move $z_1^\delta$ far away from $x_1^\delta$, because of the weaker constraint in \eqref{ADMM x_k const.}. Thus, the denoisers $H_{\sigma_k}$ need to be really effective in cleaning the heavily noised iterates $x_k^\delta$ at every step, to produce well-denoised iterates $z_k^\delta$, only then we may hope to recover a clean regularized solution. However, for weaker to moderately strong denoisers, unlike the FBS-PnP algorithm, here one might recover a noisy estimate to the solution of \eqref{Inv. prob.}, see Example \ref{Ex. ADMM_1}. 

\begin{remark}
From the above two analysis, one can notice that when $\rho$ is large the recovery process is very similar to the FBS-PnP method and one can clean the noise in the iterates $x_k^\delta$ gradually towards the ``improved" iterates $z_k^\delta$, of course, the improvement depends on the denoiser $H_{\sigma_k}$. But, the rate of descent can be very slow, similar to ISTA, and hence, it will be preferable to opt for the Fast FBS-PnP methods. Where as, if $\rho$ is small then the rate of descent corresponding to the discrepancy term $\mcal{D}(x)$ is fast, however, this doesn't ensure that the rate of recovery for a denoised solution to be fast too, since, the faster $\mcal{D}$ decreases, the faster the noise increases in the iterates. Therefore, here, it greatly depends on the effectiveness of the denoisers to obtain a cleaner recovery.
\end{remark}

\begin{remark}
In addition to the issues mentioned in the above remark, note that, in the ADMM minimization process, step \eqref{ADMM-PnP data-cons.} involves an additional minimization problem at each iteration $k$. This further increases the computational time and complexity, for example, for large scale problem the matrix inversion $(A^*A + \rho I)^{-1}$ is impossible and one has to approximate the exact minimizer $x_k^\delta$ via an iterative scheme (such as conjugate-gradient) for a certain number of fixed iterations, and hence, there is a further degradation to the proper definition of ADMM-minimization. Therefore, unless one has some prior knowledge of an ``appropriate" Lagrangian parameter $\rho$ and ``appropriate" denoisers $H_{\sigma_k}$, which can balance the speed of $\mcal{D}$-descent and the cleaning of the noisy iterates, we prefer the fast FBS-PnP algorithms, where even inappropriate denoisers can be tweaked appropriately to produce better results, as it involves a gradual recovery process. However, for appropriately strong $H_{\sigma_k}$ and $\rho$-value, one can get better performance (empirically) in ADMM-PnP than in FBS-PnP. Nevertheless, considering both these algorithms as iterative regularization processes in either family $\mcal{I}_5$ or $\mcal{I}_6$, the better one is the one that satisfies the selection criterion the best.
\end{remark}

\begin{remark}
Note that, the first iterate of the ADMM process (starting from initial $x_0^\delta = z_0^\delta = u_0^\delta \equiv 0$) is an $\ell_2$-Tikhonov solution, i.e., $x_1^\delta = \argmin_x \; \mcal{D}(x) + \rho||x||_2^2$. Hence, one can generalize it by replacing the standard $\ell_2$-Tikhonov solution with any general $\ell_2$-Tikhonov solution, i.e.,
\begin{equation}\label{ADMM-PnP x_1}
x_1^\delta \; = \; \argmin_x \; \mcal{D}(x) \; + \; \rho||Lx||_2^2,
\end{equation}
where $L$ can be either any wavelet transformation or the gradient operator. This may serve as a better initial $x_1^\delta$ for the ADMM-PnP process, than the standard $\ell_2$-solution and can lead to an improved regularized solution.
\end{remark}

\section{\textbf{Numerical Results}}
In this section we provide computational results validating the theories developed in this paper. Note that, the goal here is not to show the effectiveness of PnP-algorithms, which has been empirically shown in numerous papers, rather, we try to present some of the caveats of these algorithms, when do they fail, and how to improve them. Hence, in most of the examples we won't focus on repeating the same experiment over and over with different parameter values (i.e., tuning the denoising parameter) to attain the best solution, with the prior knowledge of the true solution. On contrary, for a fixed selection criterion $\mcal{S}_0$ and a fixed denoiser $H_{\sigma}$, we present: (1) how can one attain the best solution based on $\mcal{S}_0$, during the iterative process, and (2) how to further improve the recovery for the same $H_{\sigma}$ and $\mcal{S}_0$, for example, when the denoiser is too strong or too weak. We also compare the Fast FBS-PnP algorithm with the ADMM-PnP algorithm, i.e., the pros and cons of both these algorithms and, when do they fall in the family $\mcal{I}_3$ or in $\mcal{I}_5$.

\begin{figure}[h!]
    \centering
    \begin{subfigure}{0.495\textwidth}
        \includegraphics[width=\textwidth]{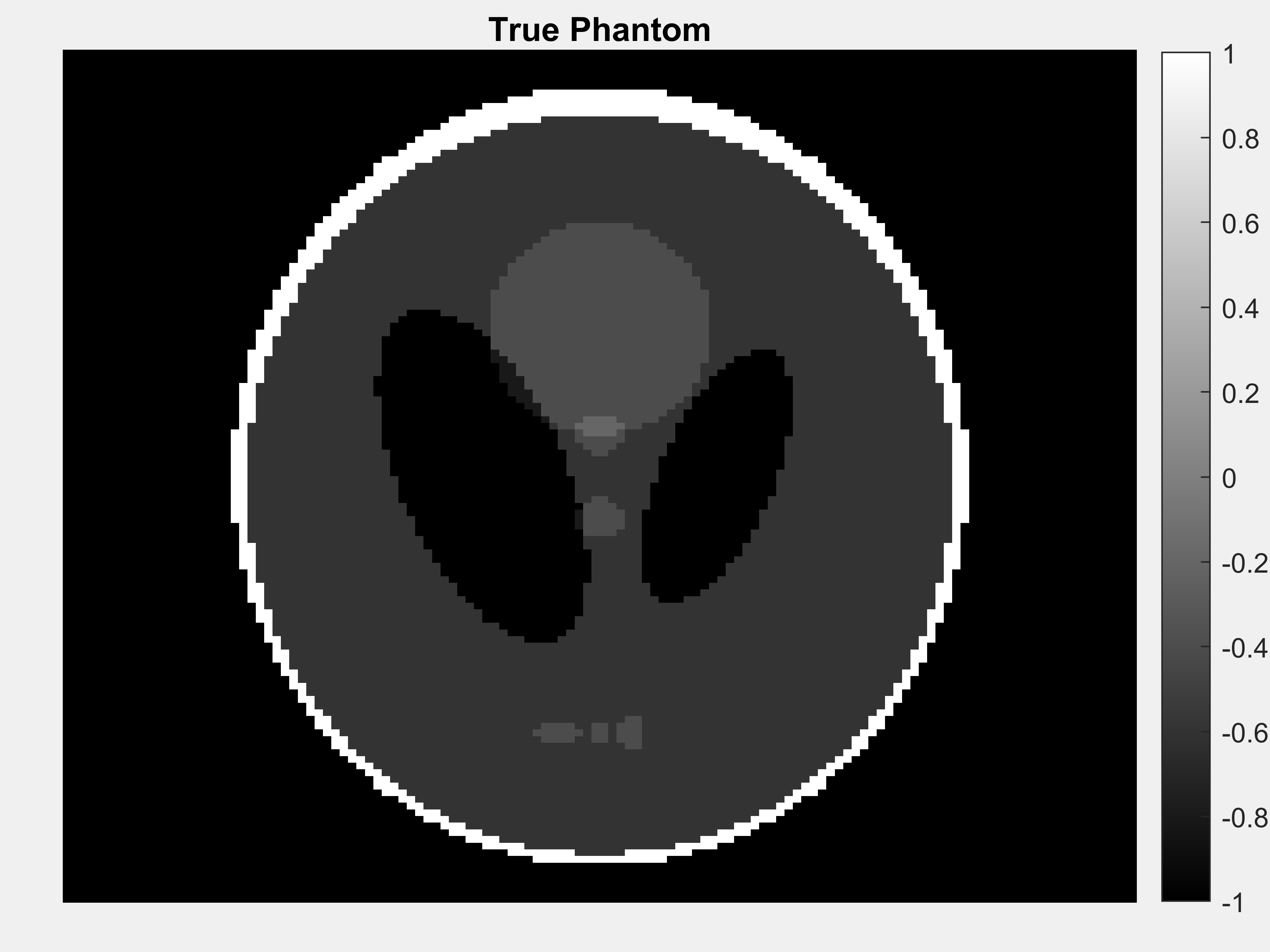}
        \caption{True Phantom.}
        \label{True Shepp-Logan Phantom}
    \end{subfigure}
    \begin{subfigure}{0.495\textwidth}
        \includegraphics[width=\textwidth]{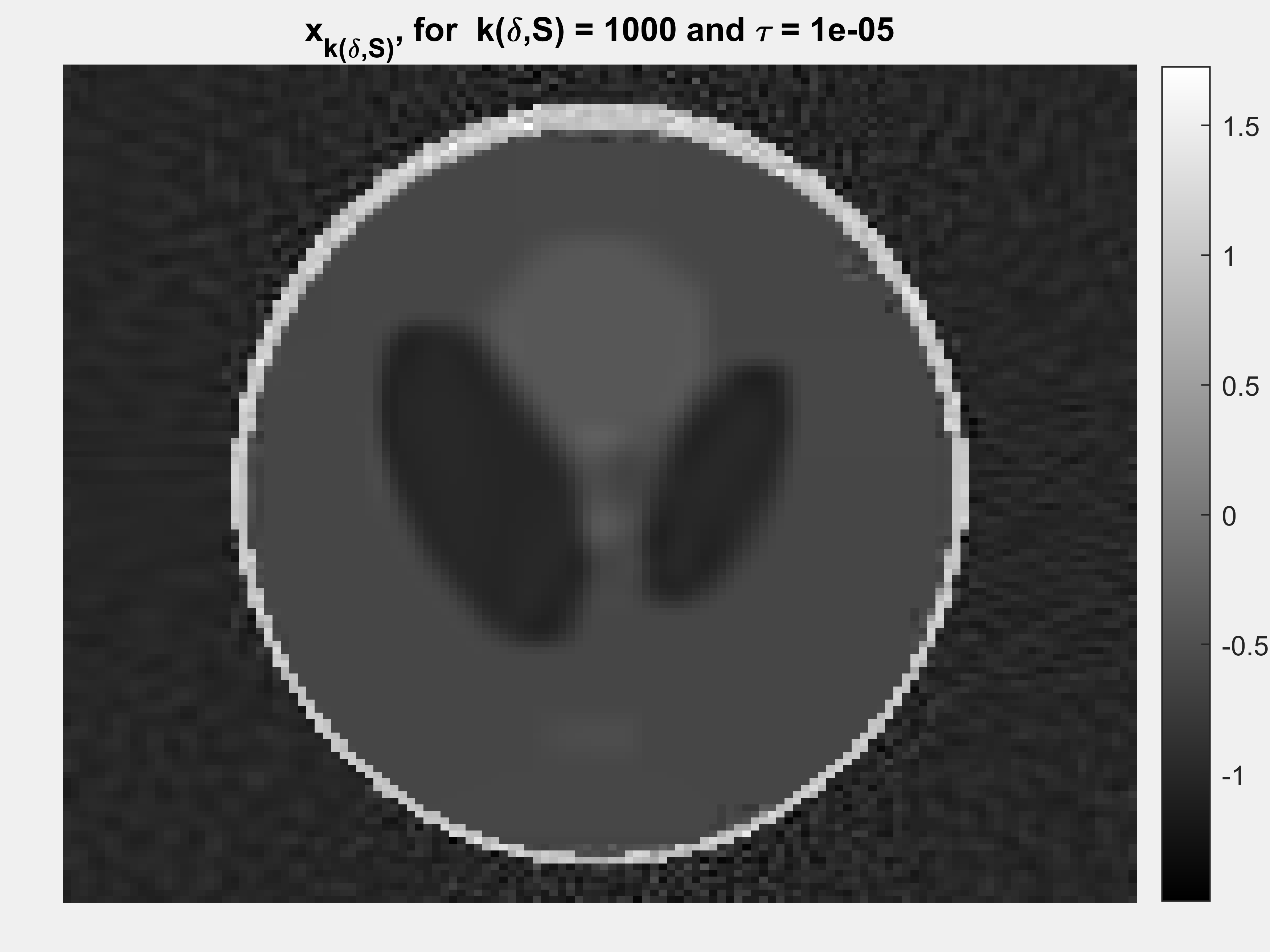}
        \caption{$z_{k(\delta,\mcal{S})}^\delta$, $k(\delta,\mcal{S})=1000$ and $\tau=10^{-5}$.}
        \label{x_k_e05}
    \end{subfigure}    
    \begin{subfigure}{0.495\textwidth}
        \includegraphics[width=\textwidth]{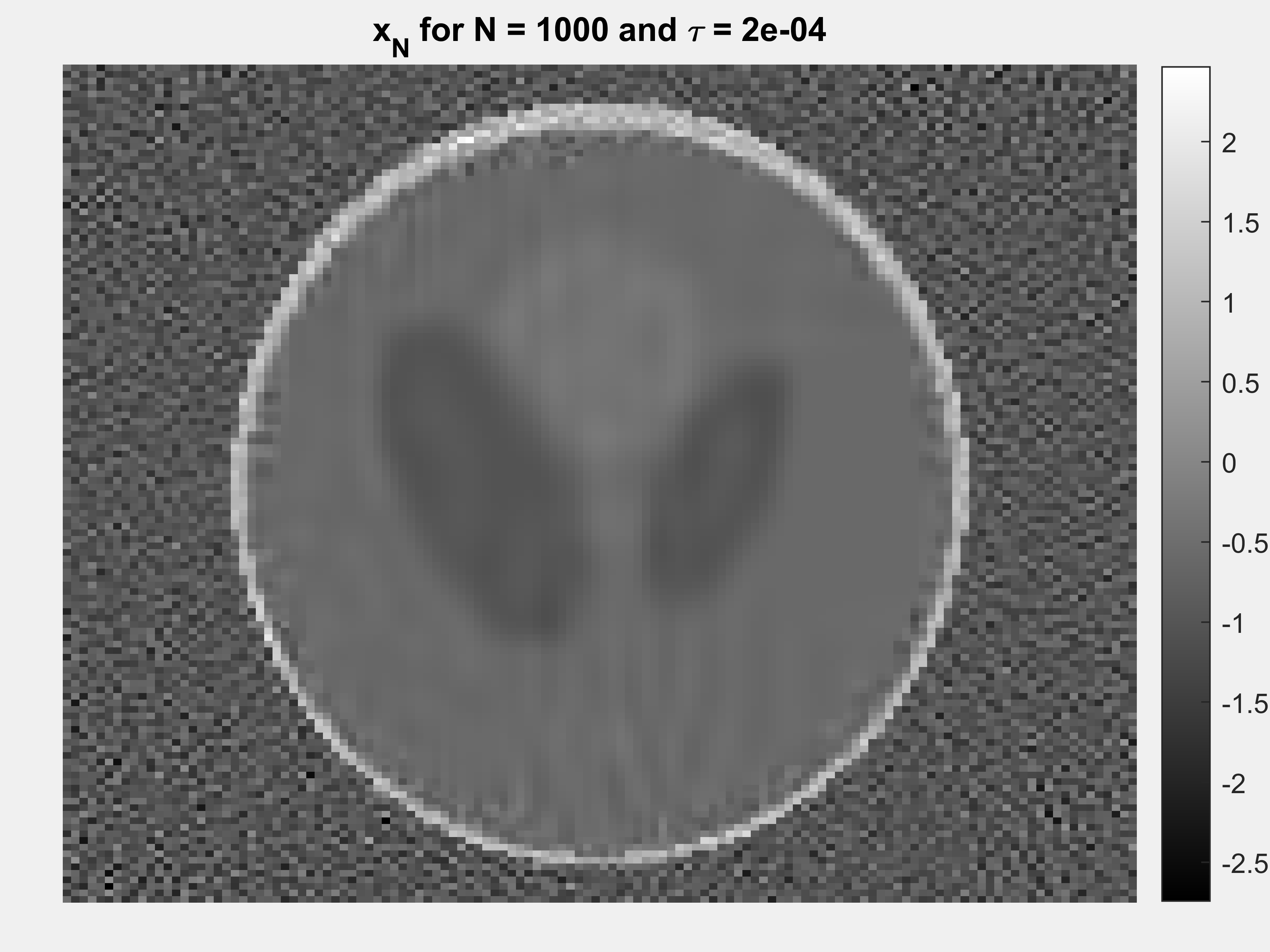}
        \caption{$z_N^\delta$, $N=1000$ and $\tau=2\times 10^{-4}$.}
        \label{x_N_2e04}
    \end{subfigure}
    \begin{subfigure}{0.495\textwidth}
        \includegraphics[width=\textwidth]{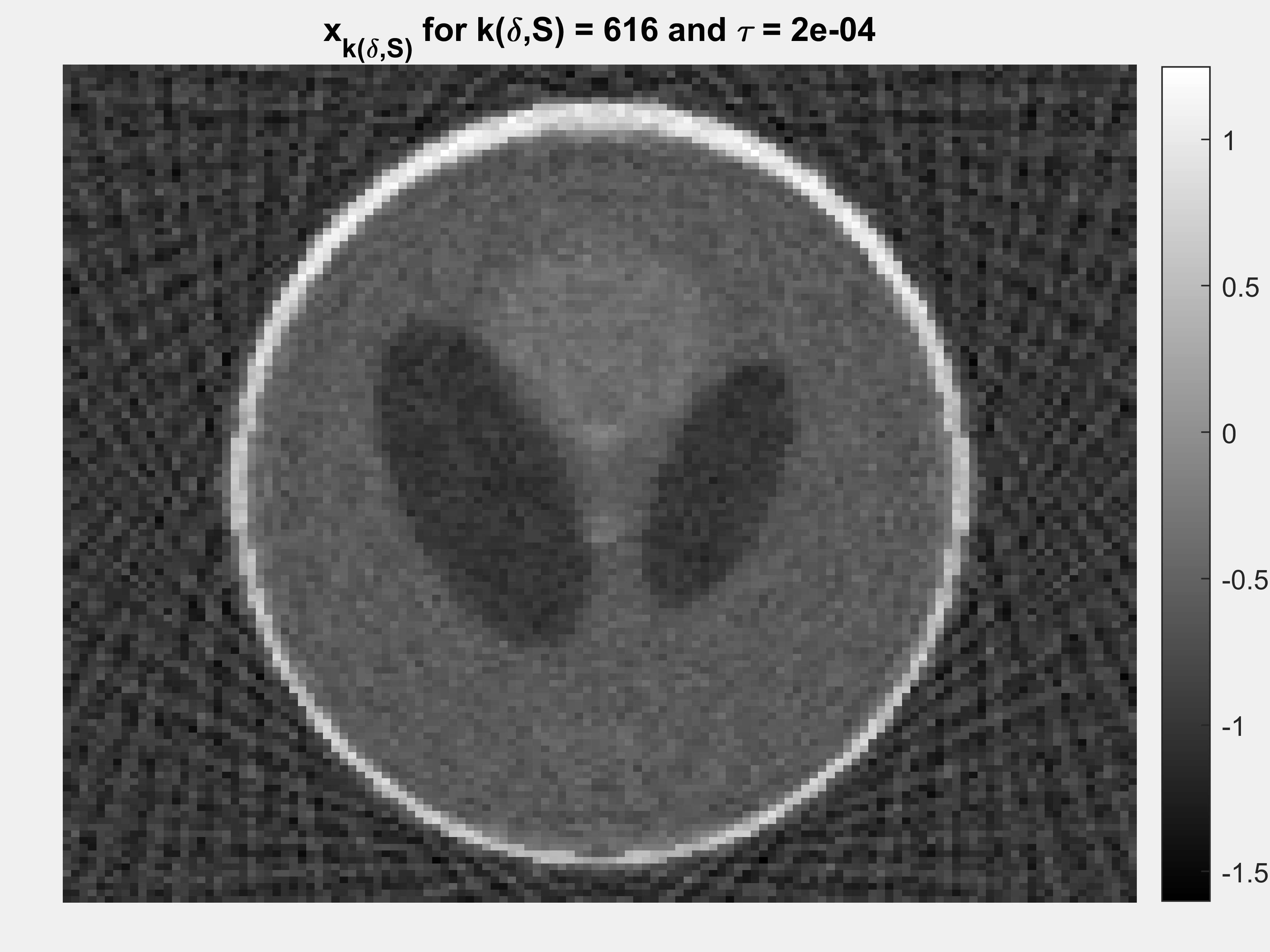}
        \caption{$z_{k(\delta,\mcal{S})}^\delta$, $k(\delta,\mcal{S})=55$ and $\tau=2\times 10^{-4}$.}
        \label{x_k_2e04}
    \end{subfigure}
    \begin{subfigure}{0.495\textwidth}
        \includegraphics[width=\textwidth]{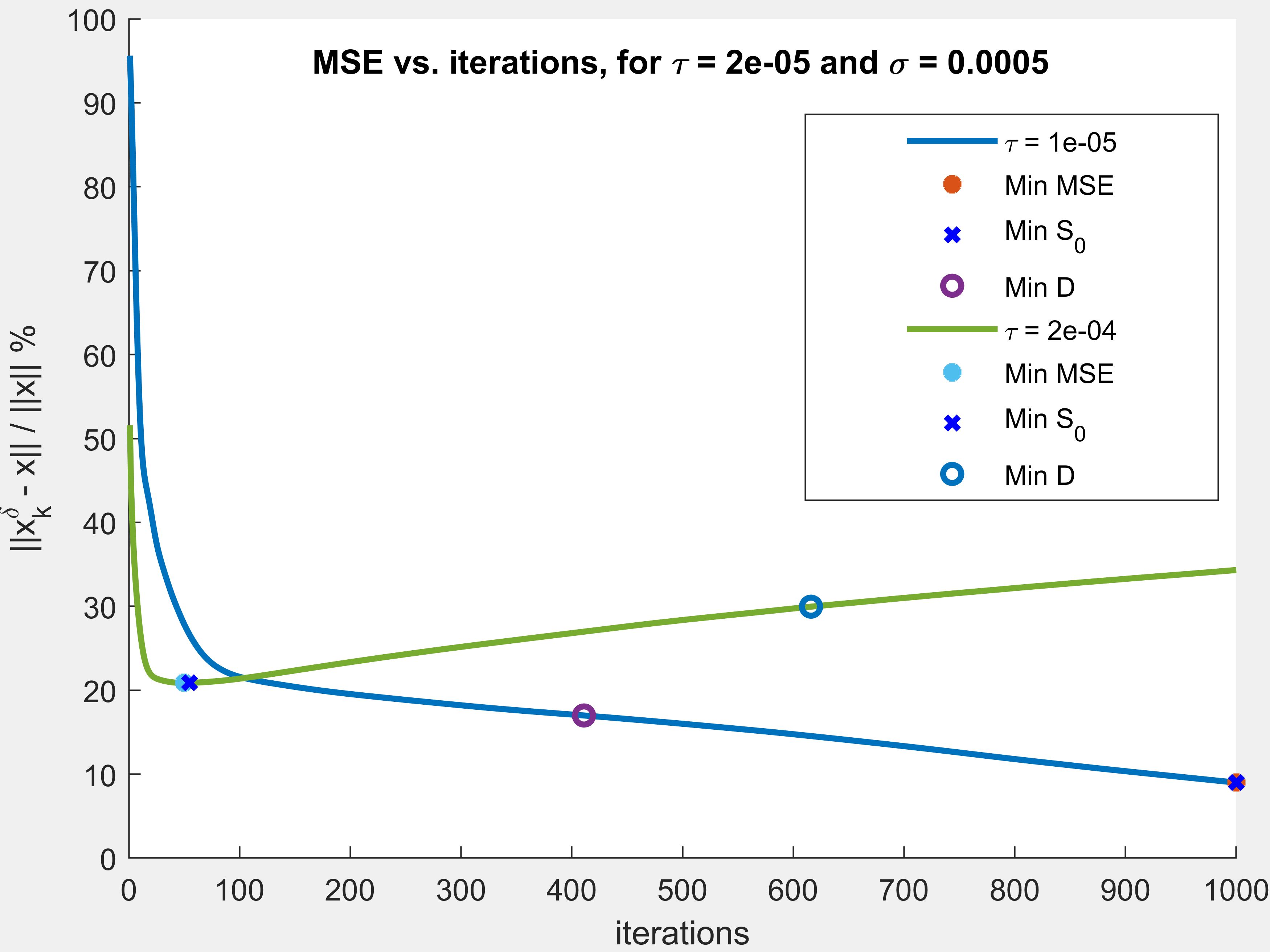}
        \caption{MSE vs. $k$, for $\tau=2\times 10^{-4}$ and $10^{-5}$.}
        \label{MSE_12}
    \end{subfigure}
    \begin{subfigure}{0.495\textwidth}
        \includegraphics[width=\textwidth]{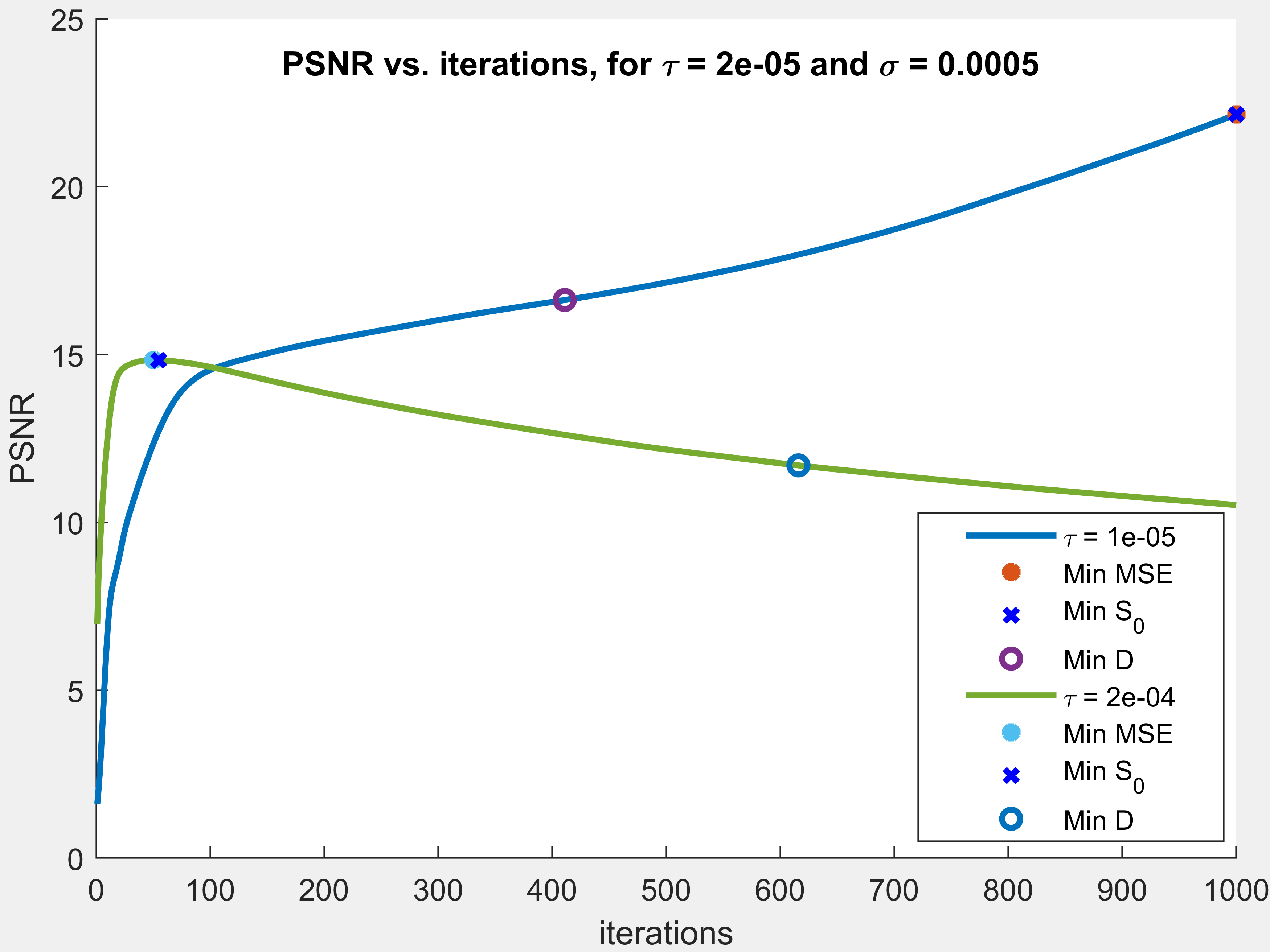}
        \caption{PSNR vs. $k$, for $\tau=2\times 10^{-4}$ and $10^{-5}$.}
        \label{PSNR_12}
    \end{subfigure}
    \caption{$z_N^\delta$ vs. $z_{k,\mcal{S}}^\delta$ and performance curves for Example \ref{Ex. weak deno.}.} 
    \label{Figure Ex. Weak deno.}
\end{figure}

\begin{figure}[h!]
    \centering
    \begin{subfigure}{0.495\textwidth}
        \includegraphics[width=\textwidth]{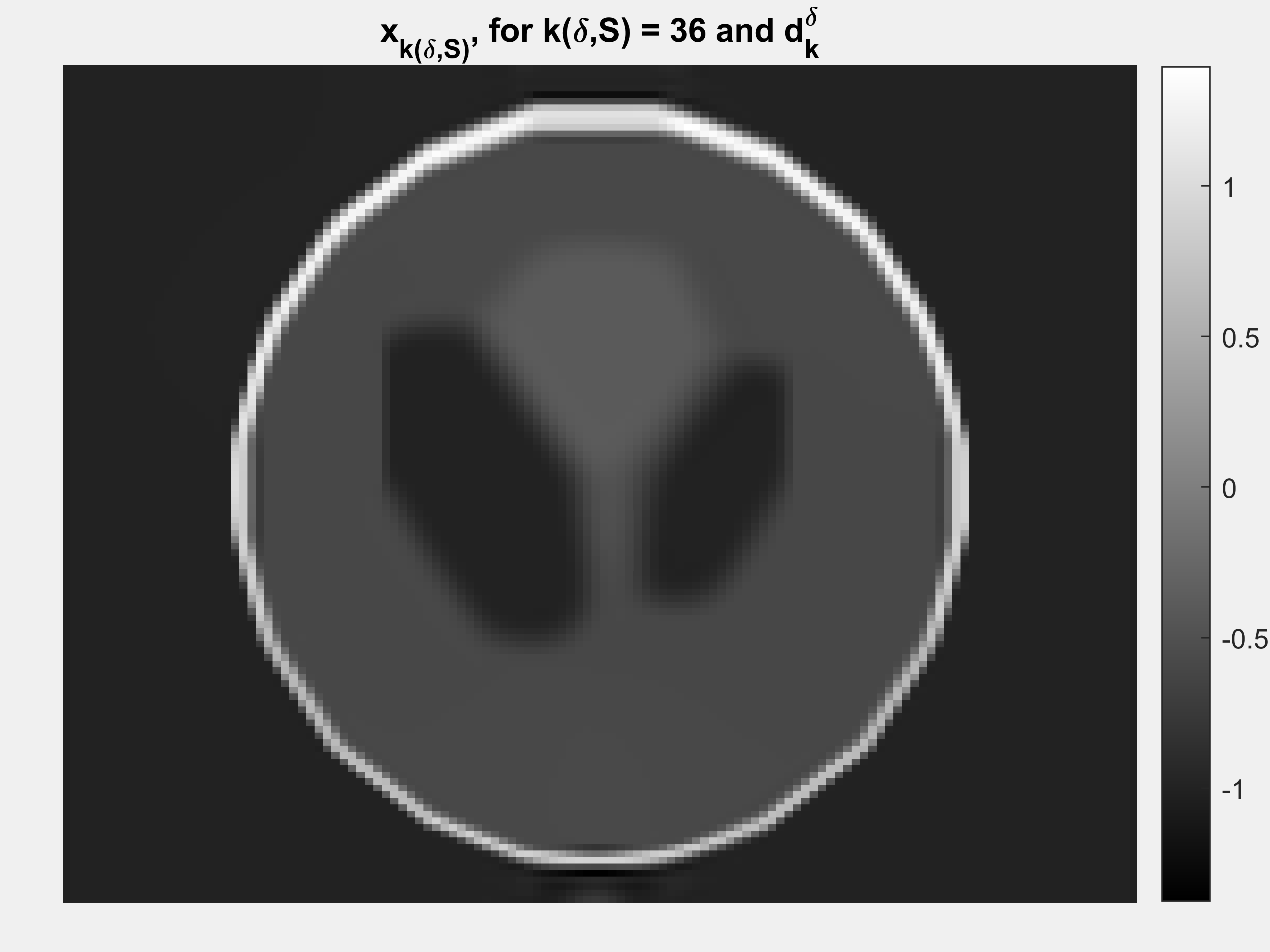}
        \caption{$z_{k(\delta,\mcal{S})}$, for $k(\delta,\mcal{S})=36$ and $d_k^\delta$}
        \label{noauto_xk}
    \end{subfigure}
    \begin{subfigure}{0.495\textwidth}
        \includegraphics[width=\textwidth]{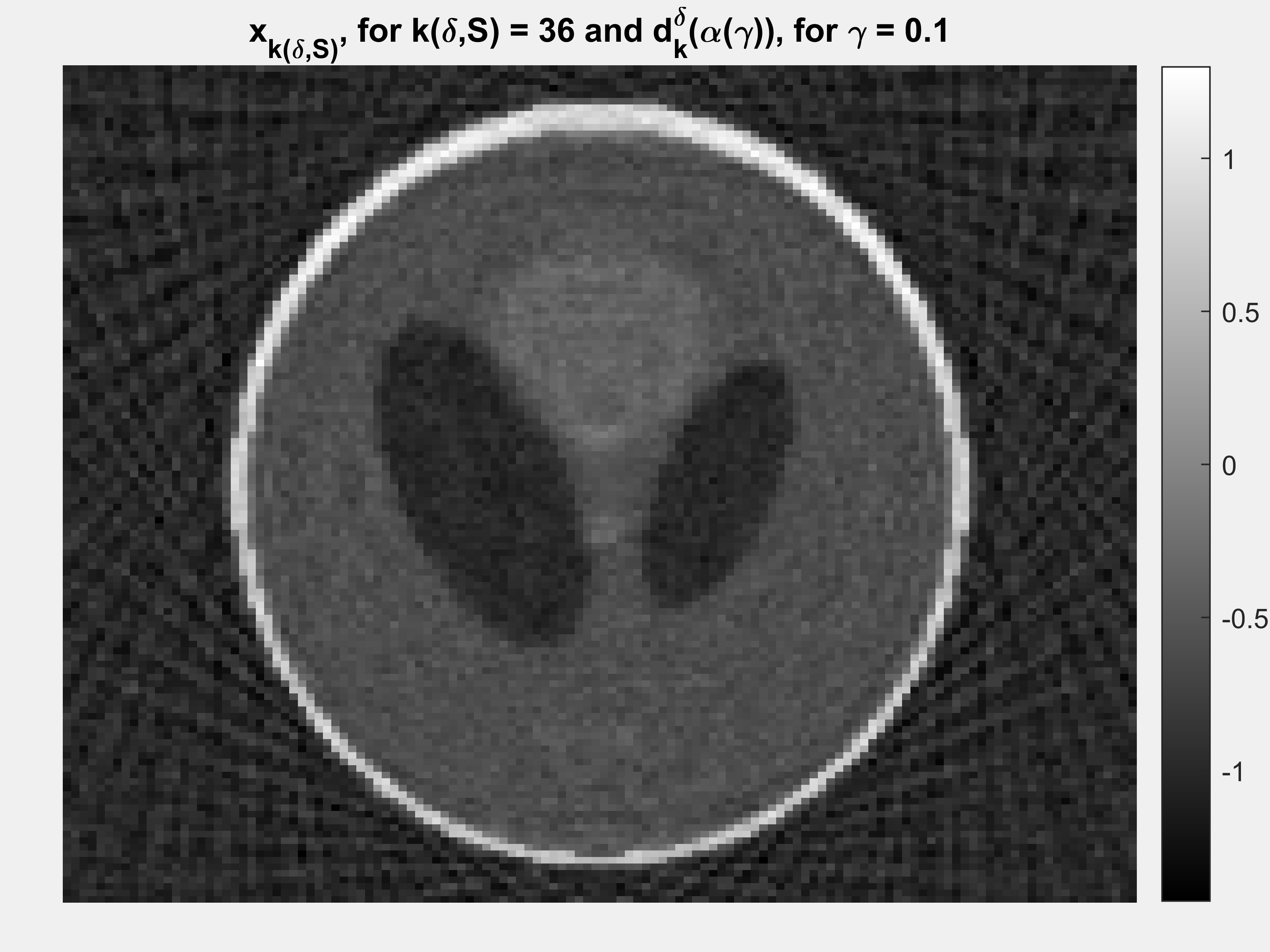}
        \caption{$z_{k(\delta,\mcal{S})}$, for $k(\delta,\mcal{S})=36$ \& $d_k^\delta(\alpha(\gamma = 0.1))$}
        \label{auto3_delpt1_xk}
    \end{subfigure}    
    \begin{subfigure}{0.495\textwidth}
        \includegraphics[width=\textwidth]{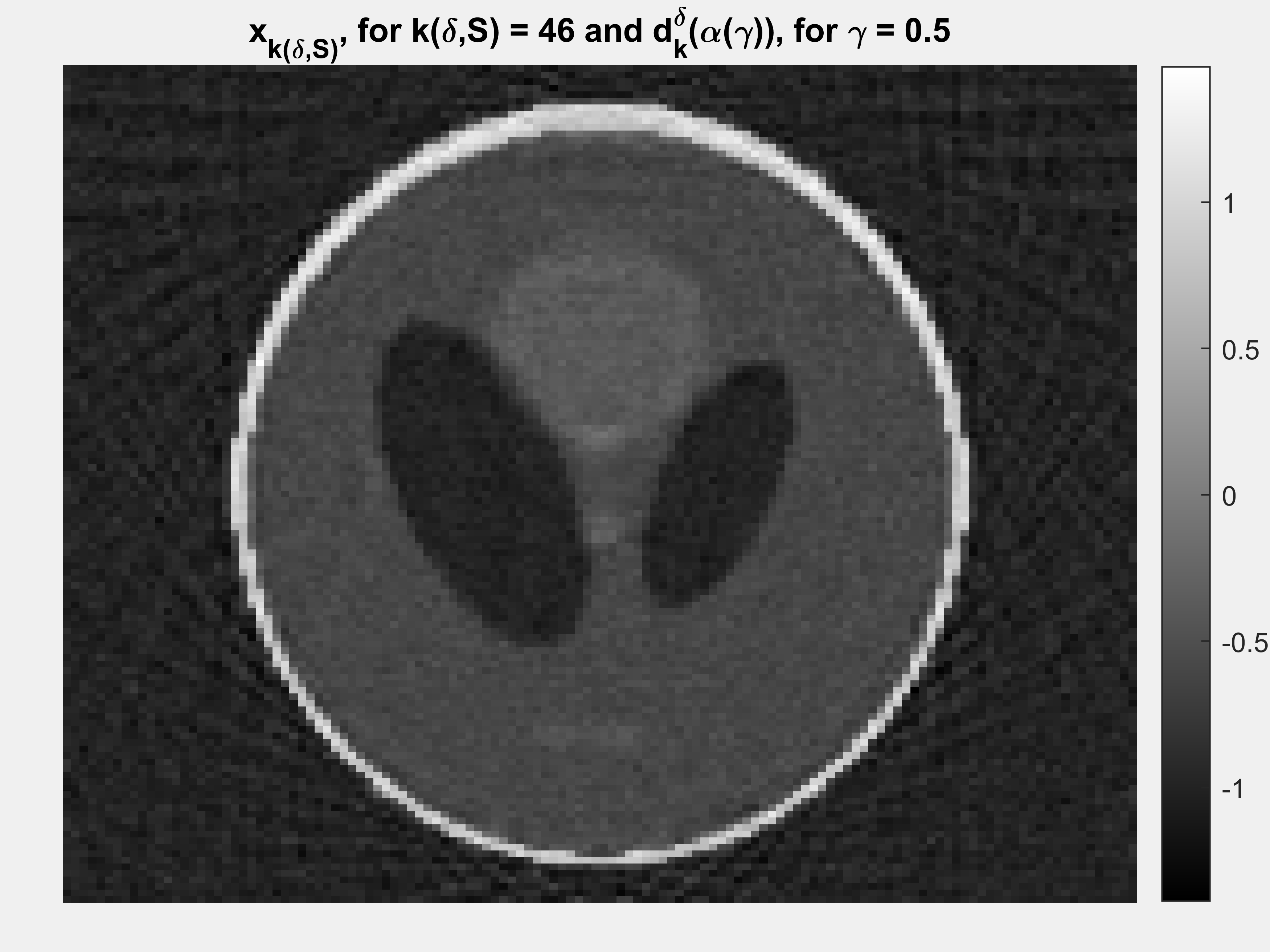}
        \caption{$z_{k(\delta,\mcal{S})}$, for $k(\delta,\mcal{S})=36$ \& $d_k^\delta(\alpha(\gamma = 0.5))$}
        \label{auto3_delpt5_xk}
    \end{subfigure}
    \begin{subfigure}{0.495\textwidth}
        \includegraphics[width=\textwidth]{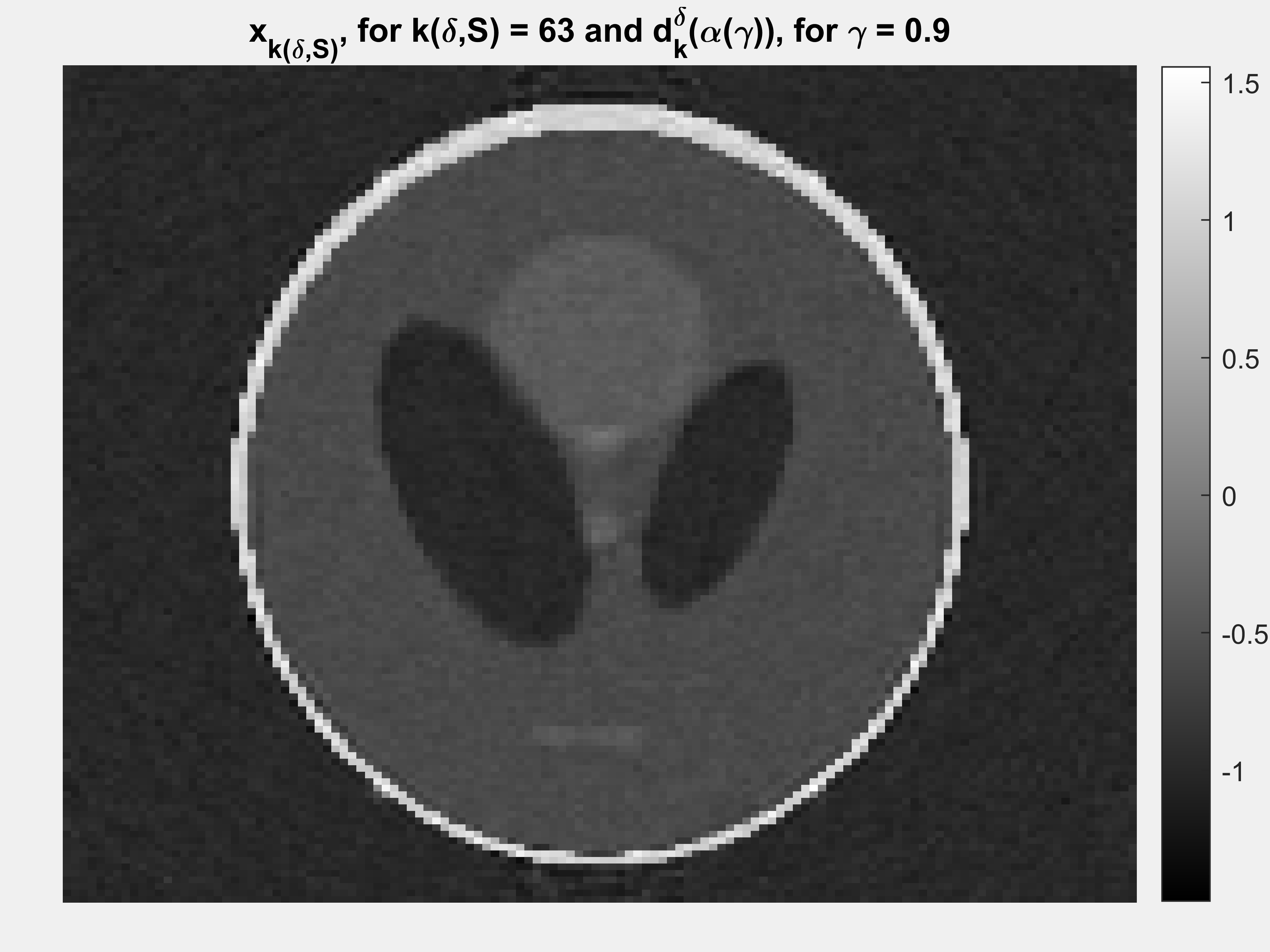}
        \caption{$z_{k(\delta,\mcal{S})}$, for $k(\delta,\mcal{S})=36$ \& $d_k^\delta(\alpha(\gamma = 0.9))$}
        \label{auto3_delpt9_xk}
    \end{subfigure}
    \begin{subfigure}{0.495\textwidth}
        \includegraphics[width=\textwidth]{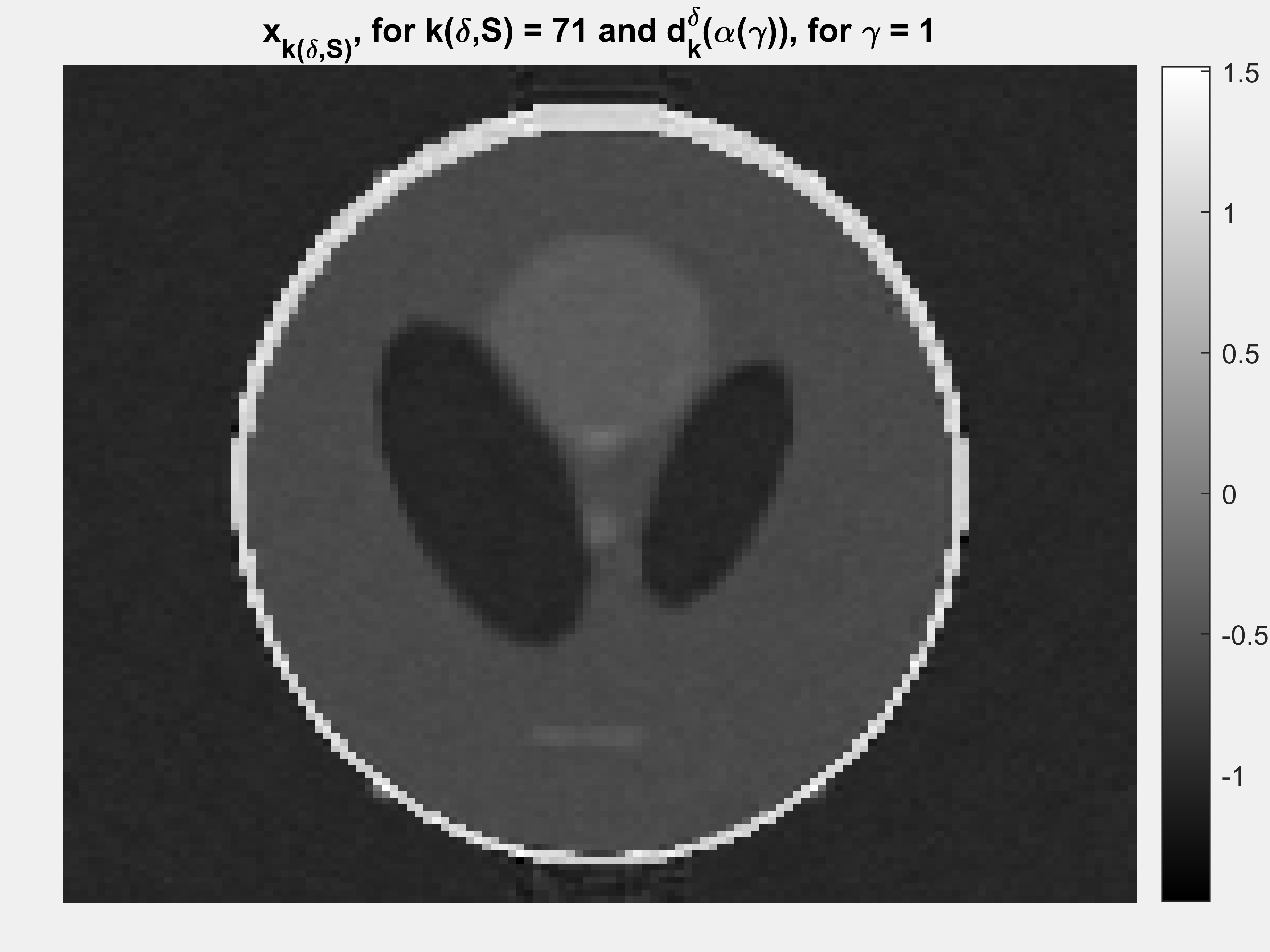}
        \caption{$z_{k(\delta,\mcal{S})}$, for $k(\delta,\mcal{S})=36$ \& $d_k^\delta(\alpha(\gamma = 1))$}
        \label{auto3_del01_xk}
    \end{subfigure}
    \begin{subfigure}{0.495\textwidth}
        \includegraphics[width=\textwidth]{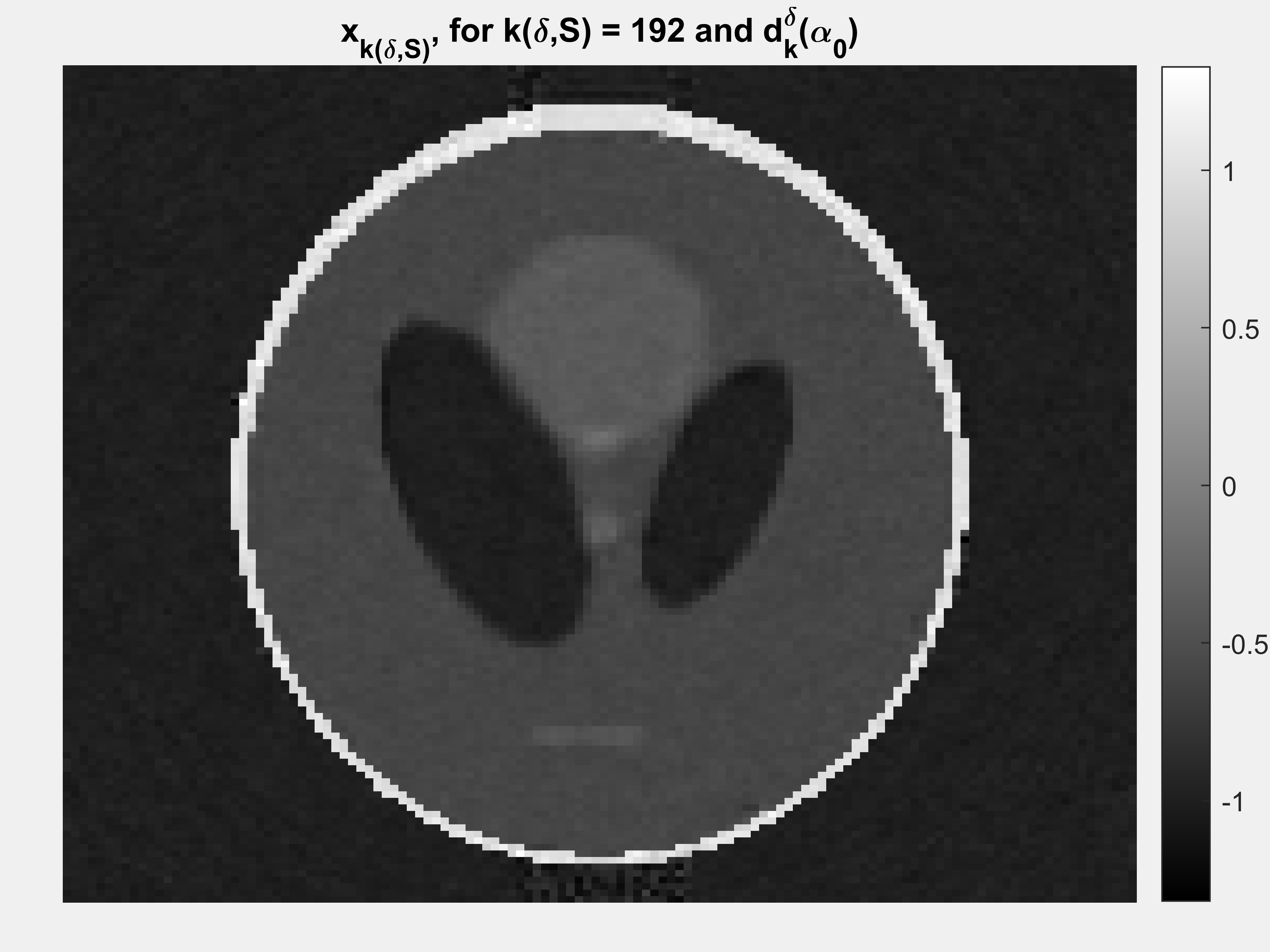}
        \caption{$z_{k(\delta,\mcal{S})}$, for $k(\delta,\mcal{S})=36$ \& $d_k^\delta(\alpha_0)$}
        \label{auto1_xk}
    \end{subfigure}
    \begin{subfigure}{0.495\textwidth}
        \includegraphics[width=\textwidth]{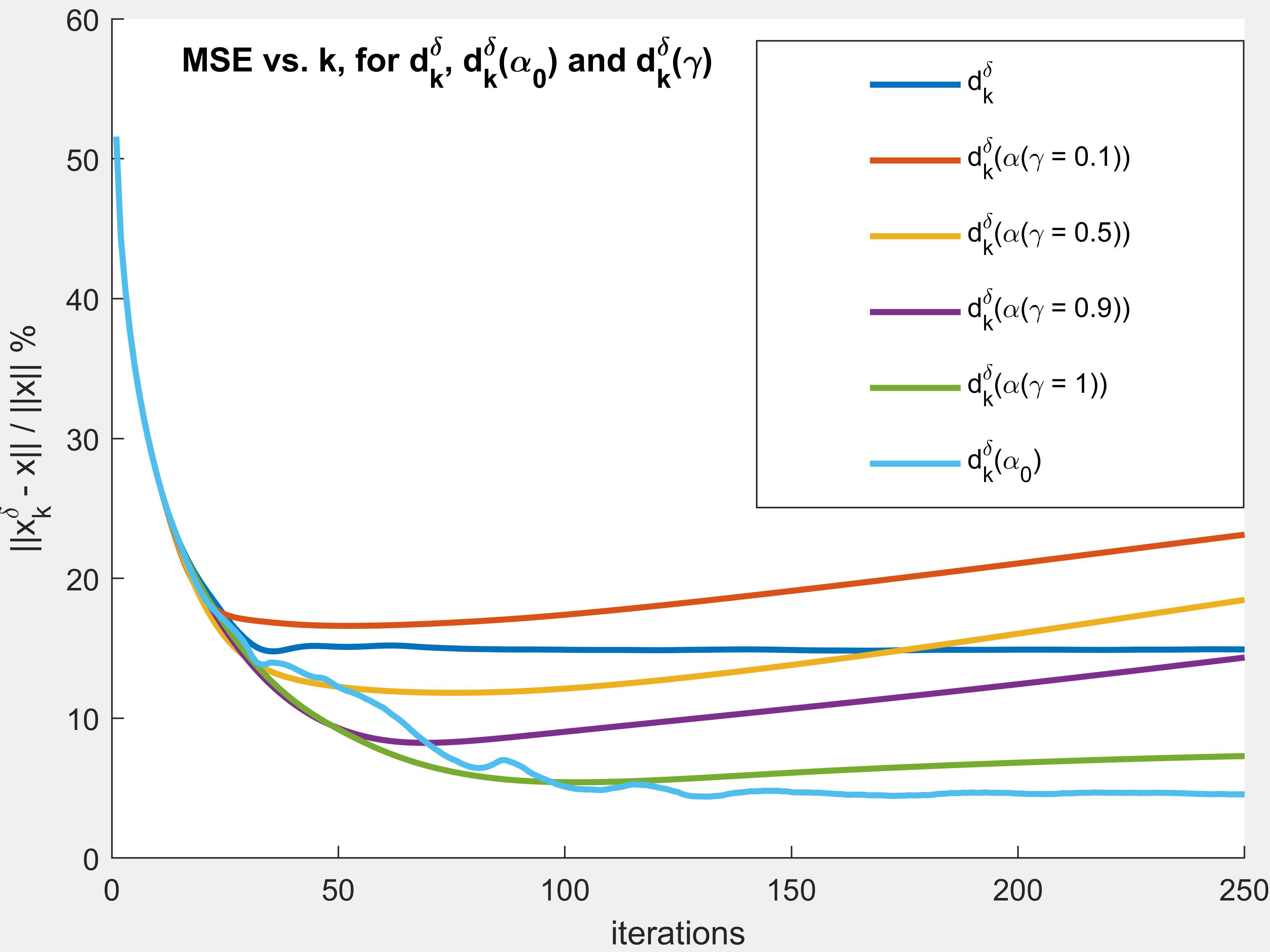}
        \caption{MSE curves for $d_k^\delta$, $d_k^\delta(\alpha_0)$ and $d_k^\delta(\gamma)$}
        \label{MSEcurves_StrongDeno}
    \end{subfigure}
    \begin{subfigure}{0.495\textwidth}
        \includegraphics[width=\textwidth]{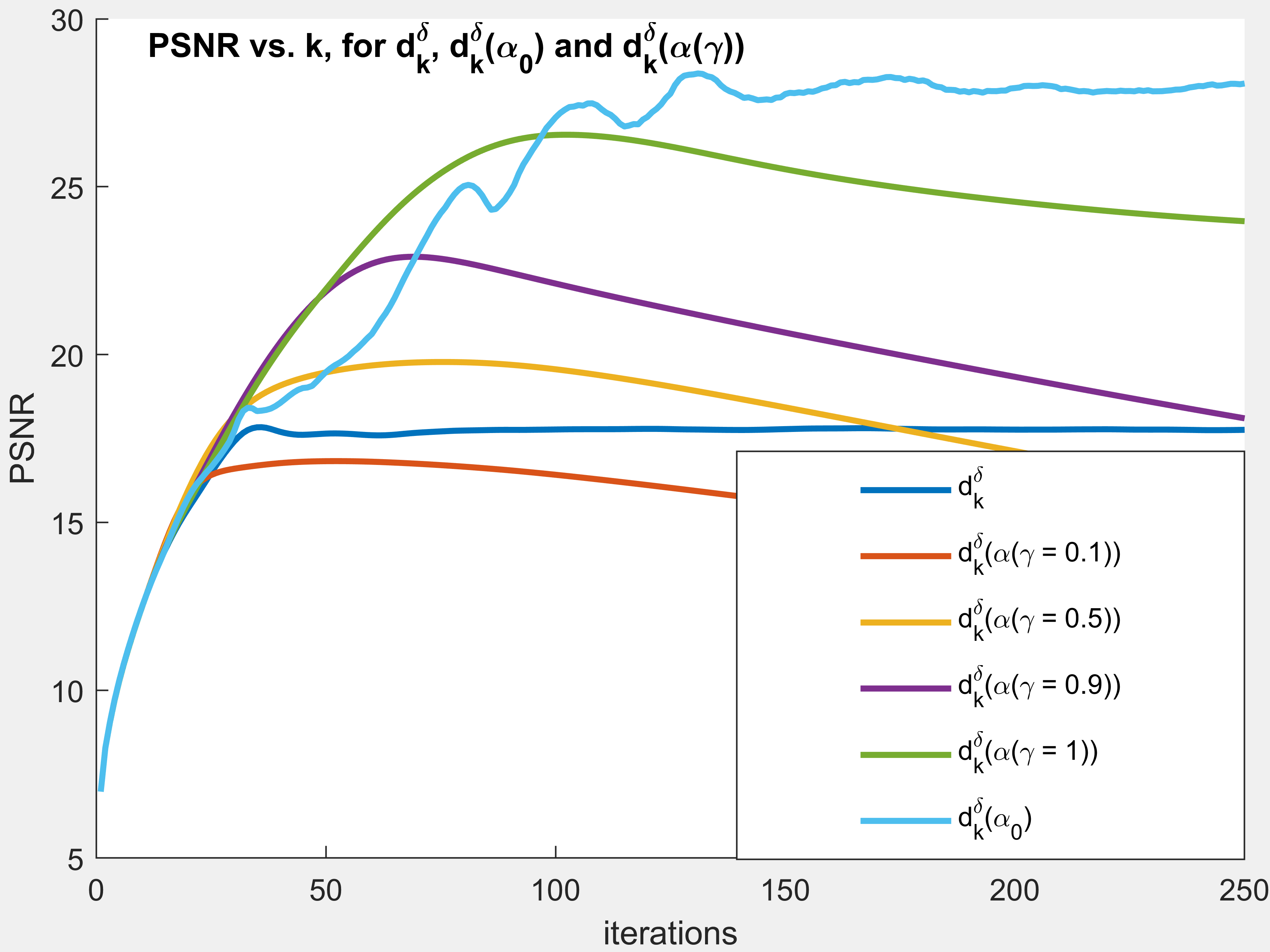}
        \caption{PSNR curves for $d_k^\delta$, $d_k^\delta(\alpha_0)$ and $d_k^\delta(\alpha(\gamma))$}
        \label{PSNRcurves_StrongDeno}
    \end{subfigure}    
    \caption{CV-recoveries and performance curves for Example \ref{Ex. Strong deno.}.} 
    \label{Figure Ex. Strong deno.}
\end{figure}

\begin{figure}[h!]
    \centering
    \begin{subfigure}{0.495\textwidth}
        \includegraphics[width=\textwidth]{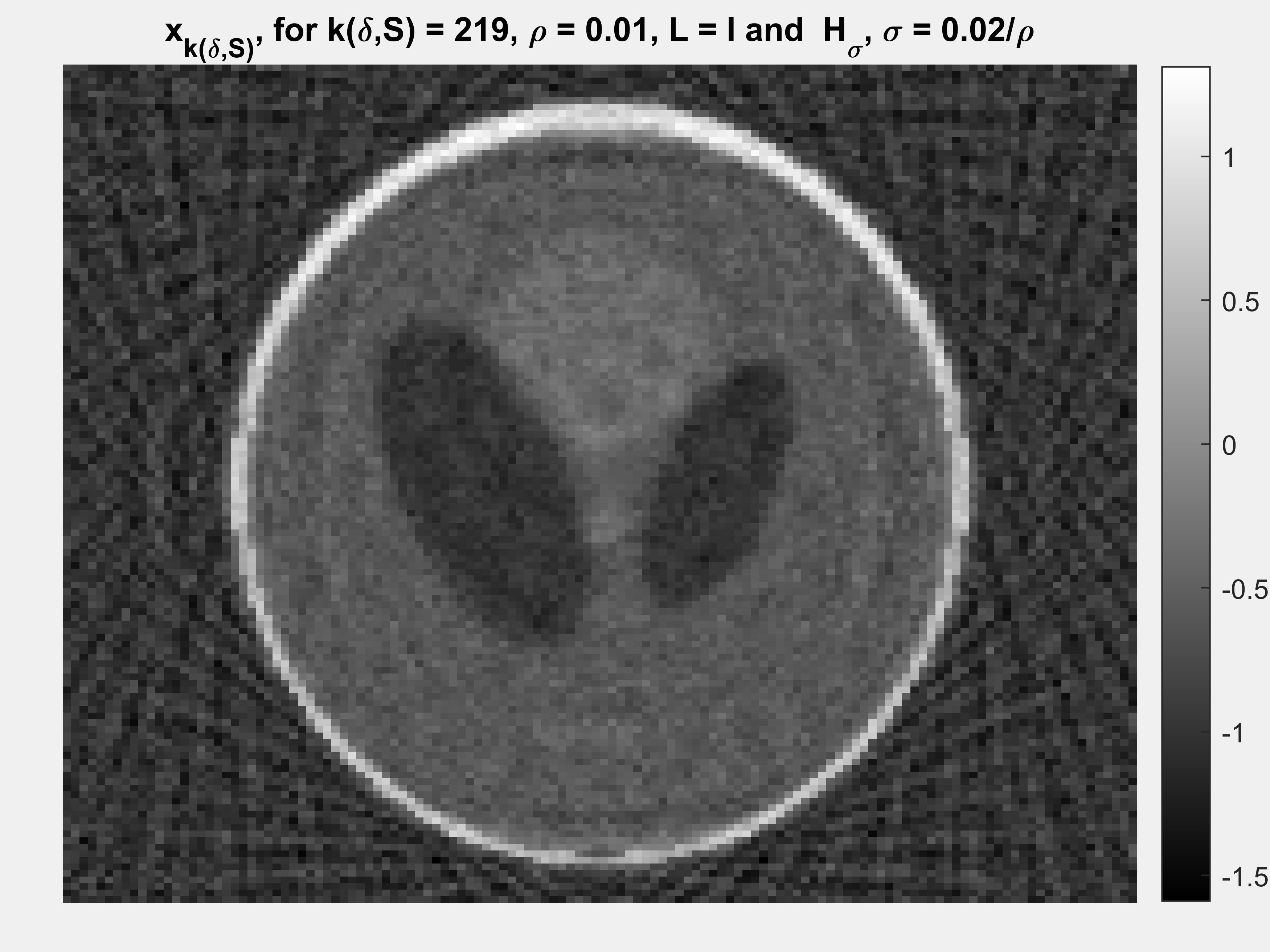}
        \caption{$x_{k(\delta)}^\delta$ for $H_{\hat{\sigma}}$ and $\rho = 0.01$}
    \end{subfigure}
    \begin{subfigure}{0.495\textwidth}
        \includegraphics[width=\textwidth]{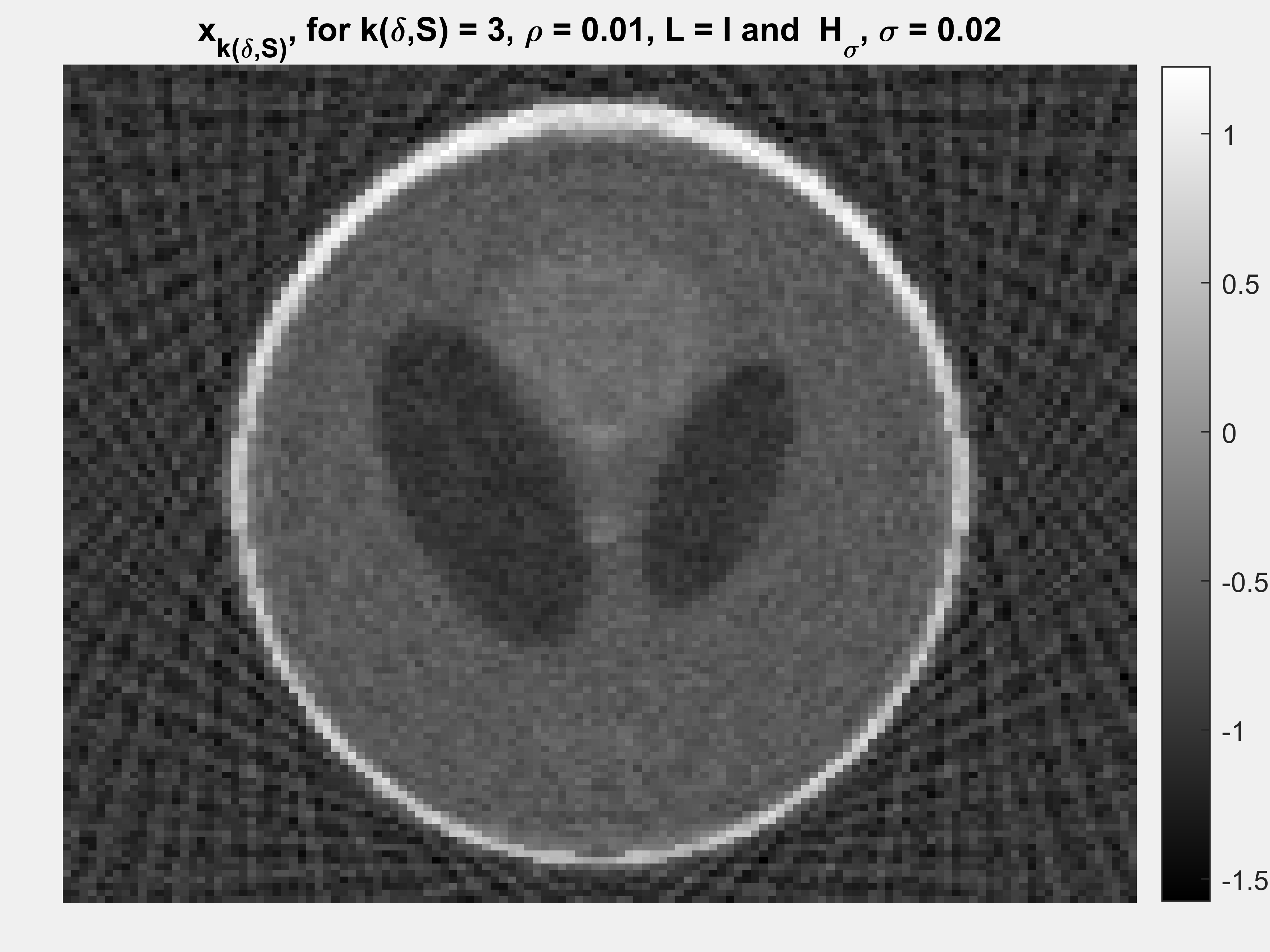}
        \caption{$x_{k(\delta)}^\delta$ for $H_{{\sigma}}$ and $\rho = 0.01$}
    \end{subfigure}    
    \begin{subfigure}{0.495\textwidth}
        \includegraphics[width=\textwidth]{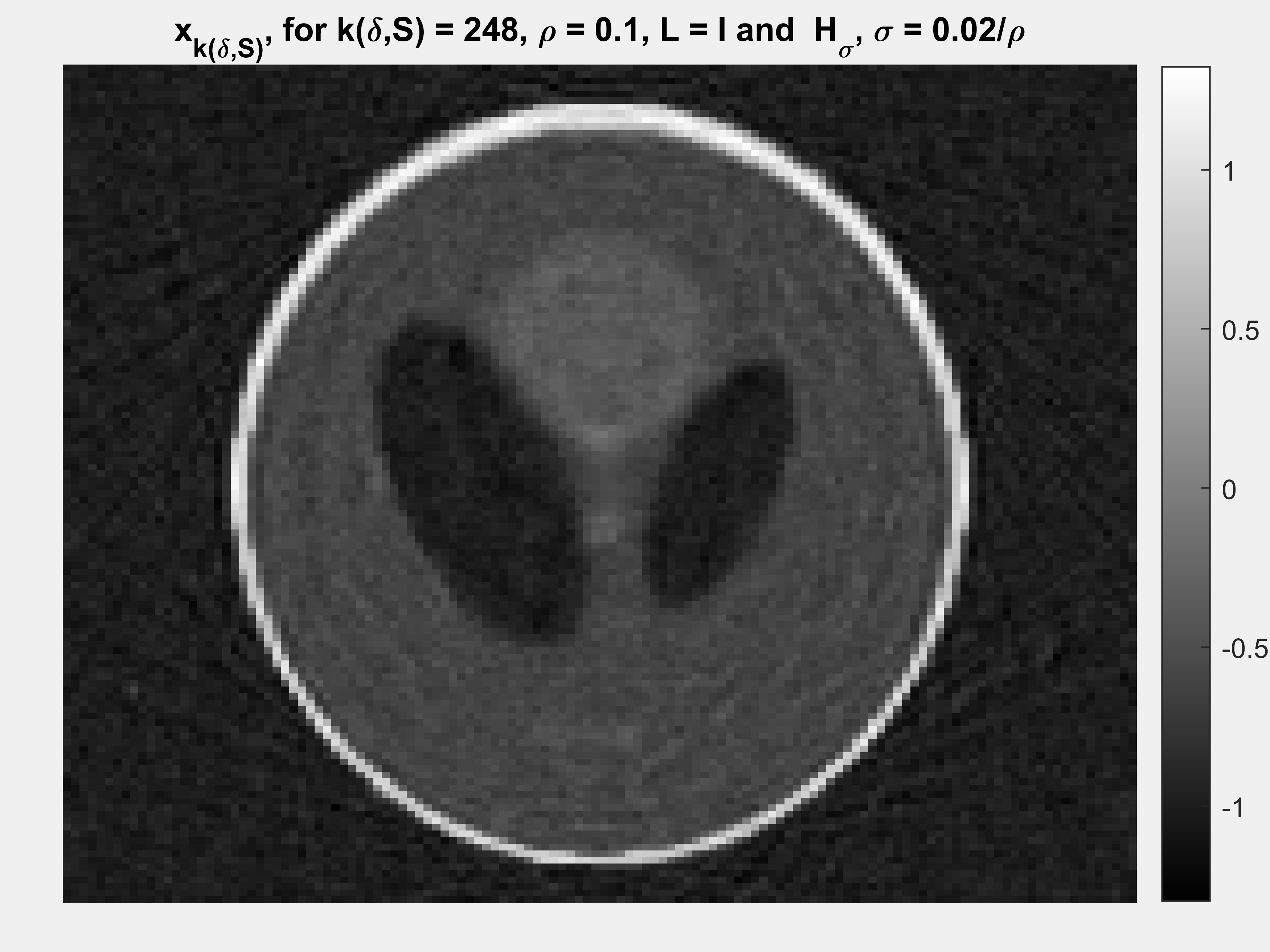}
        \caption{$x_{k(\delta)}^\delta$ for $H_{\hat{\sigma}}$ and $\rho = 0.1$}
        \label{x_CV_ADMM_2 rho update}
    \end{subfigure}
    \begin{subfigure}{0.495\textwidth}
        \includegraphics[width=\textwidth]{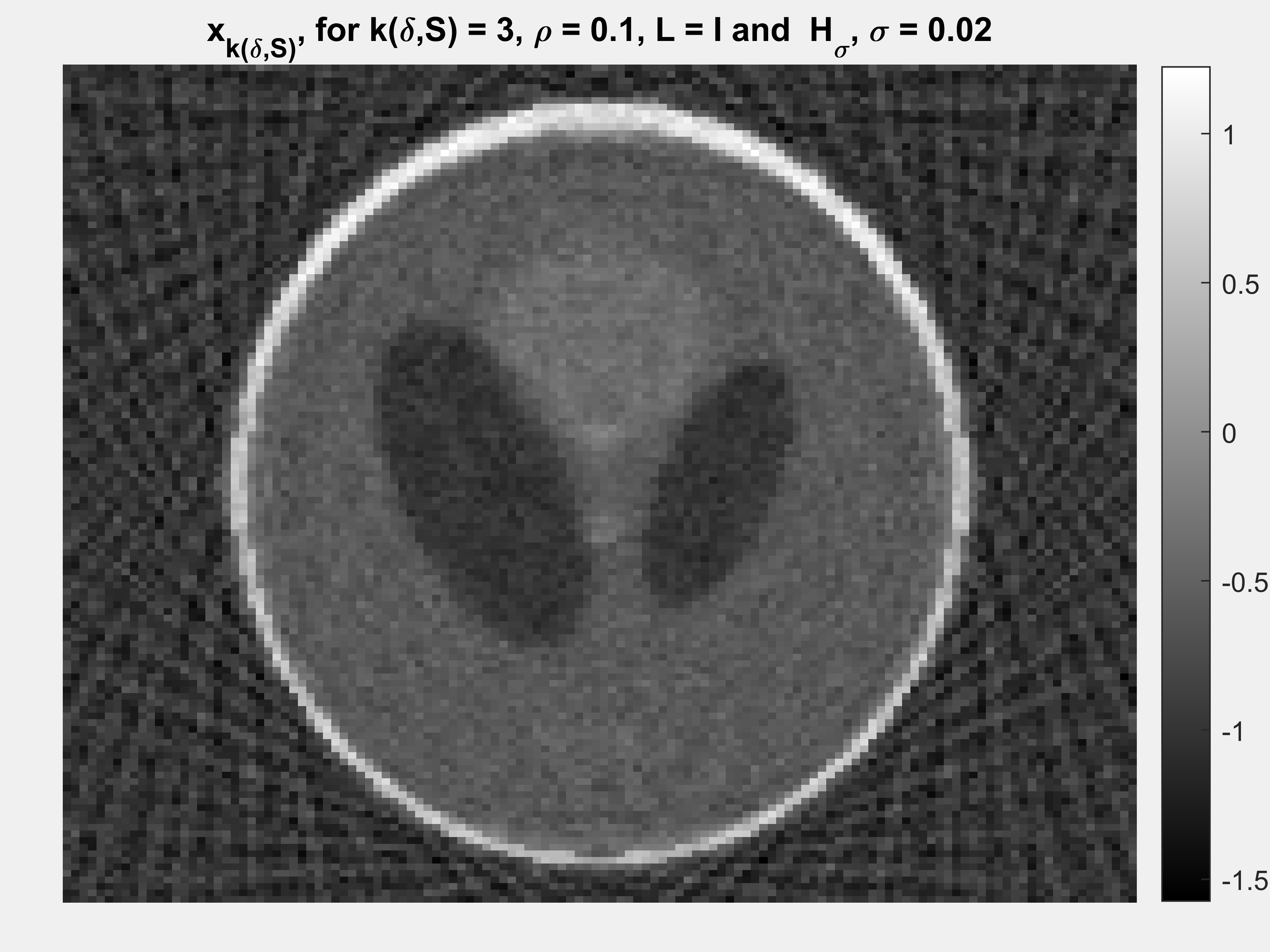}
        \caption{$x_{k(\delta)}^\delta$ for $H_{\sigma}$ and $\rho = 0.1$}
    \end{subfigure}
    \begin{subfigure}{0.495\textwidth}
        \includegraphics[width=\textwidth]{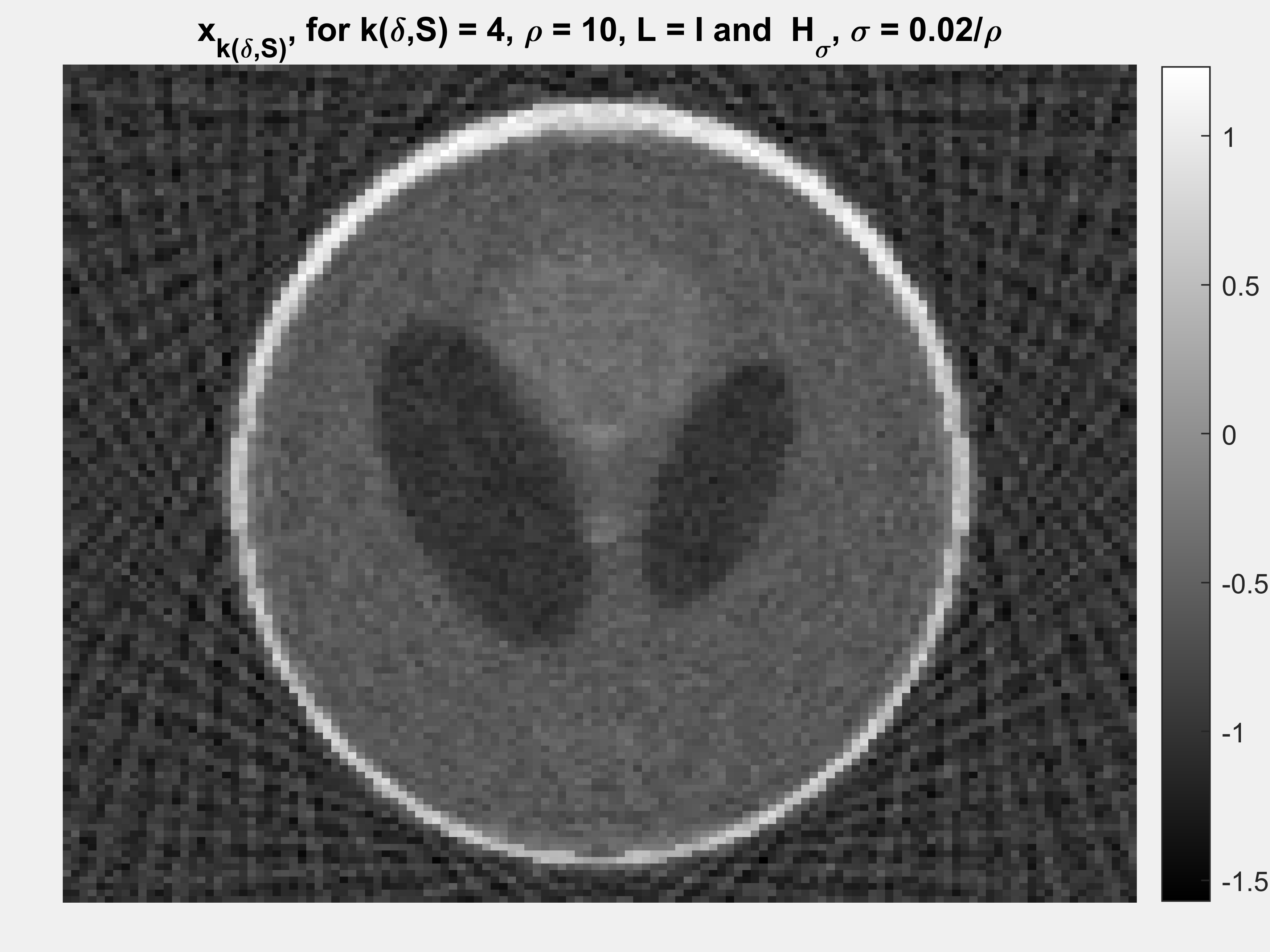}
        \caption{$x_{k(\delta)}^\delta$ for $H_{\hat{\sigma}}$ and $\rho = 10$}
        \label{x_CV_ADMM_3 rho update}
    \end{subfigure}
    \begin{subfigure}{0.495\textwidth}
        \includegraphics[width=\textwidth]{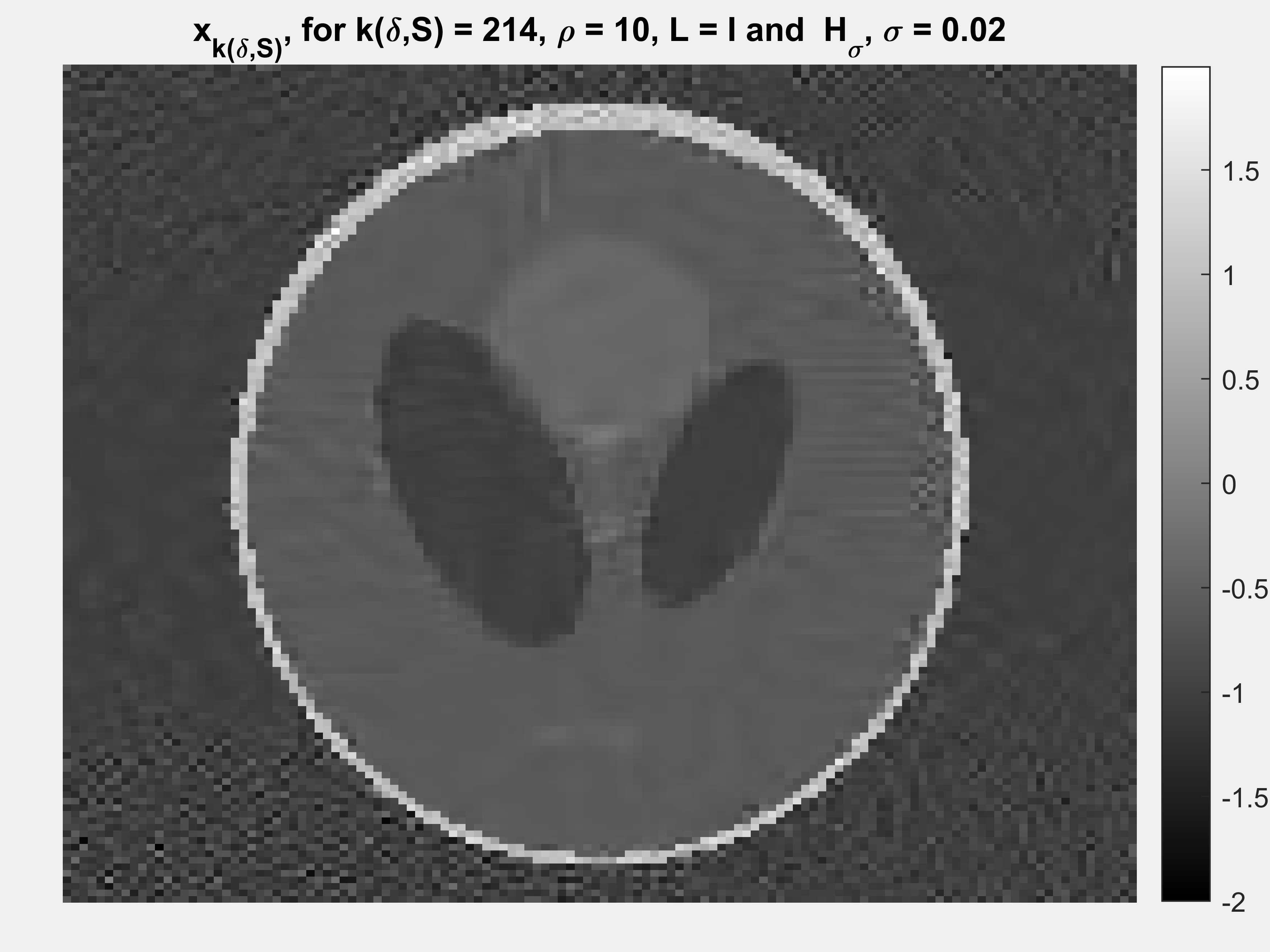}
        \caption{$x_{k(\delta)}^\delta$ for $H_{\sigma}$ and $\rho = 10$}  
    \end{subfigure} 
    \begin{subfigure}{0.495\textwidth}
        \includegraphics[width=\textwidth]{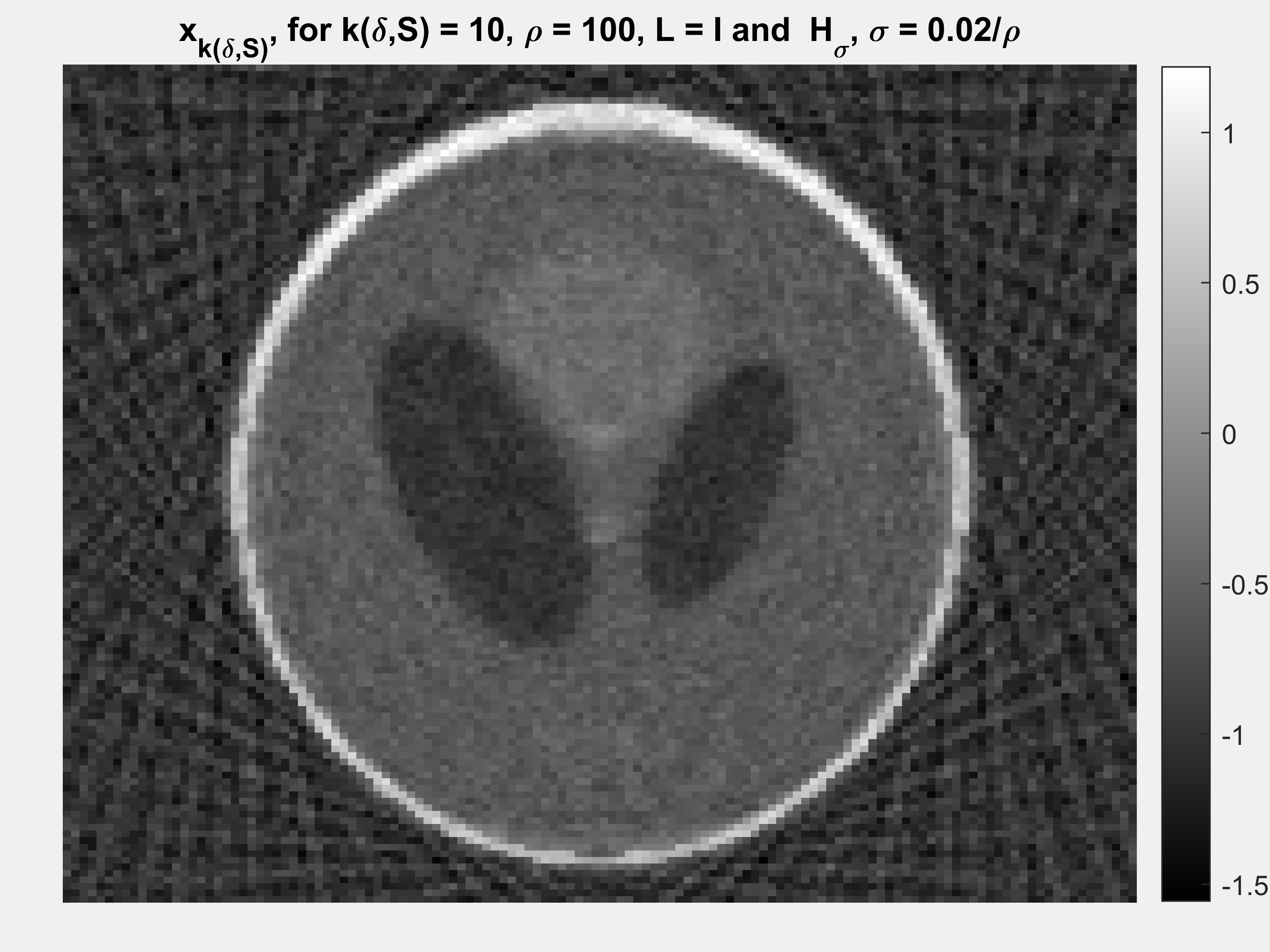}
        \caption{$x_{k(\delta)}^\delta$ for $H_{\hat{\sigma}}$ and $\rho = 100$}
        \label{Fig ADMM-PnPrho=100UpdL=I}
    \end{subfigure}
    \begin{subfigure}{0.495\textwidth}
        \includegraphics[width=\textwidth]{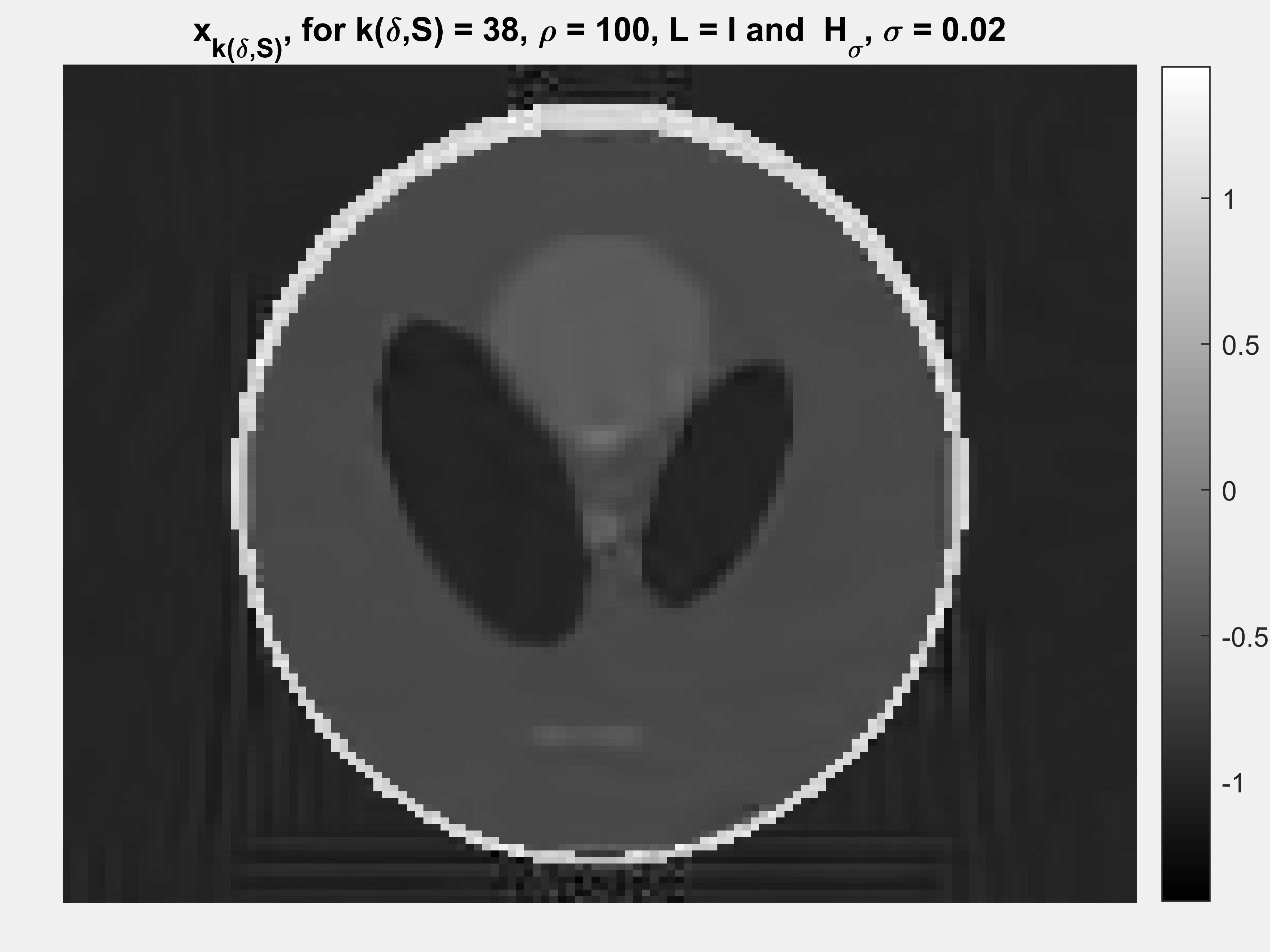}
        \caption{$x_{k(\delta)}^\delta$ for $H_{\sigma}$ and $\rho = 100$}
        \label{Fig ADMM-PnPrho=100NoUpdL=I}
    \end{subfigure}      
    \caption{$H_{\hat{\sigma}=\sigma/\rho}$ vs. $H_\sigma$ for ADMM-PnP, see Example \ref{Ex. ADMM_1}.} 
    \label{Figure ADMM-PnP 1}
\end{figure}

\begin{figure}[h!]
	\centering  
    \begin{subfigure}{0.495\textwidth}
        \includegraphics[width=\textwidth]{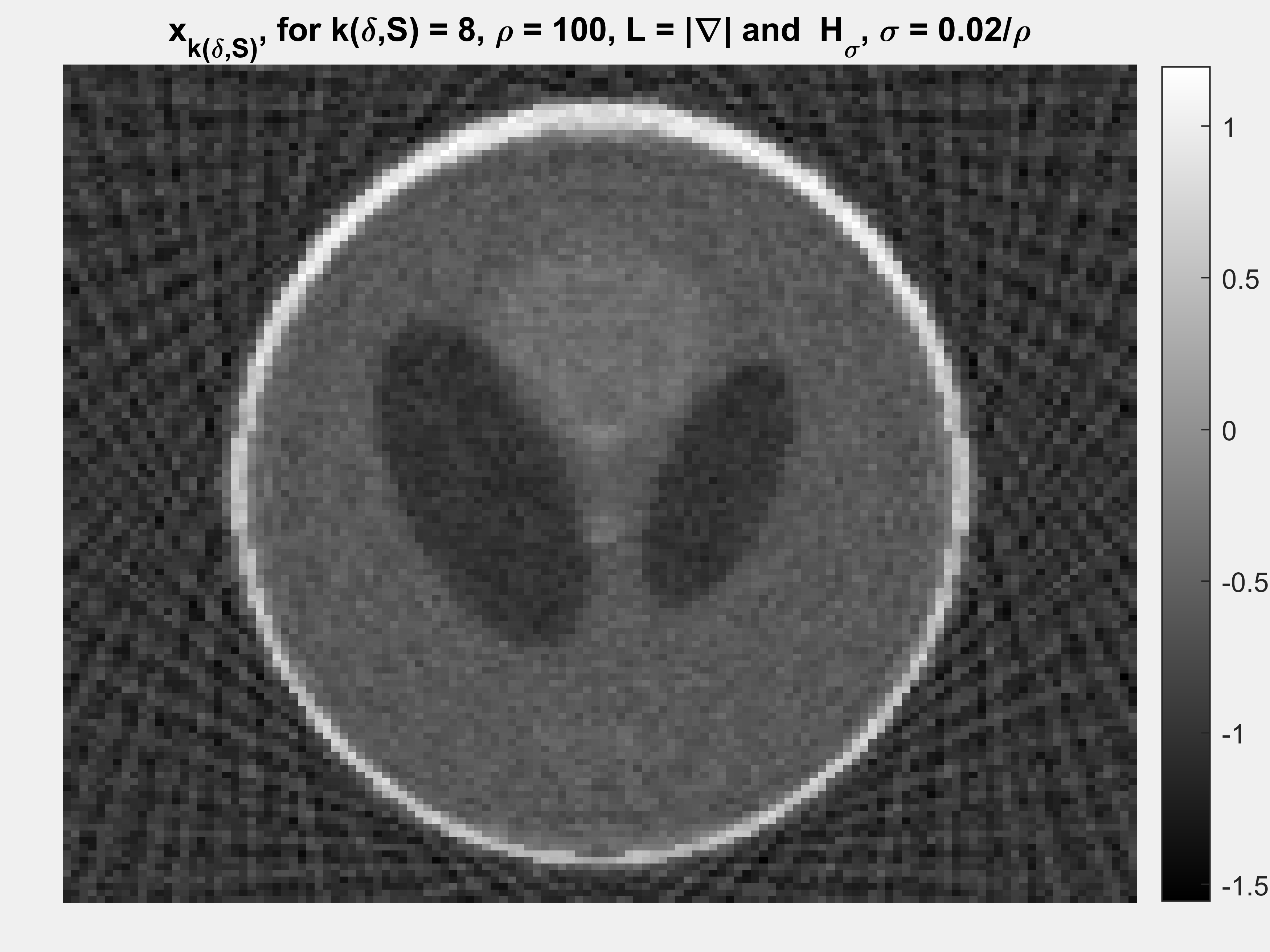}
        \caption{$x_{k(\delta,\mcal{S})}^\delta$, for $H_{\hat{\sigma}}$ and $L = |\nabla |$ for $\rho = 100$}
    \end{subfigure}
    \begin{subfigure}{0.495\textwidth}
        \includegraphics[width=\textwidth]{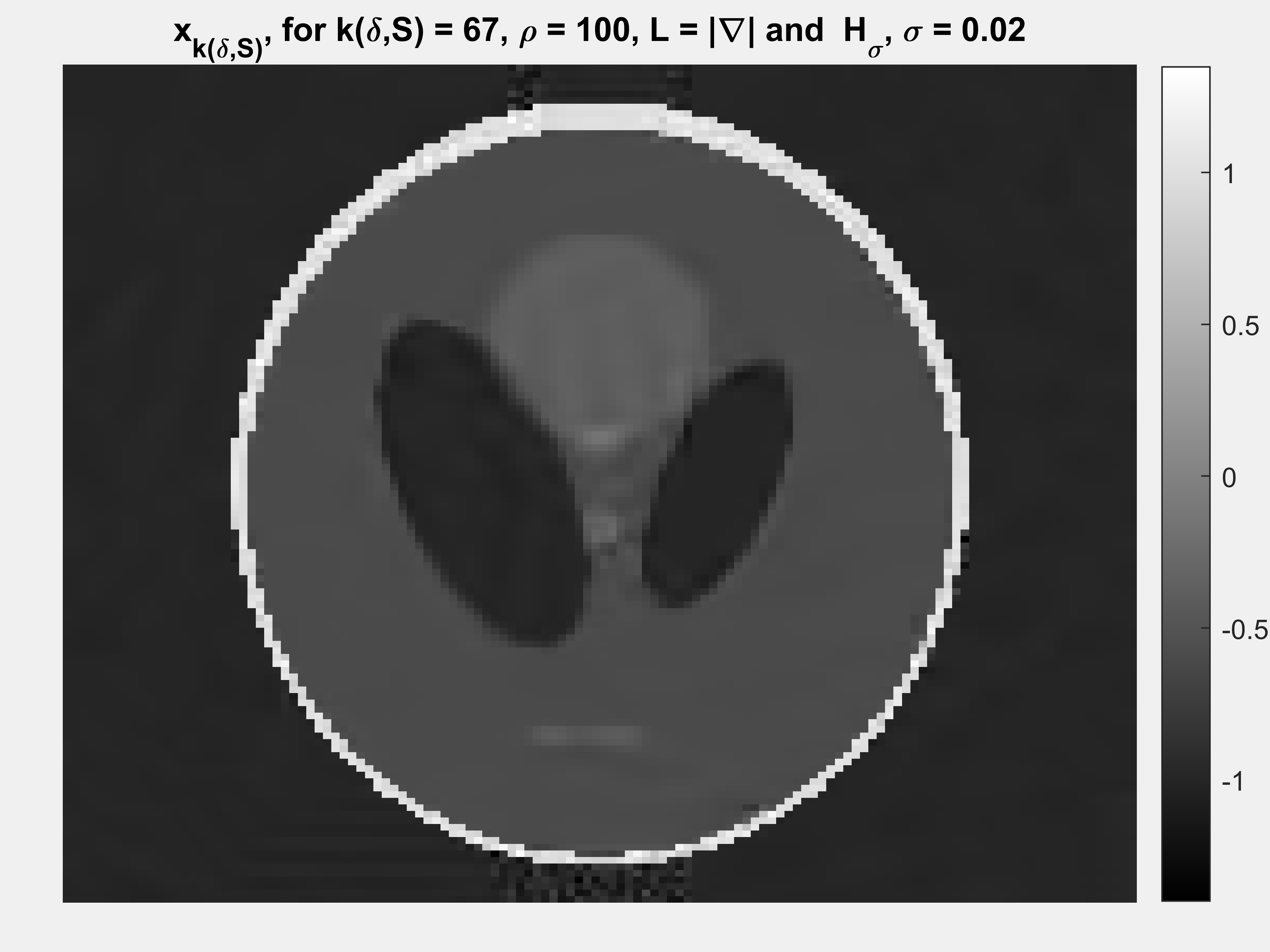}
        \caption{$x_{k(\delta,\mcal{S})}^\delta$, for $H_{{\sigma}}$ and $L = |\nabla |$ for $\rho = 100$}
    \end{subfigure}      
    \caption{$x_{k(\delta,\mcal{S})}^\delta$ with precondition matrix $L = |\nabla |$, see Example \ref{Ex. ADMM_2}.} 
    \label{Figure ADMM-PnP 2}    
\end{figure}

All the experiments are computed in MATLAB, where we consider the selection criterion ($\mcal{S}_0$) as the cross-validation criterion, for some leave out set, and the denoiser $H_{\sigma}$ as the BM3D denoiser, where the denoising parameter $\sigma >0$ corresponds to the standard deviation of the noise and, the MATLAB code for the BM3D denoiser is obtained from \href{http://www.cs.tut.fi/~foi/GCF-BM3D/}{http://www.cs.tut.fi/~foi/GCF-BM3D/}, which is based on \cite{Ymir_Azzari_Foi_1, Ymir_Azzari_Foi_2}. Here, we kept all the attributes of the code in their original (default) settings, except the standard deviation, which depends on $\sigma$. Also, the BM3D denoiser code provided at the above link works for (grayscale or color) images with intensities in the range [0,1]. Hence, for our problems when the image intensities fall outside this range, we transform the original denoiser with a simple rescaling trick, i.e., $\hat{H}_{\sigma}(x) = S^{-1}(H_{\sigma}(S(x)),x_-,x_+)$, where for any $x$ we have $x_- = min(x)$, $x_+ = max(x)$, $S(x) = \frac{x - x_-}{x_+ - x_-}$ and $S^{-1}(y,x_-,x_+) = (x_+ - x_-)y + x_-$. From now on, we denote $H_{\sigma}(x) = \hat{H}_{\sigma}(x)$, unless otherwise stated. We consider the $\ell_2$-discrepancy term, $\mcal{D}(x) = ||Ax - b_\delta||_2^2$, where the noise level ($\delta$) in the data is such that the relative error $\left( \frac{||b - b_\delta||_2}{||b||_2}\%\right)$ is around 1\%, which corresponds to a 40 dB SNR (signal to noise ratio).   
Moreover, in a typical FBS-PnP or ADMM-PnP algorithm, the denoiser is denoted as $H_{\hat{\sigma}}$, where $\hat{\sigma} = \tau \sigma$ ($\tau$ being the step-size) or $\hat{\sigma} = \frac{\sigma}{\rho}$ ($\rho$ being the Lagrangian parameter), respectively. However, to analyze the behavior of the iterative process with respect to the original denoiser, we instead fix the denoiser $H_\sigma$ in the iterative process, i.e., $H_{\hat{\sigma}} = H_{{\sigma}}$, unless otherwise stated. As for the evaluation metrics, we compare the PSNR (peak signal-to-noise ratio) values of the recovered solutions using the MATLAB inbuilt function $psnr(x_k^\delta,\hat{x},\mbox{max}(\hat{x}))$, the SSIM (structure similarity index measure) values using MATLAB function $ssim(S(x_k^\delta),S(\hat{x}))$, where $S$ is the scaling operator (defined above), since the MATLAB $ssim()$ function is well-defined for image intensities in [0,1], and the relative MSE (mean square error) computed as $\frac{||x_k^\delta - \hat{x}||_2}{||\hat{x}||_2}$. The numerical values are shown in Tables \ref{Table Weak deno.}, \ref{Table Strong deno.}, \ref{Table ADMM_1} and \ref{Table ADMM_2}, where $`\mcal{D}$-err.' stands for the discrepancy error $\left(\frac{||Ax_k^\delta - b_{D,\delta}||_2}{||b_{D,\delta}||_2}\right)$ and $`\mcal{S}$-err.' stands for the cross-validation error $\left(\frac{||Ax_k^\delta - b_{S,\delta}||_2}{||b_{S,\delta}||_2}\right)$, where $b_{S,\delta} \subset b_\delta$ is the left-out set and $b_{D,\delta}= b_\delta\backslash b_{S,\delta}$, and the recoveries are shown in Figures \ref{Figure Ex. Weak deno.}, \ref{Figure Ex. Strong deno.}, \ref{Figure ADMM-PnP 1} and \ref{Figure ADMM-PnP 2}.

In the following examples, the matrix equation \eqref{Inv. prob.} corresponds to the discretization of a radon transformation, which is associated with the X-ray computed tomography (CT) reconstruction problem, where we generate the matrix $A \in \mbb{R}^{m \times n}$, $\hat{x} \in \mbb{R}^n$ and $b \in \mbb{R}^m$ from the MATLAB codes presented in \cite{Hansen_IRtools}. The dimension $n$ corresponds to the size of a $N \times N$ image, i.e., $n = N^2$, and the dimension $m$ is related to the number of rays per projection angle and the number of projection angles, i.e., $m = M_1 \times M_2$, where $M_1$ implies the number of rays/angle and $M_2$ implies the number of angles.

\begin{example}\label{Ex. weak deno.}\textbf{[Semi-convergence + Weaker denoiser + Boosting]}\\
Here, we show the semi-convergence nature in the recovery errors shown by PnP-algorithms when an inappropriate denoiser (in this case a weak denoiser) is used during the minimization process and, how can one retrieve the ``best" possible solution in this case. We also present how to boost the denoising strength of a weak denoiser $H_{\sigma}$, without altering the $\sigma$ value, to produce better results. We start with a standard ($128\times 128$) Shepp-Logan phantom (true image $\hat{x} \in \mbb{R}^{16384}$ and $\hat{x}_i \in [0,1]$), which is then re-scaled between [-1,1]. The purpose of rescaling is that, real phantoms are not necessarily restricted to the non-negativity constraint. Moreover, the constraint ($x_i \geq 0$) minimization corresponding to the standard Shepp-Logan phantom is relatively robust to noisy data, and since the point of this example is to show the instabilities occurring in PnP-algorithms, we re-scaled the original phantom. The matrix $A\in \mbb{R}^{8064\times 16384}$ is generated using the $PRtomo()$ code from \cite{Hansen_IRtools}, corresponding to a `fancurved' CT problem  with only 45 view angles (which are evenly spread over $360^o$). The noiseless data is generated by $b := A\hat{x} \in \mbb{R}^{181\times 45 = 8064}$, which is then contaminated by Gaussian noise ($\epsilon_\delta$) to produce noisy data $b_\delta$ such that the relative error is around 1\%. Now, we leave out 1\% of the noisy data $b_\delta$ for the cross-validation error (the selection criterion $\mcal{S}_0$), and start the \textbf{Fast FBS-PnP} algorithm with the modified BM3D denoiser $H_{\sigma}$, with $\sigma = 0.0005$, for all the iterations. We perform the minimization process for two values of the step-size, i.e., (1) $\tau_1 = 2\times 10^{-4}$ and (2) $\tau_2 = 10^{-5}$, over 1000 iterations each, to check the convergence or semi-convergence of the recovery errors.

The results are shown in Table \ref{Table Weak deno.} and the figures in Figure \ref{Figure Ex. Weak deno.}. From Figure \ref{MSE_12}, one can see that the relative errors, when $\tau = 2\times 10^{-4}$, follow a semi-convergence trail, and hence, if the final iterate ($x_N^\delta$) is considered as the recovered solution, then it's much worse than the CV-solution ($x_{k(\delta,\mcal{S})}^\delta$), where $k(\delta,\mcal{S}) = 55$ in this case. Also, note that, by decreasing the step-size ($\tau=10^{-5}$), we are able to recover a much better (smoother) estimate, for the reasons explained in \S \ref{Sec. weaker deno.}. However, it leads to a much slower descent rate and the recovery errors are still following a convergence path at the final iteration, i.e., the recovery will improve upon further iterations. The reason for such a slower recovery rate is the usage of a weak denoiser, for which, one has to use a smaller step-size (to augment the denoising).

\end{example}

\begin{figure}[h!]
	\centering
    \begin{subfigure}{0.495\textwidth}
        \includegraphics[width=\textwidth]{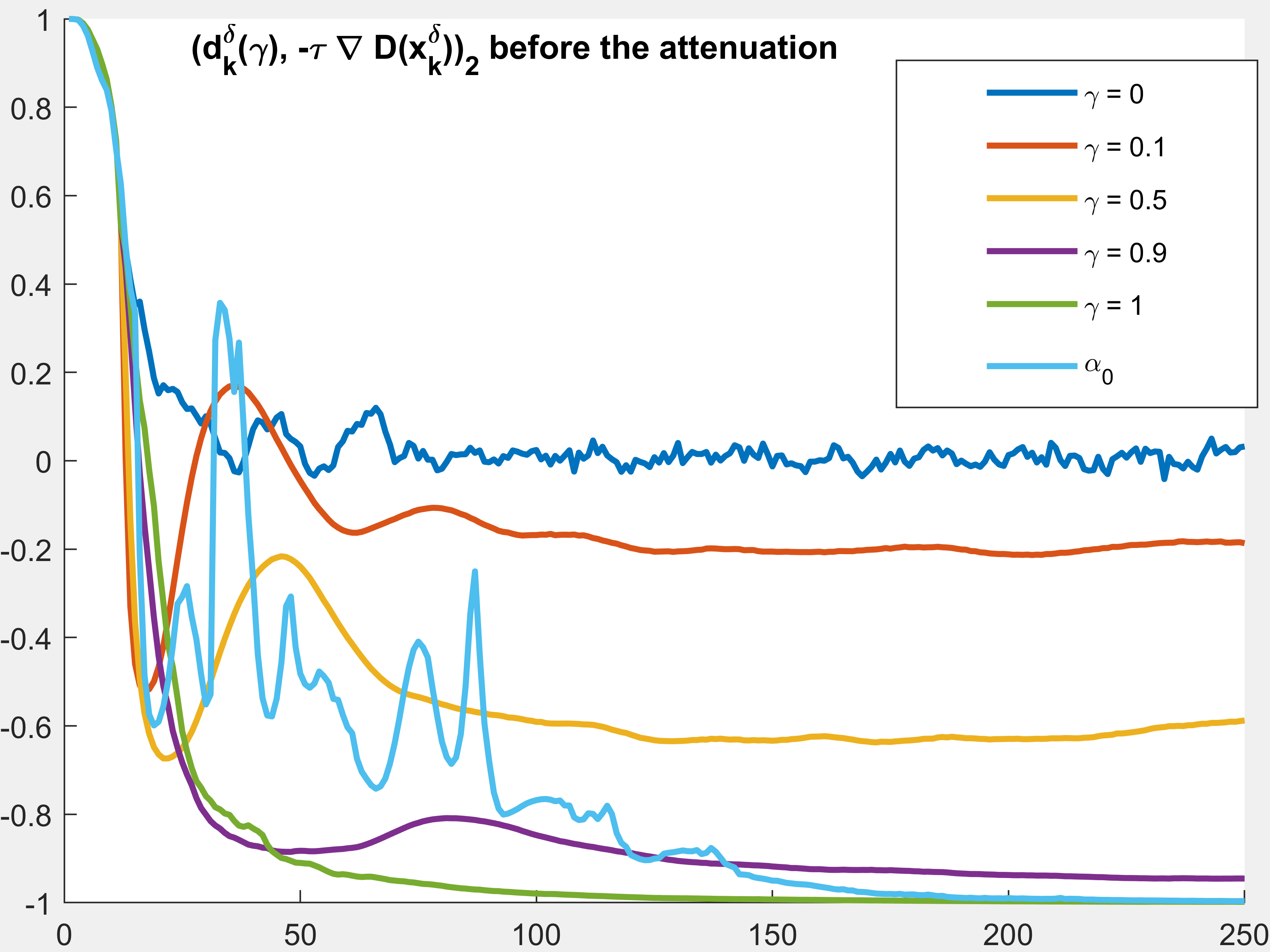}
        \caption{$\lprod{d_k^\delta(\alpha(\gamma))}{-\tau \nabla \mcal{D}(x_k^\delta)}$ before attenuating the denoiser's strength}
        \label{Fig l2prod before}
    \end{subfigure}
    \begin{subfigure}{0.495\textwidth}
        \includegraphics[width=\textwidth]{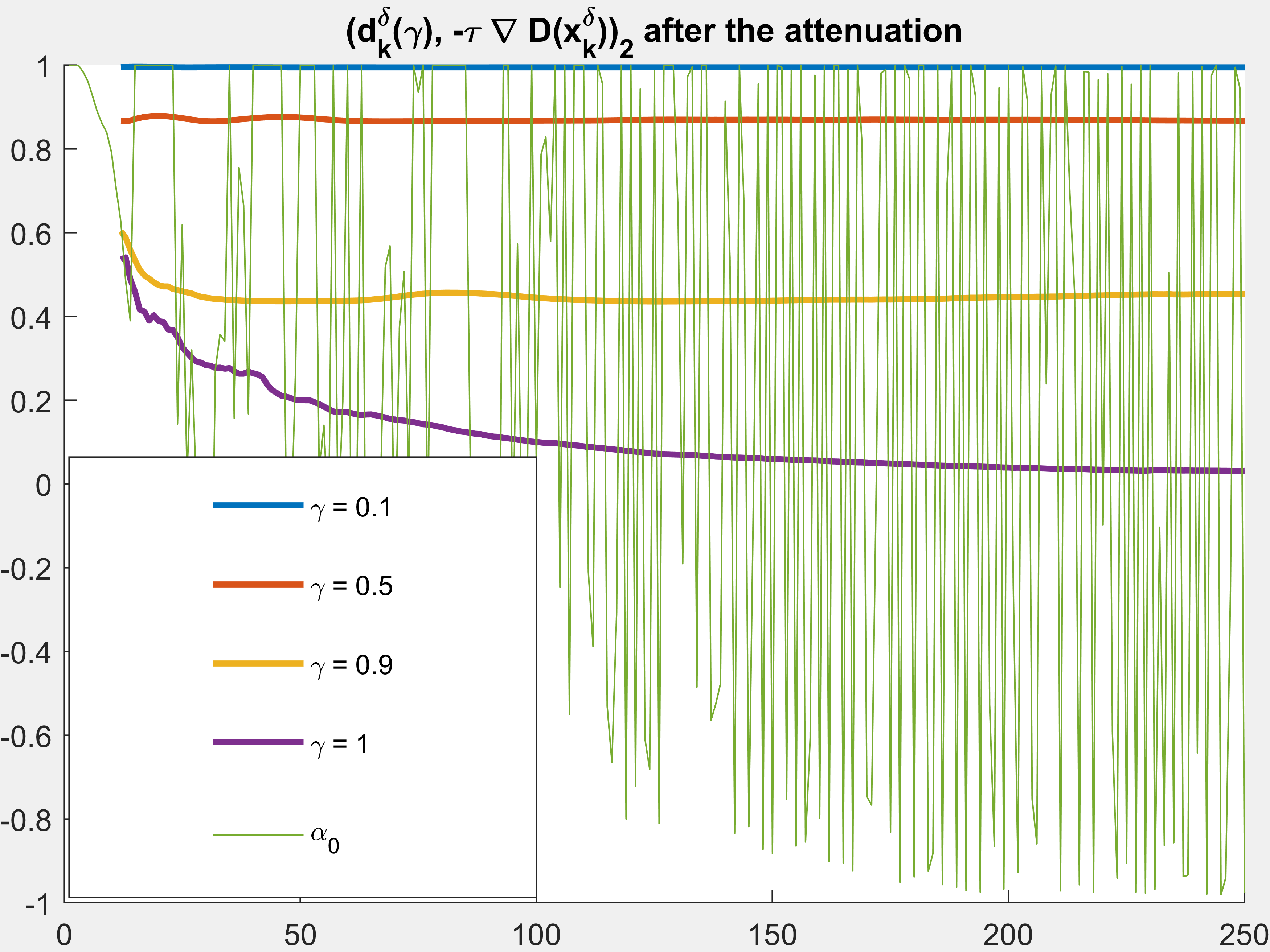}
        \caption{$\lprod{d_k^\delta(\alpha(\gamma))}{-\tau \nabla \mcal{D}(x_k^\delta)}$ after attenuating the denoiser's strength}
        \label{Fig l2prod after}
    \end{subfigure}    
    \begin{subfigure}{0.495\textwidth}
        \includegraphics[width=\textwidth]{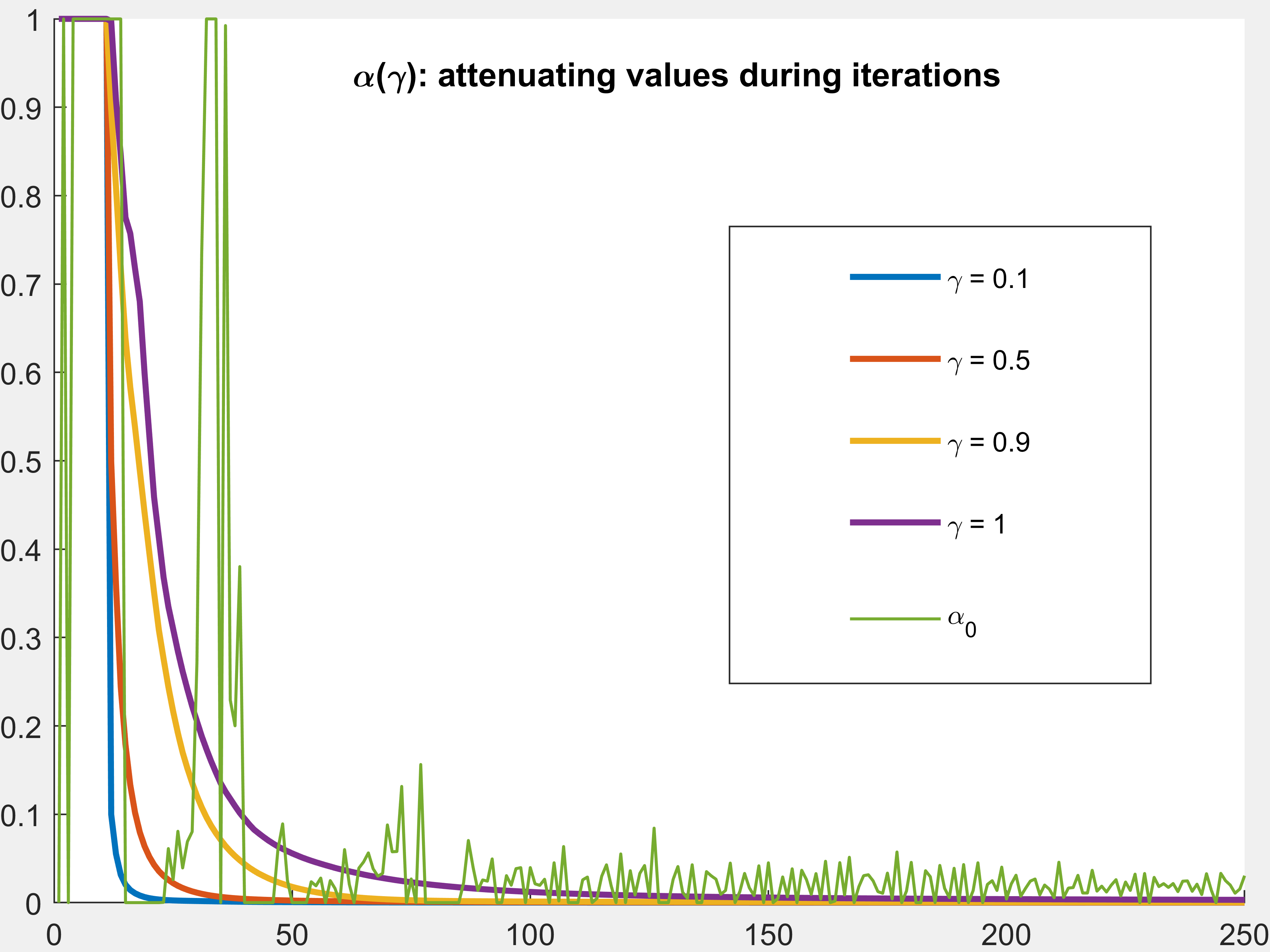}
        \caption{$\alpha(\gamma)$ (attenuating) values over iterations}
        \label{Fig alphagamma}
    \end{subfigure}
    \begin{subfigure}{0.495\textwidth}
        \includegraphics[width=\textwidth]{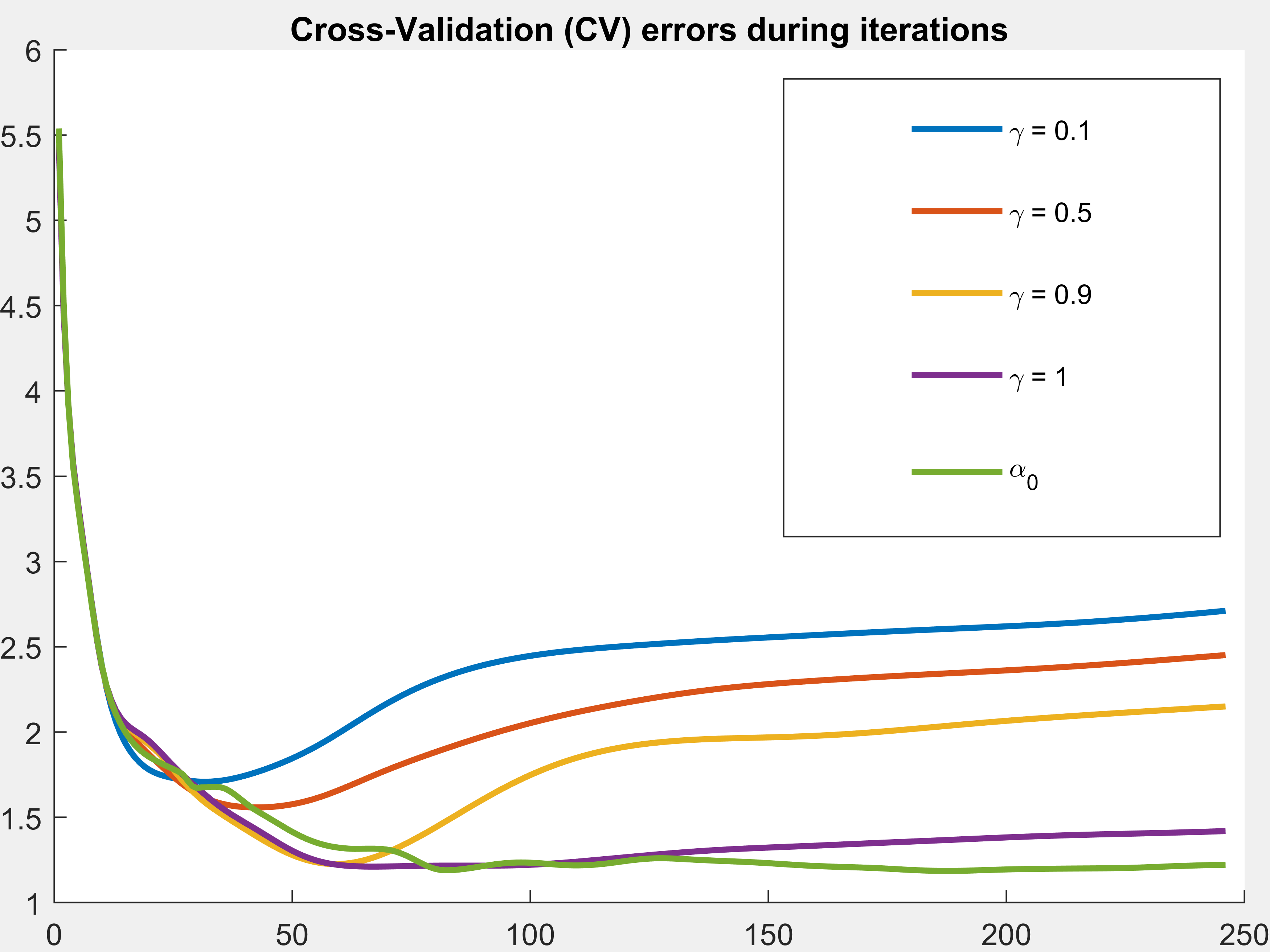}
        \caption{CV-errors for different $\gamma$ values and $\alpha_0$}
        \label{Fig CVerrors}
    \end{subfigure}      
    \caption{ Attenuation of denoising strength, see Example \ref{Ex. Strong deno.}.} 
    \label{Figure l2prod before and after}    
\end{figure}

\begin{example}\label{Ex. Strong deno.}\textbf{[Stronger denoiser + Attenuation + Family $\mcal{I}_3$ or 
$\mcal{I}_5$]}\\
Now, we repeat the experiment performed in Example \ref{Ex. weak deno.} but, using a stronger denoiser $H_{\sigma}$ instead, i.e., a larger value of $\sigma$. Here we set $\sigma = 0.02$ and $\lambda = 2\times 10^{-4}$, and perform the \textbf{Fast FBS-PnP} algorithm for 250 iterations, since $H_\sigma$ is a stronger denoiser. Then we repeat the process by attenuating the denoising strength of $H_\sigma$ through the process described in \S \ref{Sec. stronger deno.}, i.e., using \eqref{Strong deno. alpha1}, \eqref{Strong deno. alpha2} and a combination of them. Note that, while using \eqref{Strong deno. alpha1}, the descent direction is a function of $\gamma$, i.e., greater values of $\gamma$ implies greater denoising and vice-verse, as can be seen in Figures \ref{Figure Ex. Strong deno.}, in addition to, $d_k^\delta(\alpha(\gamma))$ satisfying \eqref{LS descend prop.}. Figure \ref{Fig l2prod before} shows the graph of $\lprod{d_k^\delta}{-\tau\nabla_x \mcal{D}(x_k^\delta)}$ vs. $k$, and Figure \ref{Fig l2prod after} shows the graph of $\lprod{d_k^\delta(\alpha(\gamma))}{-\tau\nabla_x \mcal{D}(x_k^\delta)}$ vs. $k$, for different values of $\gamma$. One can see that 
$d_k^\delta$ does not satisfy \eqref{LS descend prop.} for all values of $k$, where as, 
$d_k^\delta(\alpha(\gamma))$ does satisfy \eqref{LS descend prop.} for all values of $k$ and $\gamma$, i.e., without the attenuation the structured directions $d_k^\delta$ are not descent directions (due to the strong denoising), but after the attenuation, the directions $d_k^\delta(\alpha(\gamma))$ are descent directions. Where as, the same is not true for $d_k^\delta(\alpha_0)$ generated from \eqref{Strong deno. alpha2}, which can be seen in Figure \ref{Fig l2prod after}, as it does not enforce the condition \eqref{LS descend prop.} on $d_k^\delta$. However, one can see from Figure \ref{Fig l2prod before}, which shows the graph of $\lprod{d_k^\delta}{-\tau \nabla_x \mcal{D}(x_k^\delta)}$ vs. k, that the structure imposing directions are completely opposite to the data-consistency directions for larger values of $k$, 
as the inner-product is close to -1, suggesting that iterative process is close to an equilibrium state (i.e., data-consistency and data-denoising steps are negating each other), which is also visible from Figure \ref{Fig l2prod after}, which shows the graph of $\lprod{d_k^\delta(\alpha_0)}{-\tau \nabla_x \mcal{D}(x_k^\delta)}$ vs. k (highly oscillating between 1 and -1), and Figure \ref{MSEcurves_StrongDeno}, the MSE curve. In addition, note that, from Figure \ref{Fig alphagamma}, the $\alpha(\gamma)$ (attenuating) values for the $d_k^\delta(\alpha(\gamma))$ steadily decreases to zero, since the discrepancy values of the iterates $\mcal{D}(x_k^\delta)$ (and hence, the gradient values) decreases to zero; where as, the $\alpha_0$ values does decreases over the iterations but doesn't tend to zero, as it is not associated with the condition \eqref{LS descend prop.}, in fact, the oscillations of the (small) $\alpha_0$ values help the iterative process to attain an equilibrium state; in contrast, the decreasing $\alpha(\gamma)$ values lead to the semi-convergent nature, where the sharpness of the semi-convergence nature depends to the $\gamma$ values (smaller $\gamma \Rightarrow$ sharper semi-convergence). This phenomena is also reflected in the $\mcal{S}$-error curves, Figure \ref{Fig CVerrors}, where for the directions $d_k^\delta(\alpha(\gamma))$ the CV-errors are semi-convergent, but for $d_k^\delta(\alpha_0)$, it is not.

Table \ref{Table Strong deno.} shows the error metrics of the recoveries for different descent directions and different stopping iterations. Figure \ref{Figure Ex. Strong deno.} shows the different recoveries, corresponding to the $k(\delta,\mcal{S})$ iterations, the graph of the PSNR values and the MSE curve. 
As explained in \S \ref{Sec. FBS-PnP}, when using $d_k^\delta(\alpha(\gamma))$, the iterative process falls in the regularization family $\mcal{I}_3$, i.e., semi-convergence of the recovery errors, as can be seen from Figures \ref{MSEcurves_StrongDeno}, for different values of $\alpha(\gamma)$. Where as, when using $d_k^\delta(\alpha_0)$, the recovery process falls in the regularization family $\mcal{I}_5$, i.e., may not exhibit semi-convergence of the recovery errors, as can be seen from Figures \ref{MSEcurves_StrongDeno}. Therefore, it suggests that, when the recovery process belongs to the family $\mcal{I}_3$, i.e., all the $d_k^\delta$s are descent directions, then the minimization process can be terminated at a much earlier instance and a recovered solution $x_{k(\delta)}^\delta$ can correspond to the one satisfying the selection criterion the best, where as, if recovery process falls in family $\mcal{I}_5$, then one can wait for longer period of iterations and can recover better solutions, as the fluctuating/oscillating iterates $z_k^\delta$s can produce better approximations over time. 
\end{example}

\begin{figure}[h!]
	\centering
    \begin{subfigure}{0.495\textwidth}
        \includegraphics[width=\textwidth]{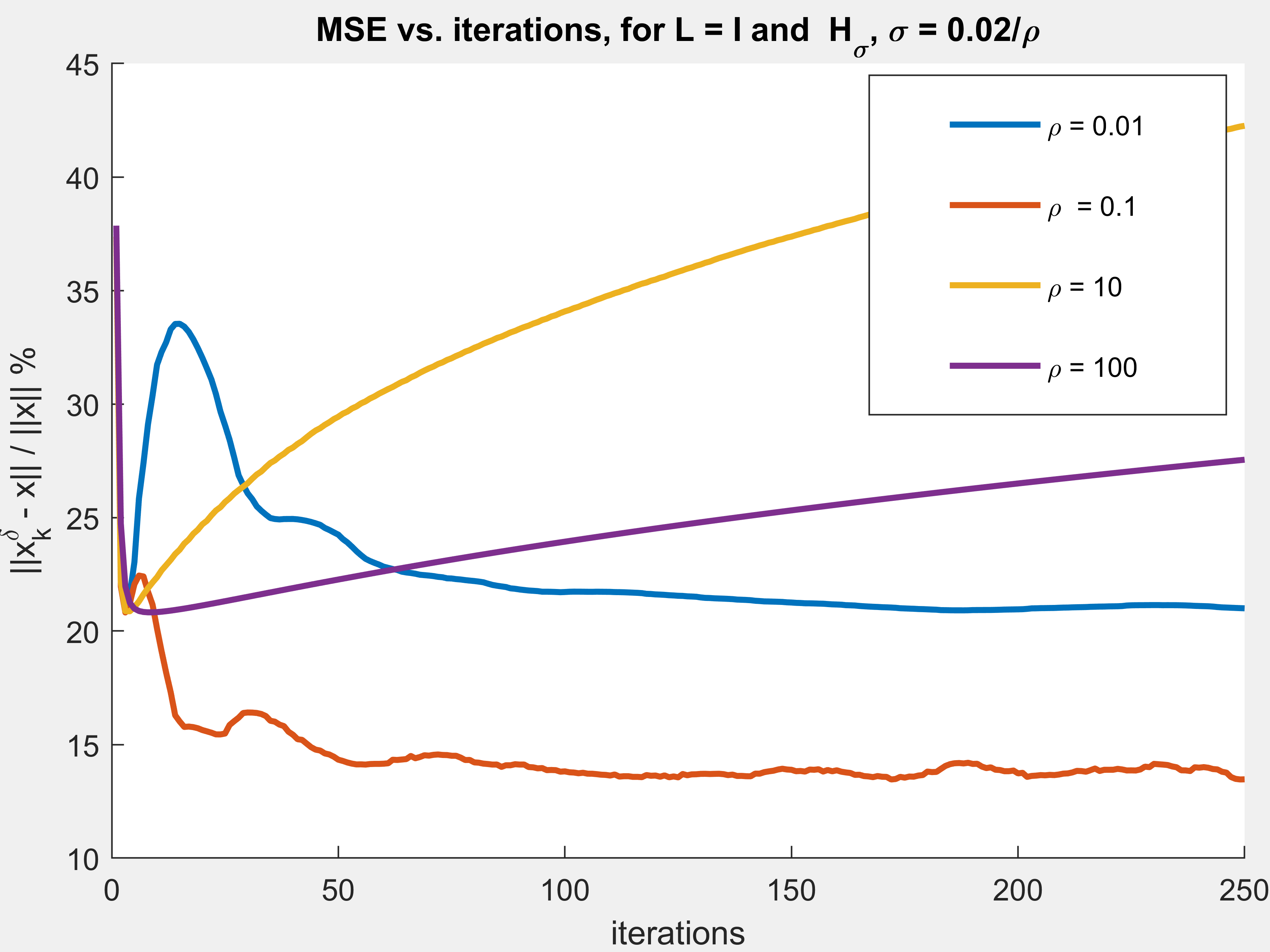}
        \caption{MSE curves for $H_{\hat{\sigma}}$, where $\hat{\sigma} = \sigma/\rho$}
    \end{subfigure}
    \begin{subfigure}{0.495\textwidth}
        \includegraphics[width=\textwidth]{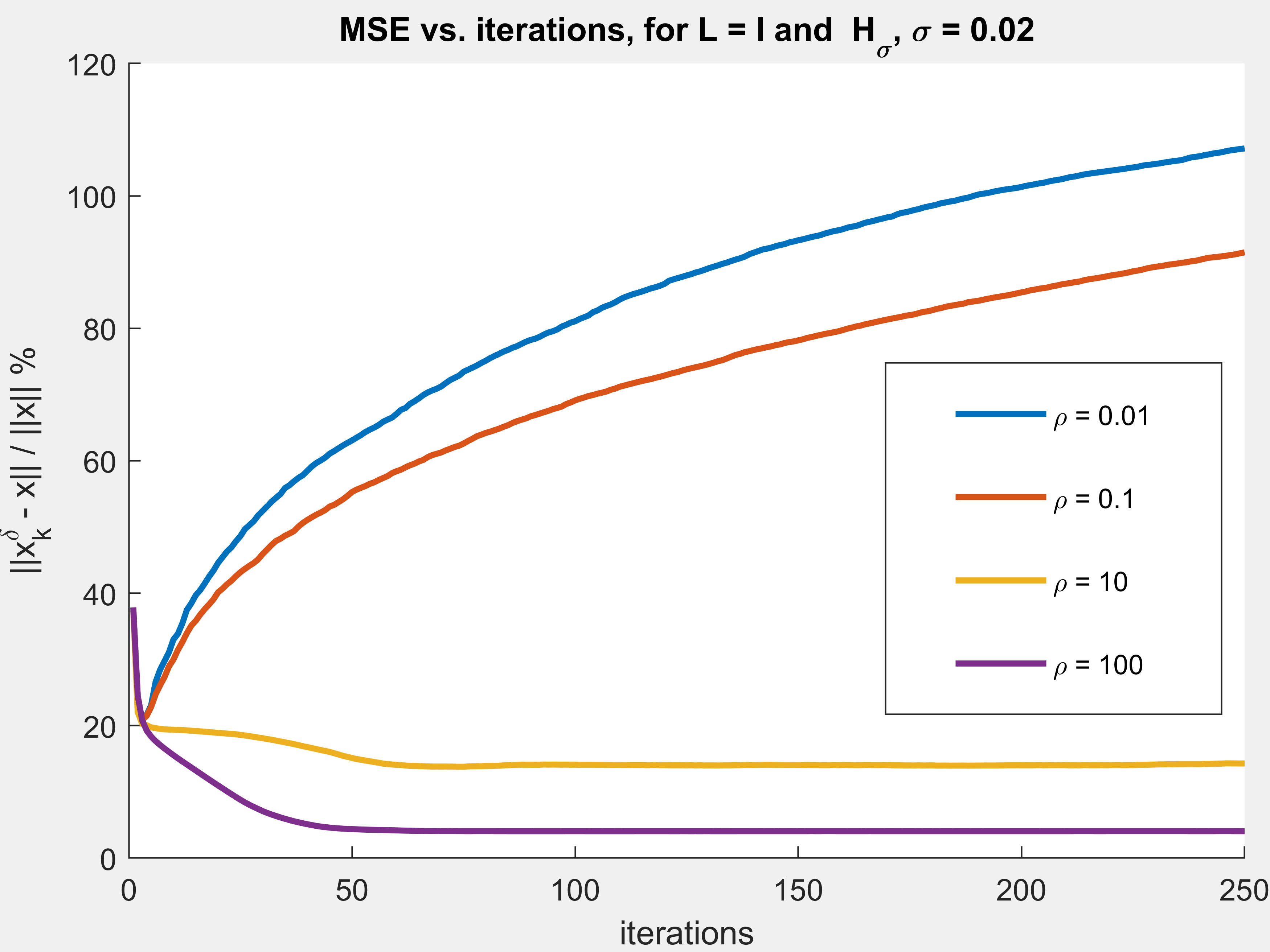}
        \caption{MSE curves for $H_{{\sigma}}$}
    \end{subfigure} 
    \caption{Recovery errors for $H_{\hat{\sigma}}$ vs. $H_{\sigma}$, see Example \ref{Ex. ADMM_1}.} 
    \label{Figure ADMM-PnP MSEcurves 1}	
\end{figure}

\begin{example}\label{Ex. ADMM_1}\textbf{[ADMM - small vs. large $\rho$ + with vs. without $\sigma$ update]}\\
In this example we show the importance of the Lagrangian parameter ($\rho$) and the influence of the update or no-update of the denoiser's denoising strength ($\sigma$-value), during the iterative process, i.e., (1st) we perform the traditional ADMM-PnP algorithm, where $H_\sigma$ is updated to $H_{\hat{\sigma}}$, where $\hat{\sigma} = \frac{\sigma}{\rho}$, during the iterative process, and (2nd) we perform the iterative process without updating $H_{\sigma}$ during the descent process. In addition, for each of the two experiments, we repeat the processes for different values of the Lagrangian parameter ($\rho$), to observe its effect on the recovered solution. Again, for the classical Tikhonov-type regularization, where the regularization function $\mcal{R}(x)$ in \eqref{Gen. Tik. fun.} is proper, closed and convex, the (1st) method is guaranteed to converge to the minimizer of \eqref{Gen. Tik. fun.}, irrespective of the $\rho$-value, although, the minimizer might not be the best estimate of the true solution $\hat{x}$ during the minimization process, for a fixed $\lambda$ value. Note that, as described in \S \ref{Sec. ADMM-PnP}, for smaller $\rho$-value the first iterate $x_1^\delta$ is close to the noisy LS-solution ($x^\dagger_\delta$), associated with the noisy data $b_\delta$. Hence, the denoiser $H_\sigma$ needs to be {particularly good} so that, after updating to $H_{\hat{\sigma}}$, it should be able to clean the large corruptions generated by the noisy data. Where as, for large $\rho$, $x_1^\delta$ is close to $x_0^\delta$, and thus, will be smooth for smoother $x_0^\delta$. Hence, one may clean the artifacts (corruptions) steadily via a {moderately strong} denoiser. However, if $H_\sigma$ is updated to $H_{\hat{\sigma}=\frac{\sigma}{\rho}}$, then the denoising ability of $H_\sigma$ is also reduced, by a factor corresponding to $\frac{1}{\rho}$, and thus, a moderately strong denoiser becomes a weaker denoiser, and may not be able to clean the corruptions. Therefore, for the fixed denoiser $H_\sigma$, (2nd) method suggests to incorporate the denoiser $H_\sigma$ {without updating} it, to preserve the denoising strength of a moderately strong denoiser and provide effective denoising. Moreover, since the convergence of ADMM-PnP algorithm for any general denoiser $H_\sigma$ is not known, it does not hurt to implement it in either way; after-all, both these processes can be thrown in either $\mcal{I}_5$ or $\mcal{I}_6$ family and the better process is the one for which the recovered solution satisfies the selection criterion ($\mcal{S}$) the best.

Here we keep the experimental settings similar to Example \ref{Ex. weak deno.}, but the iterative process follows the ADMM-PnP algorithm, using the stronger denoiser $H_\sigma$, with $\sigma=0.02$. For each of the above two strategies, we repeat the experiments for four values of the Lagrangian parameter, $\rho = 0.01, 0.1, 10, 100$, for over 250 iterations. Also, the minimizer of each $x$-update, i.e., $x_k^\delta$ as defined in \eqref{ADMM-PnP x_k}, is only estimated by a certain number of conjugate-gradient (CG) steps and not through the matrix inversion. Here, we considered 100 CG-steps, initiated from $z_{k-1}^\delta - u_{k-1}^\delta$, to approximate the k-th step x-minimizer, i.e., to generate the $x_k^\delta$ iterate. Hence, the 100 CG-steps, in combination with the 250 ADMM-steps, leads to a total of 250,000 iterations, although, there are only 250 denoising steps (which can be computationally expensive). Also, note that, the denoiser $H_\sigma$ for $\sigma = 0.02$ acted as a strong denoiser in Example \ref{Ex. Strong deno.}, where as, here it behaves like a weak denoiser, except when $\rho$ is large and $\sigma$ is not updated to $\hat{\sigma}$. Interestingly, even when $\sigma$ is updated to $\hat{\sigma} = \frac{\sigma}{\rho} >> \sigma$ (for smaller $\rho$), the denoiser $H_{\hat{\sigma}}$ is still not adequately strong to clean the corruptions in the iterates as effectively as $H_\sigma$ with larger $\rho$ value, which validates the aforementioned explanations. Table \ref{Table ADMM_1} shows the error metrics of the recoveries corresponding to the minimum CV-error, Figure \ref{Figure ADMM-PnP 1} shows the recoveries and, Figure \ref{Figure ADMM-PnP MSEcurves 1} shows the MSE curves of the recoveries. 
\end{example}

\begin{example}\label{Ex. ADMM_2}
In this example we compare the recoveries obtained via the precondition matrix $L = |\nabla|$, for the first iterate ($x_1^\delta$) in \eqref{ADMM-PnP x_1}, with the recoveries using $L=I$. We repeat the experiment in Example \ref{Ex. ADMM_1}, for $\rho=100$, but using $L=|\nabla|$ to calculate $x_1^\delta$. The CV-solution ($x_{k(\delta,\mcal{S})}^\delta$), with and without $\sigma$ update, for $L=I$ is shown in Figures \ref{Fig ADMM-PnPrho=100UpdL=I} and \ref{Fig ADMM-PnPrho=100NoUpdL=I}, and the CV-solutions for $L = |\nabla|$ is shown in Figure \ref{Figure ADMM-PnP 2}. Note that, the recoveries corresponding to $H_{\hat{\sigma}}$ is bad in both the cases and the recovery using $L=|\nabla|$ is slightly better than $L=I$, see Table \ref{Table ADMM_2}. 
\end{example}

%

\section{Conclusion and Future Research}
In this paper we tried to explain the PnP-algorithms from a different angle. We defined certain families of regularized solutions and showed that the solution of a PnP-algorithm will fall in one of them. We started with an extension of the Landweber iterations to generate structure imposing descent directions, i.e., any direction $d_k^\delta$ satisfying \eqref{LS descend prop.}. This lead to the formulation of the regularization family $\mcal{I}_3$, which contains all the semi-iterative processes. This is further generalized to the family $\mcal{I}_5$, where the directions $d_k^\delta$ satisfies \eqref{LS descend prop.} for certain number of iterations and then the restriction is relaxed, i.e., once the iterates $x_k^\delta$ have approximated the data $b_\delta$ to a certain extent ($\mcal{D}(x_k^\delta) \leq \epsilon_1(\delta)$), we don't want $d_k^\delta$ to be a descent direction anymore, as then, it will lead to semi-convergence in the recovery errors. The advantage of such a formulation is that, in addition to avoiding the semi-convergence of the recovery errors, one doesn't even have to worry about the convergence of the iterates $x_k^\delta$ anymore, i.e., the iterative process doesn't necessarily need to be associated with a penalized/constrained Variational-minimization problem, and thus, the iterates $x_k^\delta$ do not need to converge to a minimizer. Hence, it compensates the short-comings of both the classical regularization methods, i.e., the semi-iterative methods and the Tikhonov-type methods. In other words, any iterate $x_k^\delta$ during the minimization process is an ``approximate estimate" to the solution of the inverse problem ($x^\dagger$), and the appropriateness of the estimate is determined via a selection criterion $\mcal{S}_0$ for the solution of the inverse problem \eqref{Inv. prob.}. Therefore, one can terminate the iterative process at an early instance, if the error in the selection criterion $\mcal{S}_0$ starts getting worse, and hence, can reduce the computational time significantly; or can wait for longer period of iterations, if $\mcal{S}(x_k^\delta)$ values is decreasing, assuming the selection criterion $\mcal{S}_0$ is appropriate for the problem. 
Note that, one does not necessarily need a Variational formulation to regularize the solution of an inverse problem. A ``regularized solution $x^\delta$" is a solution which ``approximates $\hat{x}$", depending on the noise level $\delta$, and ``gets better" as ``noise vanishes", i.e., $\delta \rightarrow 0$ should imply $||x^\delta - \hat{x}|| \rightarrow 0$, and ``$x^\delta$ is well-defined" if there is a ``unique process" of obtaining it, for example, either via a Variational formulation or via an iterative formulation. Hence, this leads to the formulation of our largest family of regularized solution $\mcal{I}_5$, as defined in \eqref{LS. family I_6}, i.e., a solution state $x^\delta(.,t)$ corresponding to an initial value differential equation problem.

Furthermore, we showed how to improve the recoveries when dealing with an inappropriate denoiser $H_\sigma$, without altering the denoiser, i.e., when the denoiser is too strong to oversmooth the recovered solution, then one can attenuate the denoising strength via a relaxing parameter, as described in \S \ref{Sec. stronger deno.}, or, when the denoiser is very weak to remove the corruptions properly, then one can reduce the step-size to boost the denoising strength, as described in \S \ref{Sec. weaker deno.}. We also showed the importance of the Lagrangian parameter ($\rho$) in the ADMM-PnP algorithm and the influence of updating or not updating the denoiser $H_\sigma$ to $H_{\hat{\sigma}}$, where $\hat{\sigma}=\frac{\sigma}{\rho}$, during the iterative process, as well, as the significance of the preconditioned matrix $L$ (for the first iterate $x_1^\delta$) in \eqref{ADMM-PnP x_1}. The validations of these improving techniques can be seen in Examples \ref{Ex. weak deno.}, \ref{Ex. Strong deno.}, \ref{Ex. ADMM_1} and \ref{Ex. ADMM_2}.

In an upcoming paper, we extend this idea to the setting of training examples or data, i.e., for a neural network denoiser or a denoising iterative scheme (such as unrolled neural network). We believe that by incorporating this idea of regularization, i.e., monitoring the recovery process via an appropriate selection criterion, one can understand and control the instabilities arising in deep learning based reconstruction algorithms, as shown in \cite{Antun_Renna_Poon_Adcock_Hansen}. 

\begin{table}[h!]
    \centering
    \begin{tabular}{|p{1.9cm}||p{1.1cm}|p{1cm}|p{1cm}|p{1cm}|p{1cm}|p{1cm}|p{1.9cm}|}
    \hline
    \multicolumn{8}{|c|}{Fast FBS-PnP + a weak denoiser $H_\sigma$ ($\sigma = 0.0005$), for $\tau = 2\times 10^{-4}$ vs. $10^{-5}$}\\
    \hline
    step-size & iter.(k) & MSE & $\mcal{D}$-err.  & $\mcal{S}$-err. & PSNR & SSIM & Min.MSE \\
    \hline
    $\tau=2\times 10^{-4}$ & 1000(N) & 0.3432 & 0.3264 & 0.0290 & 10.51 & 0.2926 & 0.0208 (50)\\
    \hline
	 & 55 ($\mcal{S}$) & 0.2088 & 0.4891  & 0.0248 & 14.83 & 0.2540  &\\ 
    \hline  
    $\tau=10^{-5}$ & 1000(N) & 0.0899 & 0.0075 & 0.0134 & 22.15 & 0.3932 & 0.0899(1000) \\ 
    \hline  
	 & 1000($\mcal{S}$) & \textbf{0.0899} & 0.0075 & \textbf{0.0134} & \textbf{22.15} & \textbf{0.3932} &  \\ 
    \hline               
    \end{tabular}
    \caption{Comparing solution $z_N^\delta$ vs. $z_{k(\delta,\mcal{S})}^\delta$ for Example \ref{Ex. weak deno.}.}
    \label{Table Weak deno.}
\end{table}

\begin{table}[h!]
    \centering
    \begin{tabular}{|p{2.25cm}||p{1cm}|p{1cm}|p{1cm}|p{1cm}|p{1cm}|p{1cm}|p{1.8cm}|}
    \hline
    \multicolumn{8}{|c|}{Fast FBS-PnP + a strong denoiser $H_\sigma$ ($\sigma = 0.01$) +  $d_k^\delta$ vs. $d_k^\delta(\alpha(\gamma))$ vs. $d_k^\delta(\alpha_0)$}\\
    \hline
    directions & iter.(k) & MSE & $\mcal{D}$-err.  & $\mcal{S}$-err. & PSNR & SSIM & Min.MSE \\
    \hline
    \hline
    $d_k^\delta$ & 250(N) & 0.1491 & 0.0174 & 0.0182 & 17.75 & 0.3991 & 0.1477(36)\\
    \hline
		& 36($\mcal{S}$) & 0.1477 & 0.0171 & 0.0176 & 17.83 & 0.3750 & \\
    \hline
    \hline
	$d_k^\delta(\alpha(\gamma=0.1))$ & 250(N) & 0.2312 & 0.0033 & 0.0271 & 13.94 & 0.2284 & 0.1660(52)\\
    \hline         
	 & 36($\mcal{S}$) & 0.1682 & 0.0064 & 0.0171 & 16.71 & 0.2903 & \\
    \hline       
	$d_k^\delta(\alpha(\gamma=0.5))$ & 250(N) & 0.1846 & 0.0034 & 0.0245 & 15.90 & 0.2730 & 0.1181(75)\\
    \hline         
	 & 46($\mcal{S}$) & 0.1242 & 0.0055 & 0.0156 & 19.34 & 0.3395 & \\
    \hline   
	$d_k^\delta(\alpha(\gamma=0.9))$ & 250(N) & 0.1433 & 0.0036 & 0.0215 & 18.10 & 0.3046 & 0.0823(69)\\
    \hline         
	 & 63($\mcal{S}$) & 0.0832 & 0.0057 & 0.0122 & 22.82 & 0.3697 & \\
    \hline  
	$d_k^\delta(\alpha(\gamma=1))$ & 250(N) & 0.0729 & 0.0048 & 0.0142 & 23.97 & 0.3673 & 0.0542(102)\\
    \hline         
	 & 71($\mcal{S}$) & 0.0648 & 0.0073 & 0.0121 & 25.00 & 0.3928 &\\
    \hline       
    \hline    
	$d_k^\delta(\alpha_0)$ & 250(N) & 0.0455 & 0.0066 & 0.0122 & \textbf{28.07} & \textbf{0.4106} & \textbf{0.0440(131)}\\
    \hline         
	 & 192($\mcal{S}$) & 0.0467 & 0.0066 & \textbf{0.0118} & 27.83 & 0.4039 & \\
    \hline     
    \end{tabular}
    \caption{Comparing solution $z_N^\delta$ vs. $z_{k(\delta,\mcal{S})}^\delta$ for Example \ref{Ex. Strong deno.}.}
    \label{Table Strong deno.}
\end{table}

\begin{table}[h!]
    \centering
    \begin{tabular}{|p{2.3cm}||p{1cm}|p{1cm}|p{1cm}|p{1cm}|p{1cm}|p{1cm}|p{1.6cm}|}
    \hline
    \multicolumn{8}{|c|}{ADMM-PnP using $H_\sigma$, for $\sigma = 0.02$, vs. $H_{\hat{\sigma}}$, where $\hat{\sigma} = \sigma/\rho$}\\
    \hline
    denoiser and $\rho$ & iter.(k) & MSE & $\mcal{D}$-err.  & $\mcal{S}$-err. & PSNR & SSIM & Min.MSE \\
    \hline
    \hline
    $H_{\hat{\sigma}}$ \& $\rho = 0.01$ & 250(N) & 0.2100 & 0.0044 & 0.0239 & 14.78 & 0.2467 & 0.2088(3)\\
    \hline                  
	 & 219($\mcal{S}$) & 0.2108 & 0.0043 & 0.0232 & 14.75 & 0.2472 & \\
	 \hline
    $H_{{\sigma}}$ \& $\rho = 0.01$ & 250(N) & 01.07 & 0.0026 & 0.0776 & 0.6266 & 0.1307 & 0.2088(3)\\
    \hline                  
	 & 3($\mcal{S}$) & 0.2088 & 0.0048 & 0.0249 & 14.83 & 0.2533 & \\
    \hline 	 
    \hline   
    $H_{\hat{\sigma}}$ \& $\rho=0.1$ & 250(N) & 0.1347 & 0.0064 & 0.0149 & 18.64 & 0.3289 & 0.1346(172)\\
    \hline                  
	 & 248($\mcal{S}$) & 0.1349 & 0.0064 & 0.0148 & 18.63 & 0.3284 & \\
    \hline  
    $H_{{\sigma}}$ \& $\rho = 0.1$ & 250(N) & 0.9147 & 0.0026 & 0.0671 & 2.00 & 0.2115 & 0.2087(3)\\
    \hline                  
	 & 3($\mcal{S}$) & 0.2087 & 0.0048 & 0.0249 & 14.83 & 0.2537 & \\
    \hline 	 
    \hline               
    $H_{\hat{\sigma}}$ \& $\rho=10$ & 250(N) & 0.4225 & 0.0029 & 0.0373 & 08.71 & 0.1589 & 0.2087(4)\\
    \hline                  
	 & 4($\mcal{S}$) & 0.2087 & 0.0047 & 0.0247 & 14.83 & 0.2533 & \\
    \hline
    $H_{{\sigma}}$ \& $\rho = 10$ & 250(N) & 0.1428 & 0.0044 & 0.0159 & 18.13 & 0.3434 & 0.1379(74)\\
    \hline                  
	 & 214($\mcal{S}$) & 0.1403 & 0.0044 & 0.0153 & 18.28 & 0.3388 & \\
    \hline 	 
    \hline      
    $H_{\hat{\sigma}}$ \& $\rho=100$ & 250(N) & 0.2754 & 0.0032 & 0.0294 & 12.42 & 0.1954 & 0.2082(8)\\
    \hline                  
	 & 10($\mcal{S}$) & 0.2083 & 0.0049 & 0.0246 & 14.85 & 0.2548 & \\
    \hline                 
    $H_{{\sigma}}$ \& $\rho = 100$ & 250(N) & 0.0406 & 0.0087 & 0.0121 & \textbf{29.06} & \textbf{0.4130} & \textbf{0.0405(96)}\\
    \hline                  
	 & 38($\mcal{S}$) & 0.0540 & 0.0087 & \textbf{0.0120} & 26.58 & 0.4074 & \\
    \hline    
    \end{tabular}
    \caption{Comparing solution $z_N^\delta$ vs. $z_{k(\delta,\mcal{S})}^\delta$ for Example \ref{Ex. ADMM_1}.}
    \label{Table ADMM_1}
\end{table}

\begin{table}[h!]
    \centering
    \begin{tabular}{|p{2.3cm}||p{1cm}|p{1cm}|p{1cm}|p{1cm}|p{1cm}|p{1cm}|p{1.6cm}|}
    \hline
    \multicolumn{8}{|c|}{ADMM-PnP, for $\sigma=0.02$ and $\rho=100$, $H_\sigma$ vs. $H_{\hat{\sigma}}$ and $L=I$ vs. $L = |\nabla |$}\\
    \hline
    denoiser and $L$ & iter.(k) & MSE & $\mcal{D}$-err.  & $\mcal{S}$-err. & PSNR & SSIM & Min.MSE \\
    \hline
    \hline
    $H_{\hat{\sigma}}$ \& $L=I$ & 250(N) & 0.2754 & 0.0032 & 0.0294 & 12.42 & 0.1954 & 0.2082(8)\\
    \hline                  
	 & 10($\mcal{S}$) & 0.2083 & 0.0049 & 0.0246 & 14.85 & 0.2548 & \\     
    \hline
    $H_{{\sigma}}$ \& $L=I$ & 250(N) & 0.0406 & 0.0087 & 0.0121 & {29.06} & {0.4130} & {0.0405(96)}\\
    \hline                  
	 & 38($\mcal{S}$) & 0.0540 & 0.0087 & {0.0120} & 26.58 & 0.4074 & \\
    \hline   
    \hline
    $H_{\hat{\sigma}}$ \& $L=|\nabla|$ & 250(N) & 0.2750 & 0.0032 & 0.0291 & 12.44 & 0.1957 & 0.2073(8)\\
    \hline                  
	 & 8($\mcal{S}$) & 0.2073 & 0.0051 & 0.0235 & 14.89 & 0.2572 & \\     
    \hline
    $H_{{\sigma}}$ \& $L=|\nabla|$ & 250(N) & 0.0404 & 0.0087 & 0.0121 & \textbf{29.10} & \textbf{0.4131} & \textbf{0.0401(79)}\\
    \hline                  
	 & 67($\mcal{S}$) & 0.0404 & 0.0087 & \textbf{0.0120} & 29.09 & 0.4074 & \\
    \hline       
    \end{tabular}
    \caption{Comparing solution $z_N^\delta$ vs. $z_{k(\delta,\mcal{S})}^\delta$ for Example \ref{Ex. ADMM_2}.}
    \label{Table ADMM_2}
\end{table}

\bibliography{thesisref} 
\bibliographystyle{ieeetr}

\end{document}